\documentclass[a4paper,11pt]{amsart}
\usepackage{amsmath,amsthm,amssymb}
\usepackage[mathscr]{eucal}
 \usepackage{cite}
\usepackage{upgreek}
\usepackage[bookmarks,bookmarksnumbered]{hyperref}
\usepackage{enumerate}

\numberwithin{equation}{section}

\theoremstyle{plain}
\newtheorem{main}{Theorem}
\newtheorem{mcor}[main]{Corollary}
\newtheorem{mprop}[main]{Proposition}

\newtheorem{theorem}{Theorem}[section]

\newtheorem{lemma}[theorem]{Lemma}
\newtheorem{proposition}[theorem]{Proposition}
\newtheorem{corollary}[theorem]{Corollary}
\theoremstyle{definition}
\newtheorem{definition}[theorem]{Definition}
\newtheorem{example}[theorem]{Example}

\newtheorem{remark}[theorem]{Remark}

\begin{document}

\title[Orbit equivalence rigidity for translation actions]
{Strong ergodicity, property (T), and orbit equivalence rigidity for translation actions}
\author[Adrian Ioana]{Adrian Ioana}
\thanks{The author\ was partially supported by  NSF  Grant DMS \#1161047, NSF Career Grant DMS \#1253402, and a Sloan Foundation Fellowship.}
\address{Mathematics Department; University of California, San Diego, CA 90095-1555 (United States).}
\email{aioana@ucsd.edu}

\begin{abstract} 

We study equivalence relations that arise from translation actions $\Gamma\curvearrowright G$ which are associated to dense embeddings $\Gamma<G$ of countable groups  into second countable locally compact groups. Assuming that $G$ is simply connected and the action $\Gamma\curvearrowright G$ is strongly ergodic, we prove that $\Gamma\curvearrowright G$ is orbit equivalent to another such translation action $\Lambda\curvearrowright H$ if and only if there exists an isomorphism $\delta:G\rightarrow H$ such that $\delta(\Gamma)=\Lambda$. 
If $G$ is moreover a real algebraic group, then we establish analogous rigidity results for the translation actions of $\Gamma$ on homogeneous spaces of the form $G/\Sigma$, where $\Sigma<G$ is either a discrete or an algebraic subgroup. We also prove that if $G$ is simply connected and  the action $\Gamma\curvearrowright G$ has property (T), then any cocycle $w:\Gamma\times G\rightarrow\Lambda$ with values into a countable group $\Lambda$ is cohomologous to a homomorphism $\delta:\Gamma\rightarrow\Lambda$. As a consequence, we deduce that the action $\Gamma\curvearrowright G$ is orbit equivalent superrigid: any free nonsingular action $\Lambda\curvearrowright Y$ which is orbit equivalent to $\Gamma\curvearrowright G$, is necessarily conjugate to an induction of $\Gamma\curvearrowright G$.
\end{abstract}

\maketitle

\section{Introduction and statement of main results}

In the last 15 years there have been many exciting developments in the study of equivalence relations arising from nonsingular actions of countable groups  (see the surveys \cite{Po07,Fu09,Ga10}). 
The goal of the present work is to investigate the equivalence relations that are associated to dense embeddings of countable groups into locally compact groups.
More precisely, consider a locally compact second countable (abbreviated, {\it l.c.s.c.}) group  $G$ endowed with a left Haar measure  $m_G$, and a countable dense subgroup $\Gamma<G$.
 Then the left translation action $\Gamma\curvearrowright G$ is measure preserving and ergodic.  
 We denote by $\mathcal R(\Gamma\curvearrowright G)$ the {\it orbit equivalence relation} on $G$ of belonging to the same $\Gamma$-orbit: $x\sim y$ $\Longleftrightarrow$ $\Gamma x=\Gamma y$.
 This paper is motivated by the following question: to what extent does  $\mathcal R(\Gamma\curvearrowright G)$ remember the original inclusion $\Gamma<G$ it was constructed from? 

In the case $G$ is {\it compact}, this question has been investigated recently  by the author  \cite{Io08, Io13} and by A. Furman  \cite{Fu09}. Firstly, it was shown that if $\Gamma$ has property (T), then the action $\Gamma\curvearrowright (G,m_G)$ is {\it orbit equivalent superrigid}: 
 any probability measure preserving action $\Lambda\curvearrowright (Y,\nu)$ which is orbit equivalent to $\Gamma\curvearrowright (G,m_G)$ must be virtually conjugate to it. This result was obtained in \cite{Io08}, when $G$ is a profinite group,  and in \cite{Fu09}, for general compact groups $G$.   Secondly, assume that $G$ is a profinite group and that $\Gamma\curvearrowright (G,m_G)$ has spectral gap.  Then the action $\Gamma\curvearrowright (G,m_G)$ was very recently shown to satisfy the following {rigidity} statement: any translation action $\Lambda\curvearrowright (H,m_H)$ on a profinite group $H$ that is orbit equivalent to $\Gamma\curvearrowright (G,m_G)$ must be virtually conjugate to it \cite{Io13}.

On the other hand, in the case $G$ is  locally compact but {\it non-compact}, little is known about how properties of the inclusion $\Gamma<G$ reflect on the structure of $\mathcal R(\Gamma\curvearrowright G)$. In fact, even basic questions, that can be easily answered in the compact case,  are extremely challenging in the non-compact case.
This is best illustrated by the question of characterizing when is the equivalence relation $\mathcal R(\Gamma\curvearrowright G)$ amenable (or, equivalently by \cite{CFW81}, hyperfinite).
If $G$ is compact, then $\mathcal R(\Gamma\curvearrowright G)$ is amenable if and only if $\Gamma$ is amenable. Only recently,  by combining the structure theory of locally compact groups \cite{MZ55} with their topological Tits alternative, E. Breuillard and T. Gelander  were able to answer this question for arbitrary locally compact groups $G$. They showed that $\mathcal R(\Gamma\curvearrowright G)$ is amenable  if and only if there exists an open subgroup $G_0<G$ such that $\Gamma\cap G_0$ is amenable \cite{BG04}. Note that in the case $G$ is a connected Lie group, this result was established earlier by R. Zimmer  in  \cite{Zi87} (see also \cite{CG85}, for a proof in the case $G=SL_2(\mathbb R)$, and \cite{BG02}, for an alternative proof of the general case).

In this paper, we investigate $\mathcal R(\Gamma\curvearrowright G)$ for general locally compact groups $G$, and obtain several rigidity results. Roughly speaking, these results provide instances when the inclusion $\Gamma<G$ can be partially, or even entirely, recovered from  $\mathcal R(\Gamma\curvearrowright G)$.
 More precisely, if the action $\Gamma\curvearrowright (G,m_G)$ verifies certain conditions that strengthen non-amenability (strong ergodicity/property (T)), we prove that the equivalence relation  $\mathcal R(\Gamma\curvearrowright G)$ satisfies rigidity/superrigidity statements analogous to the ones obtained in \cite{Io08,Fu09,Io13} in the case when $G$ compact (see Theorems \ref{thmA}, \ref{thmB} and \ref{thmC}). Our method of proof relies on a result for untwisting cocycles $w:\Gamma\times G\rightarrow\Lambda$ with values into countable groups $\Lambda$ (see Theorem \ref{cocrig}). As such,  we are able to more generally study equivalence relations $\mathcal R(\Gamma\curvearrowright G/\Sigma)$ which arise from translation actions on homogeneous spaces $G/\Sigma$, where $\Sigma<G$ is a closed subgroup  (see Theorems \ref{thmE} and \ref{thmF}). 

Before stating our main results in detail, we first review some terminology, starting with the notions of nonsingular actions and orbit equivalence.

\subsection {Nonsingular actions and orbit equivalence}
Let $G$ be a l.c.s.c. group and $X$ be a standard Borel space together with a Borel action $G\curvearrowright X$. A  $\sigma$-finite Borel measure $\mu$ on $X$ is  {\it quasi-invariant} under the $G$-action  if $\mu(A)=0\Longrightarrow\mu(gA)=0$, for every measurable set  $A\subset X$ and  all $g\in G$. The measure $\mu$ is called {\it invariant} if $\mu(gA)=\mu(A)$, for every measurable set $A\subset X$ and all $g\in G$.  If  $\mu$ is quasi-invariant, then the action $G\curvearrowright (X,\mu)$ is called {\it nonsingular}.  If  $\mu$ is invariant, then the action $G\curvearrowright (X,\mu)$ is called {\it measure preserving}. In this case, if $\mu$ is a probability measure, we say that the action is {\it probability measure preserving} (abbreviated, {\it p.m.p.}).
A nonsingular action $G\curvearrowright (X,\mu)$ is called  {\it free} if the stabilizer in $G$ of almost every point is trivial, and {\it ergodic} if any $G$-invariant measurable set $A\subset X$ is either null or conull, i.e.  $\mu(A)(1-\mu(A))=0$.

\begin{example}
Let $G$ be a  l.c.s.c. group and $\Sigma<G$  a closed subgroup. Then there exists a Borel measure $\mu$ on $G/\Sigma$ which is quasi-invariant under the action $G\curvearrowright G/\Sigma$.  Moreover, any two such measures $\mu$, $\mu'$ are equivalent, i.e. they have the same null sets. From now on, we fix a quasi-invariant measure on $G/\Sigma$, which we denote by $m_{G/\Sigma}$.  Whenever $G/\Sigma$ admits a $G$-invariant measure, we choose $m_{G/\Sigma}$ to be $G$-invariant. In particular, $m_G$ denotes a left Haar measure of $G$.
In this paper, we study left translation actions of the form $\Gamma\curvearrowright (G/\Sigma,m_{G/\Sigma})$, where $\Gamma<G$ is a countable subgroup.
If $\Gamma<G$ is dense (as we will typically assume), then this action is ergodic.

\end{example}

\begin{definition}
Two nonsingular actions $\Gamma\curvearrowright (X,\mu)$, $\Lambda\curvearrowright (Y,\nu)$ of two countable groups $\Gamma,\Lambda$ on standard measure spaces $(X,\mu), (Y,\nu)$ are said to be {\it orbit equivalent} (OE) if there is a nonsingular isomorphism $\theta:X\rightarrow Y$ such that $\theta(\Gamma x)=\Lambda\theta(x)$, for almost every $x\in X$. The actions are called
{\it stably orbit equivalent} (SOE)  if there are non-negligible measurable sets $A\subset X$, $B\subset Y$ and a nonsingular isomorphism $\theta:A\rightarrow B$ such that $\theta(\Gamma x\cap A)=\Lambda\theta(x)\cap B$, for almost every $x\in A$.
Finally, the actions are called
 {\it conjugate} if there exist a nonsingular isomorphism $\theta:X\rightarrow Y$ and a group isomorphism $\delta:\Gamma\rightarrow\Lambda$ such that $\theta(gx)=\delta(g)\theta(x)$, for all $g\in\Gamma$ and almost every $x\in X$.
\end{definition}

\subsection{Strong ergodicity and property (T)} In order to formulate our main results, we also need to recall the notion of strong ergodicity and property (T) for nonsingular actions.

\begin{definition} \cite{CW81,Sc81}\label{strerg}
Let $\Gamma$ be a countable group and $\Gamma\curvearrowright (X,\mu)$ be  a nonsingular ergodic  action on a standard probability space $(X,\mu)$. A sequence  $\{A_n\}$ of measurable subsets of $X$ is said to be {\it asymptotically invariant} ({\it a.i.}) if
$\lim\limits_{n\rightarrow\infty}\mu(gA_n\Delta A_n)=0$, for all $g\in\Gamma$. The action $\Gamma\curvearrowright (X,\mu)$ is called {\it strongly ergodic} if any sequence $\{A_n\}$ of a.i. sets  is trivial, i.e. $\lim\limits_{n\rightarrow\infty}\mu(A_n)(1-\mu(A_n))= 0$.
\end{definition}

Strong ergodicity only depends on the measure class of $\mu$. This observation allows to extend the notion of strong ergodicity to actions on infinite measure spaces.
A nonsingular action $\Gamma\curvearrowright (X,\mu)$ on a standard (possibly infinite) measure space $(X,\mu)$  is said to be strongly ergodic if it is strongly ergodic with respect to a probability measure $\mu_0$ which is equivalent to  $\mu$.

\begin{definition}\cite{Zi81} Let $\Gamma$ be a countable group. A nonsingular action $\Gamma\curvearrowright (X,\mu)$ is said to have {\it property (T)} if any  cocycle $c:\Gamma\times X\rightarrow\mathcal U(\mathcal H)$  into the unitary group of a Hilbert space $\mathcal H$ which admits a sequence of almost invariant unit vectors necessarily has an invariant unit vector.
\end{definition}

For the notions of almost invariant and invariant unit vectors, see Section \ref{pro(T)}. For now,  note that a p.m.p. action $\Gamma\curvearrowright (X,\mu)$  has property (T) if and only if the acting group $\Gamma$ has property (T) of Kazhdan. This fact is however no longer true if the action is not p.m.p.

\subsection{Orbit equivalence rigidity for translation actions}
We are now ready to state our main results. 
Our first result is a ``locally compact analogue" of \cite[Theorem A and Corollary 6.3]{Io13} which established similar statements in the case $G$ and $H$ are profinite or connected compact groups.

\begin{main}[OE rigidity, I]\label{thmA} Let $G$ be a simply connected l.c.s.c. group and $\Gamma<G$ a countable dense subgroup. Assume that the translation action $\Gamma\curvearrowright (G,m_G)$ is strongly ergodic.   
Let $H$ be a simply connected l.c.s.c. group and $\Lambda<H$ a countable subgroup. 

Then the actions $\Gamma\curvearrowright (G,m_G)$ and $\Lambda\curvearrowright (H,m_H)$ are SOE if and only if  there exists a topological isomorphism $\delta:G\rightarrow H$ such that $\delta(\Gamma)=\Lambda$.
\end{main}

We continue with two cocycle and orbit equivalence superrigidity results for translation actions on locally compact groups.  These are analogous to the results proved in  \cite{Io08,Fu09} in the case of  translation actions  of property (T) groups $\Gamma$ on compact groups $G$. When $G$ is locally compact but not compact, the assumption that $\Gamma$ has property (T) needs to be replaced with the stronger assumption that the
(infinite measure preserving) action $\Gamma\curvearrowright (G,m_G)$ has property (T).
 The necessity of imposing a property (T) condition on the action was inspired by an analogous situation for weakly mixing s-malleable measure preserving actions $\Gamma\curvearrowright (X,\mu)$. These actions  were originally shown to be $\mathcal U_{\text{fin}}$-cocycle superrigid whenever $\Gamma$ has property (T) and $\mu(X)<+\infty$ \cite{Po05}.  Later on, cocycle superrigidity for $\Gamma\curvearrowright (X,\mu)$ was proved for possibly infinite measure spaces $(X,\mu)$, whenever the diagonal action $\Gamma\curvearrowright (X\times X,\mu\times\mu)$ has property (T) and is  weakly mixing \cite{PV08}.

\begin{main}[Cocycle superrigidity]\label{thmB}

Let $G$ be a simply connected l.c.s.c.  group and $\Gamma<G$  a countable dense subgroup. Assume that there exists a subgroup $\Gamma_1<\Gamma$ such that $g\Gamma_1 g^{-1}\cap\Gamma_1$ is dense in $G$, for all $g\in\Gamma$, and the translation  action $\Gamma_1\curvearrowright (G,m_G)$ has property (T). 
 
 Let $\Lambda$ be a countable group
and  $w:\Gamma\times G\rightarrow\Lambda$  a  cocycle. 

Then there exist a homomorphism $\delta:\Gamma\rightarrow\Lambda$ and a Borel map $\phi:G\rightarrow\Lambda$
 such that we have $w(g,x)=\phi(gx)\delta(g)\phi(x)^{-1}$, for all $g\in\Gamma$ and almost every $x\in G$.

\end{main}

As a consequence of Theorem \ref{thmB}, we are able to describe all nonsingular actions that are SOE to $\Gamma\curvearrowright (G,m_G)$. More precisely, we prove that any such action is obtained from $\Gamma\curvearrowright (G,m_G)$ by taking quotients and inducing, as explained in the following example:

\begin{example}\label{ex1}  Let $\Gamma\curvearrowright (X,\mu)$ be a nonsingular action of a countable group $\Gamma$. 
For an equivalence relation $\mathcal R$ on a set $Y$, we denote by $\mathcal R|Z:=\mathcal R\cap (Z\times Z)$ the {\it restriction} of $\mathcal R$ to a subset  $Z\subset Y$.
\begin{enumerate}
\item  Let $\Gamma_0$ be a normal subgroup of $\Gamma$ whose action on $X$ has a fundamental domain, i.e. there is a measurable set $X_0$ such that $X=\cup_{g\in\Gamma_0}gX_0$ and $\mu(gX_0\cap X_0)=0$, for all $g\in\Gamma_0\setminus\{e\}$. Then the map $\pi:X_0\rightarrow\Gamma_0\backslash X$ given by $\pi(x)=\Gamma_0x$ witnesses an isomorphism between the restriction of  $\mathcal R(\Gamma\curvearrowright X)$ to $X_0$ and $\mathcal R(\Gamma/\Gamma_0\curvearrowright\Gamma_0\backslash X)$. As a consequence, the action $\Gamma/\Gamma_0\curvearrowright \Gamma_0\backslash X$ is stably orbit equivalent to $\Gamma\curvearrowright X$. Moreover, if $\Gamma\curvearrowright (X,\mu)$ is ergodic, infinite measure preserving and $\mu(X_0)=+\infty$, then these actions are orbit equivalent.

\item Let $\Gamma_1$ be a countable group which contains $\Gamma$. Consider the action of $\Gamma$ on $\Gamma_1\times X$ given by $g\cdot (h,x)=(hg^{-1},gx)$ and the quotient space $X_1:=(\Gamma_1\times X)/\Gamma$.
 Then the {\it induced action} $\Gamma_1\curvearrowright X_1$ given by $g\cdot (h,x)\Gamma=(gh,x)\Gamma$ is stably orbit equivalent to $\Gamma\curvearrowright X$. Indeed, if $X_0:=(\Gamma\times X)/\Gamma$, then the restriction of $\mathcal R(\Gamma_1\curvearrowright X_1)$ to $X_0$ is isomorphic to $\mathcal R(\Gamma\curvearrowright X)$.
\end{enumerate}
\end{example}

\begin{main}[OE superrigidity]\label{thmC}
Let $G$ be a simply connected l.c.s.c.  group and $\Gamma<G$  a countable dense subgroup. Assume that there exists a subgroup $\Gamma_1<\Gamma$ such that  $g\Gamma_1 g^{-1}\cap\Gamma_1$ is dense in $G$, for all $g\in\Gamma$, and the translation  action $\Gamma_1\curvearrowright (G,m_G)$ has property (T). 

Let $\Lambda\curvearrowright (Y,\nu)$ be an arbitrary free ergodic nonsingular action of an arbitrary countable group $\Lambda$.

Then $\Lambda\curvearrowright Y$ is SOE to $\Gamma\curvearrowright (G,m_G)$ if and only if we can find a central discrete subgroup $\Gamma_0<G$ which is contained in $\Gamma$, a subgroup $\Lambda_0<\Lambda$ and a $\Lambda_0$-invariant measurable set $Y_0\subset Y$  such that 
$\Gamma/\Gamma_0\curvearrowright G/\Gamma_0$ is conjugate to $\Lambda_0\curvearrowright Y_0$ and the action $\Lambda\curvearrowright Y$ is induced from $\Lambda_0\curvearrowright Y_0$.
\end{main}

The above theorems show that several known rigidity phenomena for translation actions on compact groups admit analogues in the case of translation actions on locally compact non-compact  groups. This leads to the question: to what extent are  these two classes of actions related. The following result provides an answer to this question.

\begin{mprop}[Weak compactness]\label{propD}
Let $G$ be a l.c.s.c. group and  $\Gamma<G$ a countable dense subgroup. Let $A\subset G$ be a measurable set with $0<m_G(A)<+\infty$. Endow $A$ with the probability measure obtained by restricting and rescaling $m_G$. 
 
 Then the countable ergodic p.m.p. equivalence relation $\mathcal R(\Gamma\curvearrowright G)|A$ is weakly compact.
\end{mprop}

Here, we are using N. Ozawa and S. Popa's notion of weak compactness for equivalence relations \cite{OP07} (see Definition \ref{weakcomp}).

\subsection{Orbit equivalence rigidity for general translation actions}
The results stated so far apply to ``simple" translation actions $\Gamma\curvearrowright (G,m_G)$.  Next, assuming that $G$ is a real algebraic group, we present two rigidity results which apply to fairly general translation actions $\Gamma\curvearrowright (G/\Sigma,m_{G/\Sigma})$. 

\begin{main}[OE rigidity, II]\label{thmE} Let $G$ be a connected real algebraic group with trivial center, $\Sigma<G$ a discrete subgroup and $\Gamma<G$  a countable dense subgroup. Assume that the translation action $\Gamma\curvearrowright (G,m_G)$ is strongly ergodic.
Let $H$ be a connected semisimple real algebraic group with trivial center, $\Delta<H$ a discrete subgroup and $\Lambda<H$ a countable subgroup.

Then the actions $\Gamma\curvearrowright (G/\Sigma,m_{G/\Sigma})$ and $\Lambda\curvearrowright (H/\Delta,m_{H/\Delta})$ are SOE if and only if there exist a topological isomorphism $\delta:G\rightarrow H$ and $h\in H$ such that $\delta(\Gamma)=\Lambda$ and $\delta(\Sigma)=h\Delta h^{-1}$.
\end{main}

\begin{remark} In the context of Theorem \ref{thmE}, assume additionally that $G=SL_n(\mathbb R)$, for some $n\geqslant 2$, and that $\Sigma<G$ is a lattice. Then Examples \ref{ex2}-\ref{ex4} below provide many examples of countable dense subgroups $\Gamma<G$ such that the translation action $\Gamma\curvearrowright G$ is strongly ergodic. On the other hand, by \cite[Theorem D]{IS10}, the action $\Gamma\curvearrowright G/\Sigma$ is rigid, in the sense of S. Popa \cite{Po01}.
Altogether, this shows that Theorem \ref{thmE} applies to a large family of rigid actions.
 
In the proof of Theorem \ref{thmE}, we exploit the fact that the translation action $\Gamma\curvearrowright G/\Sigma$ is related to the simple translation action $\Gamma\curvearrowright G$, via the quotient  map $G\rightarrow G/\Sigma$.  Nevertheless,  from the point of view of orbit equivalence, these actions are at opposite ends. Indeed, by Proposition \ref{propD} the restriction of $\mathcal R(\Gamma\curvearrowright G)$ to any set of finite measure is weakly compact. On the other hand, since the action $\Gamma\curvearrowright G/\Sigma$ is rigid, the equivalence relation $\mathcal R(\Gamma\curvearrowright G/\Sigma)$ is not weakly compact  (this can be seen by using the ergodic-theoretic characterization of rigid actions from \cite{Io09}).
\end{remark}

\begin{main}[OE rigidity, III]\label{thmF} 
 Let $G$ and $H$ be connected real algebraic groups with trivial centers. Let $K<G$ and $L<H$  be connected real algebraic subgroups such that $\cap_{g\in G}gKg^{-1}=\{e\}$ and $\displaystyle{\cap_{h\in H}hLh^{-1}=\{e\}}$. Let $\Gamma<G$ and $\Lambda<H$ be countable dense subgroups, and
assume that the translation actions $\Gamma\curvearrowright (G,m_G)$ and $\Lambda\curvearrowright (H,m_H)$ are strongly ergodic.

Then the actions $\Gamma\curvearrowright (G/K,m_{G/K})$ and $\Lambda\curvearrowright (H/L,m_{H/L})$ are SOE if and only if there exists a topological isomorphism $\delta:G\rightarrow H$ and $h\in H$ such that $\delta(\Gamma)=\Lambda$ and $\delta(K)=hLh^{-1}$.
\end{main}

Theorems \ref{thmE} and \ref{thmF} both require that $\Gamma<G$ is a countable dense subgroup. It would be interesting to decide whether  
similar results hold under the less restrictive assumption that the action $\Gamma$ on the respective homogeneous spaces is ergodic.

\subsection{Strong ergodicity and property (T) for translation actions} In view of the above results, it is natural to wonder when are translation actions  strongly ergodic or have property (T)? The next two results provide necessary conditions for a translation action $\Gamma\curvearrowright (G,m_G)$ to be strongly ergodic or have property (T).
Recall that a p.m.p. action $G\curvearrowright^{\sigma} (X,\mu)$ of a l.c.s.c. group $G$  has {\it spectral gap} if there does not exist a sequence of unit vectors $\xi_n\in L^2_0(X,\mu):=L^2(X,\mu)\ominus\mathbb C1$ such that $\lim\limits_{n\rightarrow\infty}\sup_{g\in K}\|\sigma_g(\xi_n)-\xi_n\|_2=0$, for every compact set $K\subset G$.

\begin{mprop}\label{propG} Let $G=G_1\times G_2$ be a product of two l.c.s.c. groups and $p:G\rightarrow G_1$ be the quotient homomorphism.
Let $\Gamma<G$ be a lattice such that $p(\Gamma)<G_1$ is dense.  Consider the translation action $\Gamma\curvearrowright (G_1,m_{G_1})$ given by $g\cdot x=p(g)x$, for all $g\in\Gamma, x\in G_1$.

Then we have the following:
\begin{enumerate}
\item If the action $G_2\curvearrowright (G/\Gamma,m_{G/\Gamma})$ has spectral gap, then the translation action $\Gamma\curvearrowright (G_1,m_{G_1})$ is strongly ergodic.
\item If $G_2$ has property (T), then the translation action $\Gamma\curvearrowright (G_1, m_{G_1})$ has property (T).
\end{enumerate}
\end{mprop}

\begin{mprop}\label{propH}
Let $G=SL_n(\mathbb R)$ and $K=SO_n(\mathbb R)$, for some $n\geqslant 3$. Let  $\Gamma<G$ be a countable subgroup which is not contained in $K$. Assume that  $\Gamma\cap K$ is dense in $K$ and the translation action $\Gamma\cap K\curvearrowright (K,m_K)$ has spectral gap.

Then the translation action $\Gamma\curvearrowright (G,m_G)$ is strongly ergodic.
\end{mprop}

Next, we combine the previous propositions with known results on spectral gap and property (T) to provide concrete classes of translation actions that are strongly ergodic or have property (T). 

\begin{example}\label{ex2} Let  $\Gamma=SL_n(\mathbb Z[\sqrt{q}])$ and $G=SL_n(\mathbb R)$, where $n,q\geqslant 2$ are integers and $q$ is not a square.  
Then the translation action $\Gamma\curvearrowright (G,m_G)$ is strongly ergodic if $n\geqslant 2$, and has property (T) if $n\geqslant 3$.
Indeed, $\Gamma$ is an irreducible lattice in $G\times G=SL_n(\mathbb R)\times SL_n(\mathbb R)$ and the quotient $(G\times G)/\Gamma$ is not compact.
Then  the action $G\curvearrowright (G\times G)/\Gamma$ has spectral gap, for any $n\geqslant 2$ (see \cite[Theorem 1.12]{KM99}). Moreover, if $n\geqslant 3$, then $G$ has property (T) (see e.g. \cite{Zi84}). The assertion is now a consequence of Proposition \ref{propG}.

\end{example}

\begin{example}\label{ex3}
Let
 $\Gamma=SL_n(\displaystyle{\mathbb Z[\frac{1}{p}]})$ and $G$ be either $SL_n(\mathbb R)$ or $SL_n(\mathbb Q_p)$, for a prime $p$.
Then the translation action $\Gamma\curvearrowright (G,m_G)$ is strongly ergodic if $n\geqslant 2$ and has property (T) if $n\geqslant 3$.
To explain how this assertion follows from the literature, denote $G_1=SL_n(\mathbb R)$ and $G_2=SL_n(\mathbb Q_p)$. Then $\Gamma$ is an irreducible lattice in $\widetilde G:=G_1\times G_2$. 
Consider the unitary representation $\pi:\widetilde G\rightarrow\mathcal U(L^2_0(\widetilde G/\Gamma))$. As is well-known, 
$\pi$ is {\it strongly $L^s$}, for some $s$: the function $\widetilde G\ni g\rightarrow \langle\pi(g)\xi,\eta\rangle$ is $L^s$-integrable, 
for all $\xi,\eta$ belonging to a dense subspace of $L^2_0(\widetilde G/\Gamma)$ (see e.g. \cite[Theorem 1.11]{GMO08}). 
This implies that $\pi^{\otimes_N}$ is contained in a multiple of the left regular representation of $\widetilde G$, for all integers 
$N\geqslant\displaystyle{\frac{s}{2}}$. Since $G_1$ and $G_2$ are non-amenable, the restrictions of $\pi$ to $G_1$ and $G_2$ do not have almost invariant vectors, and therefore the actions of $G_1$ and $G_2$ on $\widetilde G/\Gamma$ have spectral gap. Moreover, 
if $n\geqslant 3$, then $G_1$ and $G_2$ have property (T). The assertion is now a corollary of Proposition \ref{propG}.
\end{example}

\begin{example}\label{ex4}
Let $G=SL_n(\mathbb R)$ and $K=SO_n(\mathbb R)$, for $n\geqslant 3$. Let $\Gamma<G$ be a countable subgroup which contains a matrix $g_0\in G\setminus K$ as well as matrices $g_1,...,g_l\in K$ that have algebraic entries and generate a dense subgroup of $K$. If $n=3$, then the work of J. Bourgain and A. Gamburd \cite{BG06} shows that the action $\Gamma\cap K\curvearrowright (K,m_K)$ has spectral gap. Moreover, the very recent work \cite{BdS14} implies that this statement holds for any $n\geqslant 3$. 
In combination with Proposition \ref{propH}, this shows that
 the translation action $\Gamma\curvearrowright (G,m_G)$ is strongly ergodic.
 
\end{example}

\subsection{Applications}
By using the above examples, one obtains many concrete families of actions to which  Theorems \ref{thmA}-\ref{thmF} apply. Next, we present a sample of applications of our main results. 
For a set of primes  $S$, we denote by $\mathbb Z[S^{-1}]$ the subring of $\mathbb Q$ consisting of rational numbers whose denominators have all prime factors from $S$.

\begin{mcor}\label{corI} Let $m,n\geqslant 2$ and $S,T$ be nonempty sets of primes. 
Then the translation actions $SL_m(\mathbb Z[S^{-1}])\curvearrowright SL_m(\mathbb R)$ and $SL_n(\mathbb Z[T^{-1}])\curvearrowright SL_n(\mathbb R)$ are SOE if and only if $(m,S)=(n,T)$.

More generally, assume that $\Sigma<SL_m(\mathbb R)$, $\Delta<SL_n(\mathbb R)$ are either 
(1) discrete subgroups, or
(2)  connected real algebraic subgroups.
 If the translation actions $SL_m(\mathbb Z[S^{-1}])\curvearrowright SL_m(\mathbb R)/\Sigma$ and $SL_n(\mathbb Z[T^{-1}])\curvearrowright SL_n(\mathbb R)/\Delta$ are SOE, then $(m,S)=(n,T)$.

\end{mcor}

\begin{remark}\label{gefter} Let us explain how in the ``higher rank case",  the first part of Corollary \ref{corI} follows from  R. Zimmer's work \cite{Zi84}, provided that $S$ and $T$ are finite. We start with a general fact.
Let $\Gamma,\Lambda$ be lattices in two products of l.c.s.c. groups  $G=G_1\times G_2$, $H=H_1\times H_2$, such that the left translation actions $\Gamma\curvearrowright G_1$, $\Lambda\curvearrowright H_1$ are SOE. Then \cite[Lemma 6]{Ge03} implies that the left translation actions $G_2\curvearrowright G/\Gamma$, $H_2\curvearrowright H/\Lambda$ are SOE. 

Let $m,n\geqslant 2$ be integers  and $S,T$ finite sets of primes such that  $SL_m(\mathbb Z[S^{-1}])\curvearrowright SL_m(\mathbb R)$ and $SL_n(\mathbb Z[T^{-1}])\curvearrowright SL_n(\mathbb R)$ are SOE. Assume that  $m\geqslant 3$ or $|S|\geqslant 2$.
Then $G_{m,S}:=\prod_{p\in S}SL_m(\mathbb Q_p)$ has rank $\geqslant 2$.
 Recall that $SL_m(\mathbb Z[S^{-1}])$ sits diagonally as a lattice in $SL_m(\mathbb R)\times G_{m,S}$. The previous paragraph then implies that the translation actions $G_{m,S}\curvearrowright (SL_m(\mathbb R)\times G_{m,S})/SL_m(\mathbb Z[S^{-1}])$ and $G_{n,T}\curvearrowright (SL_n(\mathbb R)\times G_{n,T})/SL_n(\mathbb Z[T^{-1}])$ are SOE. Applying \cite[Theorem 10.1.8]{Zi84} finally gives that $G_{m,S}\cong G_{n,T}$, and therefore $(m,S)=(n,T)$.

\end{remark}

\begin{remark}\label{zimmer} Let $m,n\geqslant 2$ be integers and $q,r\geqslant 2$ be integers that are not squares. Assume that the actions $SL_m(\mathbb Z[\sqrt{q}])\curvearrowright SL_m(\mathbb R)$ and $SL_n(\mathbb Z[\sqrt{r}])\curvearrowright SL_n(\mathbb R)$ are SOE.
Then the proof of Corollary \ref{corI} shows that $(m,q)=(n,r)$.  
If  $m\geqslant 3$,  then this conclusion also follows from R. Zimmer's strong rigidity theorem \cite{Zi80}. 
Indeed, as in the previous remark, since $SL_m(\mathbb Z[\sqrt{q}])$ sits diagonally as a lattice in $SL_m(\mathbb R)\times SL_m(\mathbb R)$,  the actions $SL_m(\mathbb R)\curvearrowright (SL_m(\mathbb R)\times SL_m(\mathbb R))/SL_m(\mathbb Z[\sqrt{q}])$, $SL_n(\mathbb R)\curvearrowright (SL_n(\mathbb R)\times SL_n(\mathbb R))/SL_n(\mathbb Z[\sqrt{r}])$ are SOE. \cite[Theorem 4.3]{Zi80} now implies that $m=n$ and the involved actions of $SL_m(\mathbb R)$ are conjugate, from which it follows that $q=r$ as well.

The novelty here consists of being able to handle the case $m=n=2$ and conclude that the actions of $SL_2(\mathbb R)$ on $(SL_2(\mathbb R)\times SL_2(\mathbb R))/SL_2(\mathbb Z[\sqrt{q}])$ are mutually non stably orbit equivalent, as $q$ varies through all the positive integers that are not squares.
\end{remark}

Let $\mathcal R$ be an ergodic countable measure preserving equivalence relation on an infinite standard measure space $(X,\mu)$.  
The automorphism group of $\mathcal R$, denoted Aut$(\mathcal R)$, consists of nonsingular isomorphisms $\theta:X\rightarrow X$ such 
that $(\theta\times\theta)(\mathcal R)=\mathcal R$, almost everywhere. Since $\mathcal R$ is measure preserving and ergodic, there is  $\text{mod}(\theta)>0$ such that $\theta_*\mu=\text{mod}(\theta)\;\mu$.
Then $\text{mod}:\text{Aut}(\mathcal R)\rightarrow\mathbb R^*_{+}$ is a homomorphism and its image $\mathcal F(\mathcal R):=\text{mod}(\text{Aut}(\mathcal R))$ is called the fundamental group of $\mathcal R$.

\begin{mcor}\label{corJ}
Let $G=SL_n(\mathbb R)$, for some $n\geqslant 2$. Let $\Gamma<G$ be a countable dense subgroup which contains the center of $G$ such that the translation action $\Gamma\curvearrowright (G,m_G)$ is strongly ergodic.

Then $\mathcal F(\mathcal R(\Gamma\curvearrowright G))=\{1\}$. In other words, any automorphism of $\mathcal R(\Gamma\curvearrowright G)$ preserves  $m_G$.

\end{mcor}

\begin{remark}
If $\Gamma=SL_n(\mathbb Q)$, then Example \ref{ex3} implies the translation action $\Gamma\curvearrowright (G,m_G)$ is strongly ergodic. Corollary \ref{corJ} therefore implies that the fundamental group of $\mathcal R(\Gamma\curvearrowright G)$ is trivial. This solves part {\it (i)} of \cite[Problem 15]{Ge03}. Now, if $H=SL_2(\mathbb Q_p)$, then Example \ref{ex3} also implies that the translation $\Gamma\curvearrowright (H,m_H)$ is strongly ergodic.
The proof of \cite[Corollary 9.2]{Io13} shows that $\mathcal R(\Gamma\curvearrowright H)$  has trivial fundamental group, which answers part {\it (ii)} of \cite[Problem 15]{Ge03}.
\end{remark}

We end the introduction with an OE superrigidity result which describes all the actions that are SOE to  $PSL_m(\mathbb Z[S^{-1}])\curvearrowright PSL_m(\mathbb R)$ (or, equivalently, to $SL_m(\mathbb Z[S^{-1}])\curvearrowright SL_m(\mathbb R)$), whenever $m\geqslant 3$ and $S$ is a nonempty set of primes. More generally, we have:

\begin{mcor}\label{corK}
Let $m\geqslant 3$ be an integer, $\Sigma<PSL_m(\mathbb R)$ a lattice, and $S$ a nonempty set of primes. If $\Lambda\curvearrowright (Y,\nu)$ is an arbitrary free ergodic nonsingular action of an arbitrary countable group $\Lambda$, then 
\begin{enumerate}
\item $\Lambda\curvearrowright Y$ is SOE to the left translation action $PSL_m(\mathbb Z[S^{-1}])\curvearrowright PSL_m(\mathbb R)$ if and only if we can find a subgroup $\Lambda_0<\Lambda$ and a finite normal subgroup $N<\Lambda_0$ such that 

\begin{itemize}
\item $\Lambda\curvearrowright Y$ is induced from some nonsingular action $\Lambda_0\curvearrowright Y_0$, and
\item  $\Lambda_0/N\curvearrowright Y_0/N$ is conjugate to $PSL_m(\mathbb Z[S^{-1}])\curvearrowright PSL_m(\mathbb R)$.
\end{itemize}

\vskip 0.05in

\item $\Lambda\curvearrowright Y$ is SOE to the left translation action $PSL_m(\mathbb Z[S^{-1}])\curvearrowright PSL_m(\mathbb R)/\Sigma$ if and only if we can find a subgroup $\Lambda_0<\Lambda$ and a finite normal subgroup $N<\Lambda_0$ such that

\begin{itemize}
\item  $\Lambda\curvearrowright Y$ is induced from some nonsingular action $\Lambda_0\curvearrowright Y_0$, and
\item $\Lambda_0/N\curvearrowright Y_0/N$ is conjugate to either $PSL_m(\mathbb Z[S^{-1}])\curvearrowright PSL_m(\mathbb R)/\Sigma$ or the left-right multiplication action $PSL_m(\mathbb Z[S^{-1}])\times\Sigma\curvearrowright PSL_m(\mathbb R)$ given by $(g,\sigma)\cdot x=gx\sigma^{-1}$.
\end{itemize}

\end{enumerate}

\end{mcor}

The second part of Corollary \ref{corK} shows that any free p.m.p. action $\Lambda\curvearrowright Y$ which is SOE to the left translation action $PSL_m(\mathbb Z[S^{-1}])\curvearrowright PSL_m(\mathbb R)/\Sigma$, must be virtually conjugate to it (in the sense of \cite[Definition 1.1]{Fu99}). This adds to the list of OE superrigid p.m.p. actions discovered recently in \cite{Fu99, Po05, Po06, Ki06, Io08, PV08}.

\subsection{Acknowledgements} I would like to thank Alireza Salehi-Golsefidy for helpful discussions on algebraic groups and in particular for pointing out Remark \ref{alireza}.

\section {Preliminaries}

\subsection{Countable nonsingular equivalence relations} Let $(X,\mu)$ be a standard measure space. An equivalence relation $\mathcal R$ on $X$ is called {\it countable nonsingular} if it satisfies the following conditions:
\begin{itemize}
  \item every equivalence class $[x]_{\mathcal R}=\{y\in X|(x,y)\in\mathcal R\}$ is countable, 
  \item $\mathcal R$  is a Borel subset of $X\times X$, and
\item the $\mathcal R$-saturation
 $\cup_{x\in A}[x]_{\mathcal R}$ of any null set $A$ is also a null set. 
\end{itemize}

If $\Gamma\curvearrowright (X,\mu)$ is a nonsingular action of a countable group $\Gamma$, then the {\it orbit equivalence relation} $$\mathcal R(\Gamma\curvearrowright X):=\{(x,y)\in X\times X|\Gamma x=\Gamma y\}$$ is a countable nonsingular equivalence relation. 
Conversely, J. Feldman and C.C. Moore proved that every countable nonsingular equivalence relation arises this way \cite{FM77}.

Let $\mathcal R$ be a countable nonsingular equivalence relation on $(X,\mu)$. We denote by $[\mathcal R]$ the {\it full group} of $\mathcal R$ consisting of all nonsingular isomorphisms $\theta:X\rightarrow X$ such that $(\theta(x),x)\in\mathcal R$, for almost every $x\in X$. We also denote by Aut$(\mathcal R)$ the {\it automorphism group} of $\mathcal R$ consisting of all nonsingular isomorphisms $\theta:X\rightarrow X$ such that $(x,y)\in\mathcal R$ if and only if $(\theta(x),\theta(y))\in\mathcal R$, for almost every $(x,y)\in\mathcal R$. Then $[\mathcal R]$ is a normal subgroup of Aut$(\mathcal R)$. The quotient group is denoted by Out$(\mathcal R)$ and called the {\it outer automorphism group} of $\mathcal R$.

Finally, we say that $\mathcal R$ is {\it measure preserving} if every $\theta\in [\mathcal R]$ preserves $\mu$. If $\mu$ is a probability measure and $\mathcal R$ is measure preserving, then we say that $\mathcal R$ is {\it probability measure preserving} ({\it p.m.p.}). 

\subsection{Strong ergodicity for countable equivalence relations} 
Let $\mathcal R$ be a countable nonsingular equivalence relation on a standard probability space $(X,\mu)$. A sequence  $\{A_n\}$ of measurable subsets of $X$ is said to be {\it asymptotically invariant} ({\it a.i.}) if
$\lim\limits_{n\rightarrow\infty}\mu(\theta(A_n)\Delta A_n)=0$, for all $\theta\in [\mathcal R]$.

\begin{definition} \cite{CW81,Sc81} A countable nonsingular equivalence relation $\mathcal R$ is called {\it strongly ergodic} if any sequence $\{A_n\}$ of a.i. sets  is trivial, i.e. satisfies $\lim\limits_{n\rightarrow\infty}\mu(A_n)(1-\mu(A_n))= 0$. 
\end{definition}
\begin{remark}
It can readily seen that strong ergodicity only depends on the measure class of $\mu$. 
Thus, we say that an equivalence relation $\mathcal R$ on a standard measure space $(X,\mu)$ is strongly ergodic if it is
 strongly ergodic with respect to a probability measure $\mu_0$ which is equivalent to $\mu$. 
It is easy to check that a nonsingular action $\Gamma\curvearrowright (X,\mu)$ is strongly ergodic (in the sense of Definition \ref{strerg}) if and only if its orbit equivalence relation $\mathcal R(\Gamma\curvearrowright X)$ is strongly ergodic.
\end{remark}

Next, we record the following result, whose proof is identical to the proof of \cite[Lemma 2.7]{Io13}, although the latter is only stated in the case of p.m.p. actions.

\begin {lemma}\label{strong}\cite{Io13}  Let $\mathcal R$ be a  countable nonsingular equivalence relation on a standard probability space  $(X,\mu)$. Assume that $\mathcal R$ is strongly ergodic.

Then for every $\varepsilon>0$, we can find $\delta>0$ and $F\subset [\mathcal R]$ finite such that  if a Borel map $\rho:X\rightarrow Y$ with values into a standard Borel space $Y$  satisfies $\mu(\{x\in X|\;\rho(\theta(x))=\rho(x)\})\geqslant 1-\delta$, for all $\theta\in F$, then there exists $y\in Y$ such that $\mu(\{x\in X|\;\rho(x)=y\})\geqslant 1-\varepsilon$.

\end{lemma}

\subsection{Strong ergodicity for actions of locally compact groups} The notion of strong ergodicity for actions of countable groups has a natural extension to actions of locally compact groups.

\begin{definition} Let $G\curvearrowright (X,\mu)$ be a nonsingular ergodic action of a l.c.s.c. group $G$ on a standard probability space $(X,\mu)$.  
\begin{itemize}
\item A sequence  $\{A_n\}$ of measurable subsets of $X$ is said to be {\it asymptotically invariant}  if it satisfies
$\lim\limits_{n\rightarrow\infty}\sup_{g\in K}\mu(gA_n\Delta A_n)=0$, for every compact set $K\subset G$. 
\item The action $G\curvearrowright (X,\mu)$ is called {\it strongly ergodic} if any asymptotically invariant sequence $\{A_n\}$  is trivial, i.e. $\lim\limits_{n\rightarrow\infty}\mu(A_n)(1-\mu(A_n))= 0$. 
\item A nonsingular action $G\curvearrowright (X,\mu)$ of a l.c.s.c. group $G$ on a standard measure space $(X,\mu)$ is strongly ergodic if it is strongly ergodic with respect to some (equivalently, to any) Borel probability measure $\mu_0$ on $X$ which is equivalent to $\mu$.
\end{itemize}
\end{definition}

Next, we provide two constructions of  strongly ergodic actions of locally compact groups.

\begin{lemma}\label{strlc} If $G$ is a l.c.s.c. group, then the left translation action $G\curvearrowright (G,m_G)$ is strongly ergodic.
\end{lemma}

{\it Proof.} Let $\mu$ be a Borel probability measure on $G$ which is equivalent to $m_G$.
Let $\{A_n\}$ be a sequence of measurable subsets of $G$ such that $\sup_{g\in K}\mu(gA_n\Delta A_n)\rightarrow 0$, for every compact $K\subset G$.
Our goal is to show that $\lim\limits_{n\rightarrow\infty}\mu(A_n)(1-\mu(A_n))=0$. 

If $\beta_n=\displaystyle{\int_{G}\mu(g^{-1}A_n\Delta A_n)\;\text{d}\mu(g)}$, then the dominated convergence theorem shows that $\lim\limits_{n\rightarrow\infty}\beta_n=0$.
Let $K_0\subset G$ be a compact set with $m_G(K_0)>0$.  By Fubini's theorem, for all $n$ we have that $$\int_{K_0}\int_{G}|1_{A_n}(gx)-1_{A_n}(x)|\;\text{d}\mu(g)\;\text{d}\mu(x)\leqslant\int_{G}\int_{G}|1_{A_n}(gx)-1_{A_n}(x)|\;\text{d}\mu(x)\;\text{d}\mu(g)=$$ $$\int_{G}\mu(g^{-1}A_n\Delta A_n)\;\text{d}\mu(g)=\beta_n.$$
 
 Thus, for every $n\geqslant 1$, we can find $x_n\in K_0$ such that $\displaystyle{\int_{G}|1_{A_n}(gx_n)-1_{A_n}(x_n)|\;\text{d}\mu(g)}\leqslant\displaystyle{\frac{\beta_n}{\mu(K_0)}}$.
 
 If $x_n\not\in A_n$, then $\mu(A_nx_n^{-1})\leqslant\displaystyle{\frac{\beta_n}{\mu(K_0)}}$. On the other hand, if $x_n\in A_n$, then $1-\mu(A_nx_n^{-1})\leqslant\displaystyle{\frac{\beta_n}{\mu(K_0)}}$.  
In either case, we get that $\mu(A_nx_n^{-1})(1-\mu(A_nx_n^{-1})\leqslant\displaystyle{\frac{\beta_n}{\mu(K_0)}},$ for all $n$. Since $\beta_n\rightarrow 0$, we conclude that $\mu(A_nx_n^{-1})(1-\mu(A_nx_n^{-1}))\rightarrow 0$.

 Since $K_0$ is compact, after passing to a subsequence, we may assume that the sequence $\{x_n\}$ converges to some $x\in K_0$. Let $m_G'$ be a right invariant Haar measure on $G$. Then $\mu$ is absolutely continuous with respect to $m_G'$, hence we can find  $f\in L^1(G,m_G')_{+}$  such that $\text{d}\mu=f\text{d}m_G'$.
 Then $|\mu(A_nx_n^{-1})-\mu(A_nx^{-1})|\leqslant\displaystyle{\int_{G}|f(yx_n)-f(yx)|\text{d}m_G'(y)}$ and therefore $|\mu(A_nx_n^{-1})-\mu(A_nx^{-1})|\rightarrow 0$.
 In combination with the last paragraph, we get that $\mu(A_nx^{-1})(1-\mu(A_nx^{-1}))\rightarrow 0$. This easily implies that $\mu(A_n)(1-\mu(A_n))\rightarrow 0$, which finishes the proof.
 \hfill$\blacksquare$

\begin{lemma}\label{lattice}
Let $G$ be a l.c.s.c. group. Let $\Gamma<G$ be a lattice and $\Gamma\curvearrowright^{\alpha} (X,\mu)$ a nonsingular action on a standard probability space $(X,\mu)$. Assume that the induced action $G\curvearrowright^{\tilde\alpha} (G/\Gamma\times X,m_{G/\Gamma}\times\mu)$ is strongly ergodic.

Then $\alpha$ is strongly egodic.
\end{lemma}

{\it Proof.} Let $\phi:G/\Gamma\rightarrow G$ be a Borel map such that $\phi(x)\Gamma=x$, for all $x\in G/\Gamma$. Let $w:G\times G/\Gamma\rightarrow\Gamma$ be the cocycle given by $w(g,x)=\phi(gx)^{-1}g\phi(x)$, for all $g\in G$ and $x\in G/\Gamma$. Then the induced action $\tilde\alpha$ is given by $g(x,y)=(gx,w(g,x)y)$, for all $g\in G, x\in G/\Gamma$ and $y\in X$. Denote $\tilde\mu=m_{G/\Gamma}\times\mu$.

Consider a sequence $\{A_n\}$ of measurable subsets of $X$ such that $\mu(gA_n\Delta A_n)\rightarrow 0$, for all $g\in\Gamma$. Define $\tilde A_n=G/\Gamma\times A_n$. Then $\sup_{g\in K}\tilde\mu(g\tilde A_n\Delta\tilde A_n)\rightarrow 0$, for every compact set $K\subset G$. Since $\tilde\alpha$ is strongly ergodic, we get that $\tilde\mu(\tilde A_n)(1-\tilde\mu(\tilde A_n))\rightarrow 0$, which implies that $\mu(A_n)(1-\mu(A_n))\rightarrow 0$.
\hfill$\blacksquare$

\subsection{Property (T) for actions and equivalence relations}\label{pro(T)} Next, we recall R. Zimmer's notion of property (T) for actions and equivalence relations. Let $(X,\mu)$ be a standard measure space.

Firstly, let $\Gamma\curvearrowright (X,\mu)$ be a nonsingular action of a countable group $\Gamma$ and let $G$ be a Polish group.
 A Borel map $c:\Gamma\times X\rightarrow G$ is called a {\it cocycle} if  it satisfies the relation $c(gh,x)=c(g,hx)c(h,x)$, for all $g,h\in\Gamma$ and almost every $x\in X$. Two cocycles $c_1,c_2:\Gamma\times X\rightarrow G$ are called {\it cohomologous} if there exists a Borel map $\phi:X\rightarrow G$ such that $c_1(g,x)=\phi(gx)c_2(g,x)\phi(x)^{-1}$, for all $g\in\Gamma$ and almost every $x\in X$.

 \begin{definition}\cite{Zi81}
 A nonsingular action $\Gamma\curvearrowright (X,\mu)$ is said to have {\it property (T)} if every  cocycle $c:\Gamma\times X\rightarrow\mathcal U(\mathcal H)$ into the unitary group of a Hilbert space $\mathcal H$ which admits a sequence of almost invariant unit vectors necessarily has an invariant unit vector, where:
 
 \begin{itemize}

\item A {\it sequence of almost invariant unit vectors} is a sequence of Borel maps $\xi_n:X\rightarrow \mathcal H$ such that for almost every $x\in X$ we have $\lim\limits_{n\rightarrow\infty}\|\xi_n(gx)-c(g,x)\xi_n(x)\|=0$, for all $g\in\Gamma$,  and $\|\xi_n(x)\|=1$, for all $n\geqslant 1$. 

\item An {\it invariant unit vector} is a Borel map $\eta:X\rightarrow \mathcal H$ which satisfies $\eta(gx)=c(g,x)\eta(x)$ and $\|\eta(x)\|=1$, for all $g\in\Gamma$ and almost every $x\in X$.
 \end{itemize}
 \end{definition}
 
Secondly, let $\mathcal R$ be a countable nonsingular equivalence relation on $(X,\mu)$ and $G$ be a Polish group. A Borel map $c:\mathcal R\rightarrow G$ is called a {\it cocycle} if
 for almost every $x\in X$ we have $c(x,y)c(y,z)=c(x,z)$, for all $y,z\in [x]_{\mathcal R}$. Two cocycles $c_1,c_2:\mathcal R\rightarrow G$ are said to be {\it cohomologous} if there exists a Borel map $\phi:X\rightarrow G$ such that $c_1(x,y)=\phi(x)c_2(x,y)\phi(y)^{-1}$, for almost every $x\in X$ and all $y\in [x]_{\mathcal R}$.

 \begin{definition}\cite{Zi81}
  A nonsingular equivalence relation $\mathcal R$ on $(X,\mu)$  is said to have {\it property (T)} if every  cocycle $c:\mathcal R\rightarrow\mathcal U(\mathcal H)$ into the unitary group of a Hilbert space $\mathcal H$ which admits a sequence of almost invariant unit vectors necessarily has an invariant unit vector, where

 \begin{itemize}
 \item A {\it sequence of almost invariant unit vectors} is a sequence of Borel maps $\xi_n:X\rightarrow \mathcal H$ such that for almost every $x\in X$ we have that  $\lim\limits_{n\rightarrow\infty}\|\xi_n(x)-c(x,y)\xi_n(y)\|=0$,  for all $y\in [x]_{\mathcal R}$,  and $\|\xi_n(x)\|=1$,  for every $n\geqslant 1$. 
 \item An {\it invariant unit vector} is a Borel map $\eta:X\rightarrow \mathcal H$ which satisfies  $\eta(x)=c(x,y)\eta(y)$ and $\|\eta(x)\|=1$, for almost every $x\in X$ and all $y\in [x]_{\mathcal R}$.
 \end{itemize}
 
 \end{definition}
 
 \begin{remark}
Let $\mathcal R$ be an ergodic countable nonsingular equivalence relation on $(X,\mu)$. Then the following hold:
 
 \begin{enumerate}
\item  If $\mathcal R$ has property (T), then it is of type II, i.e. there is an $\mathcal R$-invariant measure $\nu$ on $X$ which is equivalent to $\mu$. This follows from the fact that any cocycle $\alpha:\mathcal R\rightarrow\mathbb R_{+}^*$ is trivial (see \cite[Theorem 9.1.1]{Zi84}). 
\item If  $A\subset X$ is a non-null measurable subset, then $\mathcal R|A$ has property (T) if and only if $\mathcal R$ has property (T). 
\item If  $A\subset X$ is a non-null measurable subset, then $\mathcal R|A$ is strongly ergodic if and only if $\mathcal R$ is strongly ergodic.
\item If $\mathcal R$ has property (T), then it is strongly ergodic. In the case when $\mathcal R$ is of type II$_1$, this can be deduced easily from \cite[Th\'{e}or\`{e}me 20]{Pi04}. In general, one can reduce to this case by using the previous three facts.
 \end{enumerate} 
 \end{remark}

 In the proof of Theorem \ref{thmB}, we will need the following result asserting that if a p.m.p. equivalence relation $\mathcal R$ has property (T), then any  almost invariant vector is  close to an invariant vector. 
 
 \begin{proposition}\cite{Po86,Pi04}\label{almost}
 Let $\mathcal R$ be a countable ergodic p.m.p. equivalence relation on a standard probability space $(X,\mu)$. Assume that $\mathcal R$ has property (T). 
 
 Then there exist a constant $\kappa>0$ and a finite set $F\subset [\mathcal R]$ such that the following holds:
 
Let $c:\mathcal R\rightarrow\mathcal  U(\mathcal H)$ be a cocycle, where $\mathcal H$ is a Hilbert space, and $\xi:X\rightarrow \mathcal H$ a Borel map such that $\|\xi(x)\|=1$, for almost every $x\in X$.
Then there is an invariant unit vector $\eta:X\rightarrow \mathcal H$ such that $$\int_{X}\|\eta(x)-\xi(x)\|\;\text{d}\mu(x)\leqslant\kappa\sum_{\theta\in F}\int_{X}\|\xi(\theta(x))-w(\theta(x),x)\xi(x)\|\;\text{d}\mu(x).$$

 \end{proposition} 
Proposition \ref{almost} is an immediate consequence of  \cite[Th\'{e}or\`{e}me 20]{Pi04}.
Proposition \ref{almost}  also follows from \cite[Lemma 4.1.5]{Po86}, after noticing that a countable ergodic p.m.p. equivalence relation $\mathcal R$ has property (T) if and only if its von Neumann algebra $L(\mathcal R)$ has property (T) relative to $L^{\infty}(X)$ (in the sense of \cite[Definition 4.1.3]{Po86}).

\subsection{Algebraic groups and smooth actions} The purpose of this subsection is to establish following result about real algebraic groups that we will need in the proof of Theorem \ref{thmE}.

\begin{lemma}\label{norm}
Let $G$ be a real algebraic group, $H<G$ an $\mathbb R$-subgroup, and denote $N=\cap_{g\in G}\;gHg^{-1}$. 

Then there are an integer $n\geqslant 1$, a $G$-invariant  open conull subset $\Omega\subset  (G/H)^n$ and a Borel map $\pi:\Omega\rightarrow G/N$ such that $\pi(gx)=g\pi(x)$, for all $g\in G$ and every $x\in\Omega$.

Moreover, for every $g\in G\setminus N$, the set $\{x\in G|gxH=xH\}$ has measure zero.

\end{lemma}

This result is most likely known, but for the lack of a reference, we include a proof.
The proof of Lemma \ref{norm} relies on the fact that algebraic actions of algebraic groups on varieties are smooth. 
We begin by recalling the following (see \cite{Zi84}):

\begin{definition} 
A Borel space $X$ is called {\it countably separated} if there exists a sequence of Borel sets which separate points. A Borel action of a l.c.s.c. group $G$ on a standard Borel space $X$ is called {\it smooth} if the quotient Borel space $X/G$ is countably separated.
\end{definition}

The next well-known theorem is due to Borel and Serre (see \cite[Theorem 3.1.3]{Zi84}).
\begin{theorem}\label{smooth}
 If a real algebraic group $G$ acts algebraically on an $\mathbb R$-variety $V$, then the action of $G$ on $V$ is smooth, where $V$ is endowed with its natural Borel structure.
\end{theorem}

{\it Proof of Lemma \ref{norm}.} Let $K=\mathbb C$. By a theorem of Chevalley (see \cite[Proposition 3.1.4]{Zi84}) we can find a regular homomorphism $\pi:G\rightarrow GL_r(K)$ and a point $x\in K^r$ such that $H$ is the stabilizer of  $[x]\in\mathbb P^{r-1}(K)$ in $G$. We may  clearly assume that the set  $\{\pi(g)x|g\in G\}$ spans $K^r$. Thus, we can find $h_1,...,h_r\in G$ such that the vectors $\{\pi(h_1)x,...,\pi(h_r)x\}$ are linearly independent. 

Let $\Omega_0$ be the set of $(g_1,...,g_r)\in G^r$ such that the vectors $\{\pi(g_1)x,...,\pi(g_r)x\}$ are linearly independent. Then $\Omega_0$ is a non-empty, Zariski open subset of $G^r$. In particular, we get $m_{G^r}(G^r\setminus\Omega_0)=0$, where $m_{G^r}$ denotes the Haar measure on $G^r$ obtained by taking the $r$-fold product of $m_G$ with itself. Let $n=r+1$ and $\Omega_1$ be the set of $(g_1,...,g_n)\in G^n$ such that the vectors $\{\pi(g_i)x|1\leqslant i\leqslant n,i\not=j\}$ are linearly independent, for any $j\in\{1,...,n\}$. Then $\Omega_1$ is a $G$-invariant Zariski open subset of $G^n$. Moreover, since $m_{G^r}(G^r\setminus\Omega_0)=0$, we get that $m_{G^n}(G^n\setminus\Omega_1)=0$.

Next, we argue that  if $(g_1,...,g_n)\in\Omega_1$, then $\cap_{i=1}^n g_iHg_i^{-1}=N$. Let $h\in\cap_{i=1}^ng_iHg_i^{-1}$. Notice that $\pi(h)$ stabilizes $[\pi(g_1)x],...,[\pi(g_n)x]$. In other words, $\{\pi(g_1)x,...,\pi(g_n)x\}$ are $n=r+1$ eigenvectors for $\pi(h)\in GL_r(K)$. Since any $r$ vectors from the set $\{\pi(g_1)x,...,\pi(g_n)x\}$ are linearly independent, we get that $\pi(h)=\alpha I$, for some $\alpha\in K\setminus\{0\}$. If $g\in G$, then $\pi(g^{-1}hg)x=\alpha x$, hence $g^{-1}hg\in H$. This shows that $h\in N$, as claimed.

Let $\Omega=\{(g_1H,...,g_nH)|(g_1,...,g_n)\in\Omega_1\}$. Since $\Omega_1$ is a $G$-invariant, open and conull subset of $G^n$, we get that $\Omega$ is a $G$-invariant,  open and conull subset of $(G/H)^n$.
Note that the action of $G$ on $(G/H)^{n}$ descends to an action of $G/N$.
Then the previous paragraph shows that the action $G/N\curvearrowright\Omega$ is free, i.e. $gx\not=x$, for every $g\in G/N$, $g\not=e$ and all $x\in\Omega$.

Now, if $\widetilde H<\widetilde G$ are $\mathbb R$-groups, then $\widetilde G/\widetilde H$ can be endowed with a natural $\mathbb R$-variety structure on which $\widetilde G$ (hence, any $\mathbb R$-subgroup of $\widetilde G$) acts $\mathbb R$-algebraically  (see \cite[Proposition 3.1.4]{Zi84}). In particular,  $(G/H)^n=G^n/H^n$ is an $\mathbb R$-variety on which $G$ acts $\mathbb R$-algebraically. By applying Theorem \ref{smooth} we derive that the action $G\curvearrowright (G/H)^n$ is smooth and hence it admits a {\it Borel selector}. More precisely, there exists a Borel map $s:(G/H)^n\rightarrow (G/H)^n$ satisfying $s(x)\in Gx$ and $s(gx)=s(x)$, for all $g\in G$ and  $x\in (G/H)^n$ (see\cite[Exercise 18.20 iii]{Ke95}).

Since the action $G/N\curvearrowright\Omega$ is free, for every $x\in\Omega$, there is a unique $\pi(x)\in G/N$ such that $x=\pi(x)s(x)$. The map $\pi:\Omega\rightarrow G/N$ clearly satisfies $\pi(gx)=g\pi(x)$, for all $g\in G/N$ and every $x\in\Omega$. In order to finish the proof of the first assertion, it remains to prove that $\pi$ is Borel. To this end, let $F\subset G/N$ be a closed set.

 Let $d$ be a metric on $(G/H)^n$ which gives the Hausdorff  topology. Then  $f:(G/H)^n\rightarrow [0,\infty)$ defined by $f(x)=\inf_{g\in F}d(x,gs(x))$ is Borel. 
 Let $x\in\Omega$ such that $f(x)=0$. Then there is a sequence $\{g_{m}\}_{m\geqslant 1}$ in $F$ such that $g_{m}s(x)\rightarrow x$, as $m\rightarrow\infty$.
  Since the action $G\curvearrowright (G/H)^n$ is smooth and the stabilizer of $x$ is equal to $N$, the map $G/N\ni g\rightarrow gx\in Gx$ is a homeomorphism \cite[Theorem 2.1.14]{Zi84}. Thus, we can find $g\in F$ such that $g_{m}\rightarrow g$, as $m\rightarrow\infty$. This implies that $gs(x)=x$ and hence $\pi(x)=g\in F$.
 Altogether, it follows that $\{x\in\Omega|\pi(x)\in F\}=\Omega\cap\{x\in (G/H)^n|f(x)=0\}$ is a Borel set. Since this holds for any closed set $F\subset G/N$, we get that  $\pi$ is Borel. This finishes the proof of the first assertion.
 
Finally, note that $\cap_{i=1}^ng_iHg_i^{-1}=N$, for almost every $(g_1,...,g_n)\in G^n$. It follows that if $A\subset G$ has positive measure, then $\cap_{g\in A}gHg^{-1}=N$, which implies the moreover assertion.
 \hfill$\blacksquare$

\subsection{Extensions of homomorphisms} 
We end this section by recording a result about extending homomorphisms from a dense subgroup of a l.c.s.c. group to the whole group:

\begin{lemma}\label{ext}  Let $G$ be  a l.c.s.c. group endowed with a Haar measure $m_G$, and $H$ be a Polish group. Let $\Gamma<G$ be  a dense subgroup and $\delta:\Gamma\rightarrow H$ be a homomorphism. Assume that  $\alpha:G\rightarrow H$ is a Borel map such that for all $g\in\Gamma$ we have that $\alpha(gx)=\delta(g)\alpha(x)$, for almost every $x\in G$.

Then $\delta$ extends to a continuous homomorphism $\delta:G\rightarrow H$ and we can find $h\in H$ such that $\alpha(g)=\delta(g)h$, for almost every $g\in G$. 
\end{lemma}

For a proof, see the proof of \cite[Lemma 2.8]{Io13}.

\section{Cocycle rigidity}

The main goal of this section is to establish the following criterion for a cocycle for a translation action $\Gamma\curvearrowright G$ with values in a countable group $\Lambda$ to be cohomologous to a homomorphism $\delta:\Gamma\rightarrow\Lambda$. 

\begin{theorem}\label{cocrig} 
Let $G$ be a  simply connected l.c.s.c. group and  $A$ a Borel subset with $0<m_G(A)<+\infty$.
Let $\Gamma<G$ be a countable dense subgroup and denote $\mathcal R:=\mathcal R(\Gamma\curvearrowright G)$. 

Let $\Lambda$ be a countable group and $w:\mathcal R\rightarrow\Lambda$ be a cocycle.
Suppose that there exist a constant
 $C\in (\displaystyle{\frac{31}{32}},1)$ and a neighborhood $V$ of the identity in $G$ such that 
\begin{equation}\label{uniform}m_G(\{x\in G\;|\; w(\alpha(x)t,xt)=w(\alpha(x),x)\})\geqslant C\;m_G(A),\;\;\text{for all}\;\;\alpha\in [\mathcal R|A]\;\;\text{and every}\;\; t\in V.\end{equation}
Then we can find a homomorphism $\delta:\Gamma\rightarrow\Lambda$ and a Borel map $\phi:G\rightarrow\Lambda$ such that 
we have $w(gx,x)=\phi(gx)\delta(g)\phi(x)^{-1}$, for all $g\in\Gamma$ and almost every $x\in G$.

\end{theorem}

The proof of this result is an adaptation of A. Furman's proof of \cite[Theorem 5.21]{Fu09}. A main ingredient in the proof of Theorem \ref{cocrig} is the following immediate consequence of  \cite[Lemma 2.1]{Io08} which provides a necessary condition for two cocycles to be cohomologous. 

\begin{lemma}\label{coh}
Let $\mathcal R$ be a countable ergodic p.m.p. equivalence relation on a standard probability space $(X,\mu)$. Let $\Lambda$ be a countable group and $w_1,w_2:\mathcal R\rightarrow\Lambda$ be two cocycles. 

Let $C\in (\displaystyle{\frac{31}{32}},1)$ and assume that $\mu(\{x\in X\;|\;w_1(\alpha(x),x)=w_2(\alpha(x),x)\})\geqslant C$, for every $\alpha\in [\mathcal R].$

Then there exists a Borel map $\phi:X\rightarrow\Lambda$ such that we have $\mu(\{x\in X|\phi(x)=e\})>\frac{3}{4}$ and $w_1(x,y)=\phi(x)w_2(x,y)\phi(y)^{-1}$, for almost every $(x,y)\in\mathcal R$.
\end{lemma}

{\it Proof}. 
Let $\Gamma<[\mathcal R]$ be a countable subgroup which generates $\mathcal R$ \cite{FM77}. Since $\mathcal R$ is ergodic, the action $\Gamma\curvearrowright (X,\mu)$ is ergodic.
For $i\in\{1,2\}$,  define a cocycle $v_i:\Gamma\times X\rightarrow\Lambda$ by letting $v_i(g,x)=w_i(gx,x)$.
Since $C>\frac{7}{8}$ and the action $\Gamma\curvearrowright (X,\mu)$ is ergodic, \cite[Lemma 2.1]{Io08} implies that we can find a Borel map $\phi:X\rightarrow\Lambda$ such that $v_1(g,x)=\phi(gx)v_2(g,x)\phi(x)^{-1}$, for all $g\in\Gamma$ and almost every $x\in X$. Thus,  $w_1(x,y)=\phi(x)w_2(x,y)\phi(y)^{-1}$, for almost every $(x,y)\in\mathcal R$.

Moreover,  a close inspection of the proof of \cite [Lemma 2.1]{Io08} shows that $\phi$ verifies the following: there exists $\eta\in L^2(X\times\Lambda,\mu\times c)$
such that  $\phi(x)$ is the unique  $\lambda\in\Lambda$ satisfying 
$|\eta(x,\lambda)|>\frac{1}{2}$, and $\|\eta-1_{X\times\{e\}}\|_2\leqslant\sqrt{2-2C}<\frac{1}{4}$. Here, $c$ denotes the counting measure on $\Lambda$.
 Since $$\mu(\{x\in X|\;|\eta(x,e)|\leqslant\frac{1}{2}\})\leqslant 4\int_{X}|\eta(x,e)-1|^2\;\text{d}\mu(x)\leqslant 4\|\eta-1_{X\times\{e\}}\|_2^2<\frac{1}{4},$$
we conclude that $\mu(\{x\in X|\phi(x)=e\})>\frac{3}{4}$. \hfill$\blacksquare$

\subsection{Proof of Theorem \ref{cocrig}}
As in the proof of \cite[Theorem 5.21]{Fu09}, for  any $t\in G$, we define a new cocycle $w_t:\mathcal R\rightarrow\Lambda$ by letting $w_t(x,y)=w(xt^{-1},yt^{-1})$. 
Since $C>\frac{31}{32}$,  equation \ref{uniform} implies that  the restrictions of $w$, $w_t$ to $\mathcal R|A$ satisfy the assumptions of Lemma \ref{coh}, for any $t\in W:=\{t^{-1}|t\in V\}$.
Thus, by applying Lemma \ref{coh}, for any $t\in W$, we can find a Borel map $\phi_t:A\rightarrow\Lambda$ satisfying  $m_G(\{x\in A|\phi_t(x)=e\})>\frac{3}{4}m_G(A)$ and $w(xt^{-1},yt^{-1})=\phi_t(x)w(x,y)\phi_t(y)^{-1}$, for almost every $(x,y)\in\mathcal R|A$.

Now, since $\Gamma<G$ is dense, we can find a Borel map $\psi:G\rightarrow A$ such that $\psi(x)=x$, for all $x\in A$, and $\psi(x)\in\Gamma x$, for almost every $x\in G$. For $t\in W$, we extend $\phi_t:A\rightarrow\Lambda$ to a map $\phi_t:G\rightarrow\Lambda$ by letting $\phi_t(x)=w(xt^{-1},\psi(x)t^{-1})\phi_t(\psi(x))w(\psi(x),x)$. Then it is easy to check that \begin{equation}\label{coc}w(xt^{-1},yt^{-1})=\phi_t(x)w(x,y)\phi_t(y)^{-1},\;\;\;\text{for almost every}\;\;\; (x,y)\in\mathcal R\;\;\;\text{and every}\;\;\;t\in W.\end{equation}

Next, we claim that whenever $t,s,ts\in W$ we have that \begin{equation}\label{cocycle} \phi_{ts}(x)=\phi_t(xs^{-1})\phi_s(x),\;\;\;\text{for almost every $x\in G$}.\end{equation}
To see this, let $t,s\in W$ such that $ts\in W$, and define $F(x)=\phi_{ts}(x)^{-1}\phi_t(xs^{-1})\phi_s(x)$. Then equation \ref{coc} implies that $F(x)w(x,y)F(y)^{-1}=w(x,y)$, for almost every $(x,y)\in\mathcal R$. This further gives that the set $B=\{x\in G|F(x)=e\}$ is $\mathcal R$-invariant. 
Since  $m_G(\{x\in A|\phi_t(x)=e\})>\frac{3}{4}m_G(A)$, for every $t\in W$, we deduce that $m_G(\{x\in A|F(x)=e\})>0$ and hence $m_G(B)>0$.  Since $\mathcal R$ is ergodic, we conclude that $B=G$, almost everywhere, which proves the claim.

Since $G$ is  simply connected,  the second part of the proof of \cite[Theorem 5.21]{Fu09} shows that  we can find a family of measurable maps $\{\phi_t:G\rightarrow\Lambda\}_{t\in G}$ which extends the family $\{\phi_t:G\rightarrow\Lambda\}_{t\in W}$ defined above in such a way that the identity \ref{cocycle} holds for every $t,s\in G$. 

By arguing exactly as in the end of the proof of \cite[Theorem 5.21]{Fu09} it follows that we can find a measurable map $\phi:G\rightarrow\Lambda$ such that \begin{equation}\label{phi_t}\phi_t(x)=\phi(xt^{-1})\phi(x)^{-1},\;\;\;\text{for almost every $(x,t)\in G\times G$}. \end{equation}
By combining equations \ref{coc} and \ref{phi_t} we get that $\phi(xt^{-1})^{-1}w(xt^{-1},yt^{-1})\phi(yt^{-1})=\phi(x)^{-1}w(x,y)\phi(y)$, for almost every $(x,y)\in\mathcal R$ and almost every $t\in G$.
Let $g\in\Gamma$ and define $L_g:G\rightarrow\Lambda$ by letting $L_g(x)=\phi(gx)^{-1}w(gx,x)\phi(x)$. Then we have that $L_g(xt)=L_g(x)$, for almost every $(x,t)\in G\times G$.
This implies that we can find $\delta(g)\in\Lambda$ such that $L_g(x)=\delta(g)$, for almost every $x\in G$. But then $\delta:\Gamma\rightarrow\Lambda$ must be a homomorphism and  the proof is finished.
\hfill$\blacksquare$

We continue with the following consequence of Theorem \ref{cocrig} which will be a key ingredient in the proof of Theorem \ref{thmA}.

\begin{corollary}\label{cocrig2}
 Let $G$ be a simply connected l.c.s.c.  group and $\Gamma<G$  be a countable dense subgroup such that the action $\Gamma\curvearrowright G$ is strongly ergodic. Let  $\Lambda$ be a countable subgroup of a Polish group $H$ and $w:\Gamma\times G\rightarrow\Lambda$ be a cocycle. Assume that there exists a Borel map $\theta:G\rightarrow H$ such that $w(g,x)=\theta(gx)\theta(x)^{-1}$, for all $g\in\Gamma$ and almost every $x\in G$.
 
 Then there exist a homomorphism $\delta:\Gamma\rightarrow\Lambda$ and a Borel map $\phi:G\rightarrow\Lambda$ such that  we have $w(g,x)=\phi(gx)\delta(g)\phi(x)^{-1}$, for all $g\in\Gamma$ and almost every $x\in G$.
\end{corollary}

Corollary \ref{cocrig2} is a ``locally compact analogue" of \cite[Theorem 4.1]{Io13} (see also \cite[Remark 4.2]{Io13}), with a very similar proof.

{\it Proof.} Denote $\mathcal R:=\mathcal R(\Gamma\curvearrowright G)$ and
let $v:\mathcal R\rightarrow \Lambda$ be the cocycle given by $v(gx,x)=w(g,x)$, for all $g\in\Gamma$ and $x\in G$. 
Fix a Borel set $A\subset G$  with $0<m_G(A)<+\infty$ and let $\mu$ be the probability measure on $A$ given by $\mu(B)=m_G(A)^{-1}m_G(B)$, for every $B\subset A$. Also, fix $\varepsilon \in (0,\frac{1}{64})$.

Since  $\Gamma\curvearrowright G$ is strongly ergodic, $\mathcal R$ and hence $\mathcal R|A$ is strongly ergodic.  By Lemma \ref{strong}, we can find $F\subset [\mathcal R|A]$ finite and $\delta>0$ such that if a Borel map $\rho:A\rightarrow Y$ into a standard Borel space $Y$ satisfies $\mu(\{x\in A|\rho(\alpha(x))=\rho(x)\})\geqslant 1-\delta$, for all $\alpha\in F$, then there exists $y\in Y$ such that $\mu(\{x\in A|\rho(x)=y\})\geqslant 1-\varepsilon$.

Let $\alpha\in F$. Then there is a Borel map $\gamma:A\rightarrow\Gamma$ such that $\alpha(x)=\gamma(x)x$, for all $x\in A$. Hence, for all $x\in A$ and $t\in G$ we have that $v(\alpha(x)t,xt)=v(\gamma(x)xt,xt)=w(\gamma(x),xt)$. Since $\gamma$ and $w$ take countably many values, it is easy to see that $\lim\limits_{t\rightarrow e}\mu(\{x\in A\;|\;w(\gamma(x),xt)=w(\gamma(x),x)\})=1$. Equivalently, we have that $\lim\limits_{t\rightarrow e}\mu(\{x\in A\;|\;v(\alpha(x)t,xt)=v(\alpha(x),x)\})=1$.
We can therefore find a neighborhood $V$ of the identity in $G$ such that \begin{equation}\label{u1}\mu(\{x\in A\;|\;v(\alpha(x)t,xt)=v(\alpha(x),x)\})\geqslant 1-\delta,\;\;\;\text{for all}\;\;\; \alpha\in F\;\;\;\text{and every}\;\;\; t\in V. \end{equation}

For every $t\in V$, we define a Borel map $\rho_t:G\rightarrow H$ by letting $\rho_t(x)=\theta(x)^{-1}\theta(xt)$. Since $\rho_t(\alpha(x))=\rho_t(x)\Longleftrightarrow v(\alpha(x),x)=v(\alpha(xt),xt)$, equation \ref{u1} rewrites as \begin{equation}\label{u2} \mu(\{x\in A\;|\;\rho_t(\alpha(x))=\rho_t(x)\})\geqslant 1-\delta,\;\;\;\text{for all}\;\;\; \alpha\in F\;\;\;\text{and every}\;\;\;t\in V. \end{equation}

By combining \ref{u2} and the above consequence of strong ergodicity, for every $t\in V$, we can find $y_t\in H$ such that we have $\mu(\{x\in G\;|\;\rho_t(x)=y_t\})\geqslant 1-\varepsilon$. Hence, if $\alpha\in [\mathcal R|A]$, then since $\alpha$ preserves $\mu$, we get that $\mu(\{x\in G\;|\;\rho_t(\alpha(x))=\rho_t(x)\})\geqslant 1-2\varepsilon$. From this we further derive that  \begin{equation}\label{u3}\mu(\{x\in A\;|\;v(\alpha(x)t,xt)=v(\alpha(x),x)\})\geqslant 1-2\varepsilon\;\;\;\text{for all}\;\;\; \alpha\in [\mathcal R|A]\;\;\;\text{and every}\;\;\; t\in V. \end{equation}

Since $1-2\varepsilon\in (\frac{31}{32},1)$,   applying Lemma \ref{uniform} yields a homomorphism $\delta:\Gamma\rightarrow\Lambda$ and a Borel map $\phi:G\rightarrow\Lambda$ such that $v(gx,x)=\phi(gx)\delta(g)\phi(x)^{-1}$, for all $g\in\Gamma$ and almost every $x\in G$. This clearly implies the conclusion.
\hfill$\blacksquare$

\subsection{Approximately trivial cocycles}
We end this section by recalling \cite[Lemma 4.2]{Io13}. Note that although \cite[Lemma 4.2]{Io13} is only stated for p.m.p. actions, its proof applies verbatim to nonsingular actions.

\begin{lemma}\label{untwist}\cite{Io13} Let $\Gamma\curvearrowright (X,\mu)$ be a strongly ergodic nonsingular action of a countable group $\Gamma$ on a standard probability space $(X,\mu)$. Let  $H$ be a Polish group and $w:\Gamma\times X\rightarrow H$  a cocycle.
Assume that there exists a sequence of Borel maps $\{\theta_n:X\rightarrow H\}_{n\geqslant 1}$, such that for all $g\in\Gamma$ we have that $\lim\limits_{n\rightarrow\infty}\mu(\{x\in X|w(g,x)=\theta_n(gx)\theta_n(x)^{-1}\})=1$.

Then there exists a Borel map $\theta:X\rightarrow H$ such that $w(g,x)=\theta(gx)\theta(x)^{-1}$, for all $g\in\Gamma$ and almost every $x\in X$.
\end{lemma}

\section{Proof of Theorem  \ref{thmA}}

\subsection{A generalization of Theorem \ref{thmA}}
We begin this section by proving the following more general version of Theorem \ref{thmA}:

\begin{theorem}\label{OE}

 Let $G$ be a connected l.c.s.c. group and $\Gamma<G$ a countable dense subgroup such that the action $\Gamma\curvearrowright (G,m_G)$ is strongly ergodic.   
Let $H$ be a connected l.c.s.c. group and $\Lambda<H$ a countable subgroup. Let $A\subset G, B\subset H$ be non-negligible measurable sets and $\theta:A\rightarrow B$ be a nonsingular isomorphism such that $\theta(\Gamma x\cap A)=\Lambda\theta(x)\cap B$, for almost every $x\in A$.

Suppose that $\widetilde G$ and $\widetilde H$  are simply connected l.c.s.c. groups  together with continuous onto homomorphisms $p:\widetilde G\rightarrow G$ and $q:\widetilde H\rightarrow H$ such that $\ker(p)<\widetilde G$ and $\ker(q)<\widetilde H$ are discrete subgroups. Denote $\widetilde\Gamma=p^{-1}(\Gamma)$ and $\widetilde\Lambda=q^{-1}(\Lambda)$.

Then we can find a topological isomorphism $\delta:\widetilde G\rightarrow \widetilde H$ such that $\delta(\widetilde\Gamma)=\widetilde\Lambda$, a Borel map $\phi:\widetilde G\rightarrow\Lambda$, and $h\in H$ such that $\theta(p(x))=\phi(x)q(\delta(x))h$, for almost every $x\in p^{-1}(A)$.

Moreover, if $G$ and $H$ have trivial centers, then we can find a topological isomorphism $\bar{\delta}:G\rightarrow H$ such that $\bar{\delta}(\Gamma)=\Lambda$, a Borel map $\phi:G\rightarrow\Lambda$, and $h\in H$ such that $\theta(x)=\phi(x)\bar{\delta}(x)h$, for almost every $x\in A$.
\end{theorem}

{\it Proof.}
Since the action $\Gamma\curvearrowright G$ is ergodic, we can extend $\theta$ to a measurable map $\theta:G\rightarrow H$ such that $\theta(\Gamma x)\subset\Lambda\theta(x)$, for almost every $x\in G$.
Define $\tilde\theta:\widetilde G\rightarrow H$ by letting $\tilde\theta(x)=\theta(p(x))$.
Let $w:\widetilde\Gamma\times\tilde G\rightarrow \Lambda$ be the cocycle given by the relation $\tilde\theta(gx)=w(g,x)\tilde\theta(x)$.

 By Example \ref{ex1} (1), the actions $\Gamma\curvearrowright G$ and $\widetilde\Gamma\curvearrowright\widetilde G$ are stably orbit equivalent. Since the action $\Gamma\curvearrowright G$ is strongly ergodic, we deduce that the action $\widetilde\Gamma\curvearrowright\widetilde G$ is also strongly ergodic.
Since  applying Corollary \ref{cocrig2}, we can find a homomorphism $\rho:\widetilde\Gamma\rightarrow\Lambda$ and a Borel map $\phi:\widetilde G\rightarrow\Lambda$ such that $w(g,x)=\phi(gx)\rho(g)\phi(x)^{-1}$, for all $g\in\widetilde\Gamma$ and almost every $x\in\widetilde G$.

Define $\hat\theta:\widetilde G\rightarrow H$ by letting $\hat\theta(x)=\phi(x)^{-1}\tilde\theta(x)$. Then $\hat\theta(gx)=\rho(g)\hat\theta(x)$, for all $g\in\widetilde\Gamma$ and almost every $x\in\widetilde G$. By Lemma \ref{ext}, $\rho$ extends to a continuous homomorphism $\rho:\widetilde G\rightarrow H$ and we can find $h\in H$ such that $\hat\theta(x)=\rho(x)h$, for almost every $x\in\widetilde G$.
From this we get that \begin{equation}\label{u4}\tilde\theta(x)=\phi(x)\rho(x)h,\;\;\;\text{for  almost every}\;\;\; x\in\widetilde G.\end{equation}

We claim that $\ker(\rho)$ is discrete in $\widetilde G$.  Otherwise, we can find a sequence $\{g_n\}$ in $\ker(\rho)\setminus\{e\}$ such that $\lim\limits_{n\rightarrow\infty}g_n=e$. Since $\ker(p)$ is discrete in $\widetilde G$, we may assume that $p(g_n)\not=e$, for all $n$.
By using \ref{u4} we derive that $\tilde\theta(g_nx)=\phi(g_nx)\phi(x)^{-1}\tilde\theta(x)$, for  almost every $x\in\widetilde G$. Since $\Lambda$ is countable, $m_{\widetilde G}(\{x\in\widetilde G|\phi(g_nx)=\phi(x)\;\text{and}\;p(x),p(g_nx)\in A\})>0$, for $n$ large enough. We would thus get that $m_{\widetilde G}(\{x\in\widetilde G|\tilde\theta(g_nx)=\tilde\theta(x)\;\text{and}\;p(x),p(g_nx)\in A\})>0$, for some $n$, contradicting the fact that the restriction of $\theta$ to $A$ is 1-1 and $p(g_n)\not=e$.

Next, let us show that $\rho$ is an onto open map. Let $V$ be a neighborhood of $e\in\widetilde G$.  Let $W\subset\widetilde G$ be a compact subset such that $W^{-1}W\subset V$ and $m_{\widetilde G}(W)>0$. Then $\rho(W)\subset H$ is a compact subset. Moreover, since $\theta:A\rightarrow B$ is a nonsingular isomorphism and $p:\widetilde G\rightarrow G$ is countable-to-1,  we get that $m_H(\widetilde\theta(W))>0$. Since $\phi$ takes countably many values, by using \ref{u4} we deduce that $m_H(\rho(W))>0$. By \cite[Lemma B.4]{Zi84} we derive that $\rho(W)^{-1}\rho(W)$ and therefore $\rho(V)$ contains a neighborhood of $e\in H$. This shows that $\rho$ is an open map.
In particular, $\rho(\widetilde G)$ is an open subgroup of $H$. Since $H$ is connected, we deduce that $\rho(\widetilde G)=H$.

Altogether, we have that  both $q:\widetilde H\rightarrow H$ and $\rho:\widetilde G\rightarrow H$ are covering homomorphisms. Since $\widetilde H$ and $\widetilde G$ are simply connected,  by using the universality property of universal covering groups, we can find a topological isomorphism $\delta:\widetilde G\rightarrow\widetilde H$ such that $q\circ\delta=\rho$. 
Thus,  equation \ref{u4} rewrites as \begin{equation}\label{u5} \theta(p(x))=\tilde\theta(x)=\phi(x)q(\delta(x))h,\;\;\;\text{for almost every}\;\;\; x\in\widetilde G.\end{equation}

Finally, note that if $g\in\widetilde\Gamma$, then $q(\delta(g))=\rho(g)\in\Lambda$ and hence $\delta(g)\in\widetilde\Lambda$. Conversely, let $g\in\widetilde G$ such that $\delta(g)\in\tilde\Lambda$. 
Then $\theta(p(gx))=\phi(gx)q(\delta(gx))\in\Lambda q(\delta(x))=\Lambda\theta(p(x))$, for almost every $x\in\widetilde G$. 
Note that by the construction of $\theta$, for almost every $x\in G$ we have that $\theta(y)\in\Lambda\theta(x)\Rightarrow y\in\Gamma x$.
From this get we that $p(g)p(x)=p(gx)\in\Gamma p(x)$, for almost every $x\in\widetilde G$. Therefore, $p(g)\in\Gamma$ and hence $g\in\widetilde\Gamma$. This shows that $\delta(\widetilde\Gamma)=\widetilde\Lambda$ and finishes the proof of the main assertion.

For the moreover assertion, assume that $G$ and $H$ have trivial centers. Then $\ker(p)=Z(\widetilde G)$ and $\ker(q)=Z(\widetilde H)$, and therefore $\delta$ descends to a topological isomorphism $\bar{\delta}:G\rightarrow H$. It is now clear that $\phi$ factors through the map $p:\widetilde G\rightarrow G$, and the moreover assertion follows.
 \hfill$\blacksquare$

\subsection{Proof of Theorem \ref{thmA}}
 Since $G,H$ are simply connected, Theorem \ref{thmA} follows 
 by applying Theorem \ref{OE} in the case $\widetilde G=G$, $\widetilde H=H$. \hfill$\blacksquare$

\subsection{The outer automorphism of $\mathcal R(\Gamma\curvearrowright G)$}
Theorem \ref{OE} also allows us to compute the outer automorphism group of $\mathcal R(\Gamma\curvearrowright G)$. To state this precisely, for a l.c.s.c. group $G$ and a subgroup $\Gamma$, we denote by  Aut$(\Gamma<G)$ the group of topological automorphisms $\delta$  of $G$ such that $\delta(\Gamma)=\Gamma$. 

\begin{corollary} Let $G$ be a connected l.c.s.c. group and $\Gamma<G$ be a countable dense subgroup such that the action $\Gamma\curvearrowright (G,m_G)$ is strongly ergodic.   Suppose that $\widetilde G$ is a simply connected l.c.s.c. group together with an onto continuous homomorphism $p:\widetilde G\rightarrow G$ such that $\ker(p)$ is discrete in $\widetilde G$. Denote 
$\mathcal R=\mathcal R(\Gamma\curvearrowright G)$ and $\widetilde\Gamma=p^{-1}(\Gamma)$.

Consider the semidirect product $L:=\widetilde G\rtimes\text{Aut}(\widetilde\Gamma<\widetilde G)$. For $\tilde\gamma\in\widetilde\Gamma$, denote by $\text{Ad}(\tilde\gamma)\in\text{Aut}(\widetilde\Gamma<\widetilde G)$ the conjugation with $\tilde\gamma$.
Then $\Delta:=\{(g,\text{Ad}(\tilde\gamma))\;|\;g\in\widetilde G,\tilde\gamma\in\widetilde\Gamma$ with $p(g\tilde\gamma)=e\}$ is a normal subgroup of $L$ and we have the following:
\begin{enumerate}
\item $Out(\mathcal R)\cong L/\Delta$.
\item Assume additionally that every $\delta\in Aut(\widetilde\Gamma<\widetilde G)$  preserves $m_{\widetilde G}$.  Then $\mathcal F(\mathcal R)=\{1\}$.
\end{enumerate}
\end{corollary}

{\it Proof.} (1) We begin the proof of this assertion with a claim. Let $\delta\in$ Aut$(\widetilde\Gamma<\widetilde G)$.

{\bf Claim.}  There exists $\theta_{\delta}\in$ Aut$(\mathcal R)$ such that $\theta_{\delta}(p(x))\in\Gamma p(\delta(x))$, for almost every $x\in\widetilde G$.

{\it Proof of the claim.}
Since $\ker(p)$ is discrete, we can find an open neighborhood $V$ of $e\in\tilde G$ such that $m_G(V)<\infty$ and $p$ is 1-1 on $V^{-1}V\cup\delta(V^{-1}V)$. It follows that the map $\theta:p(V)\rightarrow p(\delta(V))$ given by $\theta(p(x))=p(\delta(x))$ is well-defined and 1-1. Since $\delta$ scales the Haar measure $m_G$,  $\delta$ and hence further $\theta$ is nonsingular. Moreover, for all $x,y\in p(V)$ we have that $\Gamma x=\Gamma y$ if and only if $\Gamma\theta(x)=\Gamma\theta(y)$.

Let us argue that $\theta$ extends to an automorphism $\theta_{\delta}\in$ Aut$(\mathcal R)$.  Assuming that this is the case, then since $\theta_{\delta}(p(x))=p(\delta(x))$, for all $x\in V$, it is easy to show that $\theta_{\delta}$ satisfies the claim.

To construct $\theta_{\delta}$, we consider two cases. Firstly, assume that $G$  is compact. Then $m_G$ is a finite measure, hence $\delta$ and $\theta$ preserve $m_G$. Our claim now follows from the proof of \cite[Lemma 2.2]{Fu03}. 
Secondly, suppose that $G$ is locally compact but not compact.
Then $m_G$ is an infinite measure. Since $m_G(V)<\infty$, $m_G(\delta(V))<\infty$ and $\mathcal R$ is ergodic, we can find sequences of disjoint measurable subsets $\{X_i\}_{i=0}^{\infty},\{Y_i\}_{i=0}^{\infty}$ of $G$ and of elements $\{\alpha_i\}_{i=1}^{\infty}$, $\{\beta_i\}_{i=1}^{\infty}$ in $[\mathcal R]$ such that 

\begin{itemize}
\item $X_0=p(V)$ and $Y_0=p(\delta(V))$. 
\item $X_i=\alpha_i(X_0)$ and $Y_i=\beta_i(Y_0)$, for all $i\geqslant 1$.
\item $G=\cup_{i\geqslant 0}X_i=\cup_{i\geqslant 0}Y_i$, almost everywhere.
\end{itemize}

We define $\theta_{\delta}$ by letting $\theta_{\delta}(x)=\theta(x)$, if $x\in X_0$, and $\theta_{\delta}(x)=\beta_i\theta\alpha_i^{-1}(x)$, if $x\in X_i$, for $i\geqslant 1$. 
\hfill$\square$

Next, let $\varepsilon: $ Aut$(\mathcal R)\rightarrow$ Out$(\mathcal R)$ be the quotient homomorphism. If $\delta\in$ Aut$(\widetilde\Gamma<\widetilde G)$, then $\varepsilon(\theta_{\delta})$ only depends on $\delta$ (and not on the choices made in the proof of the claim). This allows to define a homomorphism $\rho_1:$ Aut$(\widetilde\Gamma<\widetilde G)\rightarrow$ Out$(\mathcal R)$ by letting $\rho_1(\delta)=\varepsilon(\theta_{\delta})$. Further, for $g\in\widetilde G$, we define $\theta_g\in$ Aut$(\mathcal R)$ by letting $\theta_g(x)=xp(g^{-1})$, for every $x\in G$. Then we consider the homomorphism $\rho_2:\widetilde G\rightarrow$ Out$(R)$ given by $\rho_2(g)=\varepsilon(\theta_g)$. It is easy to check that $\rho_1(\delta)\rho_2(g)\rho_1(\delta)^{-1}=\rho_2(\delta(g))$, for all $\delta\in$ Aut$(\widetilde\Gamma<\widetilde G)$ and every $g\in\widetilde G$. 

We can therefore define a homomorphism $\rho:L\rightarrow$ Out$(\mathcal R)$ by letting $\rho(g,\delta)=\rho_1(g)\rho_2(\delta)=\varepsilon(\theta_g\theta_{\delta})$, for every $g\in\widetilde G$ and $\delta\in$ Aut$(\widetilde\Gamma<\widetilde G)$. Theorem \ref{OE} immediately gives that $\rho$ is onto.
 If $g\in\widetilde G$ and $\tilde\gamma\in\widetilde\Gamma$ are such that $p(g\tilde\gamma)=e$, then $\theta_g\theta_{\text{Ad}(\tilde\gamma)}(p(x))\in\Gamma p(\tilde\gamma x\tilde\gamma^{-1}g^{-1})=\Gamma p(x)$, for almost every $x\in\widetilde G$. This shows that $\Delta\subset\ker(\rho)$.

Conversely, let $(g,\delta)\in\ker(\rho)$. Thus, $\theta_g\theta_{\delta}\in [\mathcal R]$, hence $p(\delta(x)g^{-1})\in\Gamma p(x)$, for almost every $x\in\widetilde G$.
We derive that there exists $\gamma\in\Gamma$ such that $A=\{x\in\widetilde G|p(\delta(x)g^{-1})=\gamma p(x)\}$ has positive measure. Since $p(\delta(xy^{-1}))=\gamma p(xy^{-1})\gamma^{-1}$, for all $x,y\in A$, the subgroup $\{x\in\widetilde G|p(\delta(x))=\gamma p(x)\gamma^{-1}\}$  of $\widetilde G$ has positive measure. Since $\widetilde G$ is connected, we conclude that $p(\delta(x))=\gamma p(x)\gamma^{-1}$, for all $x\in\widetilde G$. Let $\tilde\gamma\in\widetilde\Gamma$ such that $p(\tilde\gamma)=\gamma$. Then $p(\delta(x))=p(\tilde\gamma x\tilde\gamma^{-1})$, for all $x\in\widetilde G$. Since $\ker(p)$ is discrete and $\widetilde G$ is connected, we deduce that $\delta=\text{Ad}(\tilde\gamma)$.
Therefore, for almost every $x\in A$ we have that $p(\tilde\gamma x\tilde\gamma^{-1}g^{-1})=p(\delta(x)g^{-1})=\gamma p(x)=p(\tilde\gamma x)$. We further get that $p(g\tilde\gamma)=e$ and hence $(g,\delta)=(g,\text{Ad}(\tilde\gamma))\in\Delta$. This completes the proof of assertion (1).

(2)  Assume  that every automorphism $\delta\in$ Aut$(\widetilde\Gamma<\widetilde G)$  preserves $m_{\widetilde G}$.  Then Ad$(\tilde\gamma)$ preserves $m_{\widetilde G}$, for every $\tilde\gamma\in\widetilde\Gamma$. Since $\widetilde\Gamma<\widetilde G$ is dense, we deduce that  Ad$(g)$ preserves $m_{\widetilde G}$, for every $g\in\widetilde G$. This implies that $\widetilde G$ is unimodular. It follows that the map $\widetilde G\ni x\rightarrow\delta(x)g\in\widetilde G$ preserves $m_{\widetilde G}$, for every $\delta\in$ Aut$(\widetilde\Gamma<\widetilde G)$ and all $g\in\widetilde G$.
Since the homomorphism $\rho:L\rightarrow$ Out$(\mathcal R)$ defined above is onto, we get that every automorphism of $\mathcal R$ preserves $m_{G}$ and assertion (2) follows.
\hfill$\blacksquare$

\subsection{Borel reducibility rigidity} We end this section with an analogue of Theorem \ref{OE} for Borel reducibility.
Let $\mathcal R,\mathcal S$ be countable Borel equivalence relations on standard Borel spaces $X,Y$. We say that $\mathcal R$ is {\it Borel reducible to} $\mathcal S$ if there exists a Borel map $\theta:X\rightarrow Y$ such that $(x,y)\in\mathcal R$ if and only if $(\theta(x),\theta(y))\in\mathcal S$.

\begin{theorem} Let $G$ be a connected l.c.s.c. group and $\Gamma<G$ be a countable dense subgroup such that the action $\Gamma\curvearrowright (G,m_G)$ is strongly ergodic.   Suppose that $\widetilde G$ is a simply connected l.c.s.c. group together with an onto continuous homomorphism $p:\widetilde G\rightarrow G$ such that $\ker(p)$ is discrete in $\widetilde G$. Denote 
 $\widetilde\Gamma=p^{-1}(\Gamma)$.
Let $H$ be a  l.c.s.c. group and $\Lambda<H$ a countable subgroup. 

Then $\mathcal R(\Gamma\curvearrowright G)$ is Borel reducible to $\mathcal R(\Lambda\curvearrowright H)$ if and only if there exists a continuous homomorphism $\delta:\widetilde G\rightarrow H$ such that  $\delta^{-1}(\Lambda)=\widetilde\Gamma$.
\end{theorem}

{\it Proof.} To see the {\it if} part, assume that $\delta:\widetilde G\rightarrow H$ is a continuous homomorphism such that  $\delta^{-1}(\Lambda)=\widetilde\Gamma$. Let $r:G\rightarrow\widetilde G$ be a Borel map such that $p(r(x))=x$, for all $x\in G$. Then it is routine to check that $\theta:G\rightarrow H$ given by $\theta(x)=\delta(r(x))$ is the desired Borel reduction.

For the {\it only if} part, assume that there exists a Borel map $\theta:G\rightarrow H$ such that $x\in\Gamma y$ if and only if $\theta(x)\in\Lambda\theta(y)$. Then the proof of Theorem \ref{OE} shows that we can find a Borel map $\phi:\widetilde G\rightarrow\Lambda$, a continuous homomorphism $\delta:\widetilde G\rightarrow H$ satisfying $\delta(\widetilde\Gamma)\subset\Lambda$, and $h\in H$ such that $\theta(p(x))=\phi(x)\delta(x)h$, for almost every $x\in\widetilde G$. In order to finish the proof, we only need to argue that $\delta^{-1}(\Lambda)\subset\widetilde\Gamma$. To this end, let $g\in G$ such that $\delta(g)\in\Lambda$. Then for almost every $x\in G$, hence for some $x\in G$, we have that $\theta(p(gx))=\phi(gx)\delta(gx)h=(\phi(gx)\delta(g)\phi(x)^{-1})\theta(p(x))\in\Lambda\theta(p(x))$. This implies that $p(g)p(x)=p(gx)\in\Gamma p(x)$ and therefore $p(g)\in\Gamma$, hence $g\in\tilde\Gamma$, as desired.
\hfill$\blacksquare$

\section{Proof of Theorem \ref{thmB}} 
The proof of Theorem \ref{thmB} follows closely the proof of  \cite[Theorem 5.21]{Fu09}. The idea behind the proof is based on S. Popa's deformation/rigidity theory.  
Roughly speaking, we exploit the tension between the rigidity coming from the property (T) of the action $\Gamma\curvearrowright G$ and the deformation associated to the right multiplication action of $G$ on itself. 

\subsection{Proof of Theorem \ref{thmB}}
 Assume first that $G$ is simply connected and the action $\Gamma\curvearrowright (G,m_G)$ has property (T). Let $w:\Gamma\times G\rightarrow\Lambda$ be a cocycle, where $\Lambda$ is a countable group. Our goal is to show that $w$ is cohomologous to a homomorphism $\delta:\Gamma\rightarrow\Lambda$.

Let $A\subset G$ be a Borel subset with $0<m_G(A)<+\infty$.
  Since $\Gamma\curvearrowright G$ has property (T), we get that $\mathcal R:=\mathcal R(\Gamma\curvearrowright G)$ has property (T) and that $\mathcal R|A$ has property (T). Note that $\mathcal R|A$ preserves the probability measure $\mu$ on $A$ given by $\mu(B)=m_G(A)^{-1}m_G(B)$, for all measurable subsets $B\subset A$.

By Proposition \ref{almost} we can find $\kappa>0$ and a finite set $F\subset [\mathcal R|A]$ such that the following holds: if $c:\mathcal R\rightarrow\mathcal  U(\mathcal H)$ is a cocycle, where $\mathcal H$ is a Hilbert space, and $\xi:A\rightarrow \mathcal H$ is a Borel map satisfying $\|\xi(x)\|=1$, for almost every $x\in A$, then there exists an invariant unit vector $\eta:A\rightarrow \mathcal H$ such that $$\int_{A}\|\eta(x)-\xi(x)\|\;\text{d}\mu(x)\leqslant\kappa\sum_{\theta_0\in F}\int_{A}\|\xi(\theta_0(x))-c(\theta_0(x),x)\xi(x)\|\;\text{d}\mu(x).$$

 Since $\eta$ is an invariant vector and $[\mathcal R|A]$ preserves $\mu$, it follows that for all $\theta\in [\mathcal R|A]$ we have that  \begin{equation}\label{uni}\int_{A}\|\xi(\theta(x))-c(\theta(x),x)\xi(x)\|\;\text{d}\mu(x)\leqslant 2\kappa\sum_{\theta_0\in F}\int_{A}\|\xi(\theta_0(x))-c(\theta_0(x),x)\xi(x)\|\;\text{d}\mu(x).\end{equation}

Next, we let $v:\mathcal R\rightarrow\Lambda$ be given by $v(gx,x)=w(g,x)$.
For $t\in G$, we define a  cocycle $v_t:\mathcal R\rightarrow\Lambda$ by letting $v_t(x,y)=v(xt^{-1},yt^{-1})$. 
Further, we define a cocycle $c_t:\mathcal R\rightarrow\mathcal U(\ell^2(\Lambda))$ by letting $$(c_t(x,y)f)(\lambda)=f(v_t(x,y)^{-1}\lambda v(x,y)),\;\;\text{for all}\;\; (x,y)\in\mathcal R,\;f\in\ell^2(\Lambda)\;\;\text{and}\;\;\lambda\in\Lambda.$$

Let $\xi:A\rightarrow\ell^2(\Lambda)$ be given by $\xi(x)=\delta_e$, for all $x\in A$. Then we have $c_t(x,y)\xi(y)=\delta_{v_t(x,y)v(x,y)^{-1}}$, for all $x,y,t\in G$, where $\{\delta_{\lambda}\}_{\lambda\in\Lambda}$ denotes the usual orthonormal basis of $\ell^2(\Lambda)$. Thus, by applying equation \ref{uni} we get that for all $\theta\in [\mathcal R|A]$ and every $t\in G$ we have that \begin{equation}\label{uni2}\int_{A}\|\delta_e-\delta_{v_t(\theta(x),x)v(\theta(x),x)^{-1}}\|\;\text{d}\mu(x)\leqslant 2\kappa\sum_{\theta_0\in F}\int_{A}\|\delta_e-\delta_{v_t(\theta_0(x),x)v(\theta_0(x),x)^{-1}}\|\;\text{d}\mu(x) \end{equation}

Let $\varepsilon\in (0,\displaystyle{\frac{\sqrt{2}}{32}})$.
Since $\Lambda$ is countable,  $\lim\limits_{t\rightarrow e}\mu(\{x\in A|v_t(\theta(x),x)=v(\theta(x),x)\})=1$, for any $\theta\in [\mathcal R|A]$ (see the proof of Corollary \ref{cocrig2}). Therefore, we can find a neighborhood $V$ of $e\in G$ such that \begin{equation}\label{uni3}\int_{A}\|\delta_e-\delta_{v_t(\theta_0(x),x)v(\theta_0(x),x)^{-1}}\|\;\text{d}\mu(x)\leqslant\frac{\varepsilon}{2\kappa |F|},\;\;\text{for every}\;\; t\in V\;\;\text{and all}\;\;\theta_0\in F.\end{equation}

By combining \ref{uni2} and \ref{uni3} we derive that for all $t\in V$ and every $\theta\in [\mathcal R|A]$, we have that $$\sqrt{2}\;\mu(\{x\in A\;|\;v_t(\theta(x),x)\not=v(\theta(x),x)\})=\int_{A}\|\delta_e-\delta_{v_t(\theta(x),x)v(\theta(x),x)^{-1}}\|\;\text{d}\mu(x)\leqslant\varepsilon.$$  

Thus, we conclude that $\mu(\{x\in A|v_t(\theta(x),x)=v(\theta(x),x)\})\geqslant \displaystyle{1-\frac{\varepsilon}{\sqrt{2}}}$, for all $t\in V$ and  $\theta\in [\mathcal R|A]$.
Since $\displaystyle{1-\frac{\varepsilon}{\sqrt{2}}>\frac{31}{32}}$, the conclusion follows by applying Theorem \ref{cocrig}.

This finishes the proof in the case the action $\Gamma\curvearrowright (G,m_G)$ has property (T). In general, the conclusion is obtained by combining this case with the following lemma.
\hfill$\blacksquare$

\begin{lemma}\label{ergodique}
Let $G$ be a connected l.c.s.c. group, $\Gamma<G$ a countable subgroup, $\Gamma_1<\Gamma$ a subgroup, and $g\in\Gamma$  such that $g\Gamma_1g^{-1}\cap\Gamma_1$ is dense in $G$. Let $\Lambda$ be a countable group, $\delta:\Gamma_1\rightarrow\Lambda$ a homomorphism, and $w:\Gamma\times G\rightarrow\Lambda$ a cocycle such that $w(h,x)=\delta(h)$, for all $h\in\Gamma_1$ and almost every $x\in G$.

Then there is $\lambda\in\Lambda$ such that $w(g,x)=\lambda$, for almost every $x\in G$.
\end{lemma}

{\it Proof.} Let $\Gamma_2=g\Gamma_1 g^{-1}\cap\Gamma_1$ and $\alpha:\Gamma_2\rightarrow\Gamma_1$ given by $\alpha(h)=g^{-1}hg$. If $h\in\Gamma_2$, then $g\alpha(h)=hg$ and the cocycle relation yields that $w(g,\alpha(h)x)\delta(\alpha(h))=w(g\alpha(h),x)=w(hg,x)=\delta(h)w(g,x)$, for almost every $x\in G$.
Let $S$ be the set of $(x,y)\in G\times G$ such that $w(g,x)=w(g,y)$. Then the last identity implies that $S$ is invariant under the diagonal action of $\alpha(\Gamma_2)$ on $G\times G$. Since $\alpha(\Gamma_2)$ is dense in $G$, $S$ must be invariant under the diagonal action of $G$ on $G\times G$. 

Therefore, since $S$ is non-negligible, there is a non-negligible measurable set $T\subset G$ such that
$S=\{(z,zt)|z\in G,t\in T\}$, almost everywhere. As a consequence, the set $T_0$ of all $t\in G$ such that $(z,zt)\in S$, for almost every $z\in G$, is non-negligible.
Since $T_0<G$ is a  subgroup and  $G$ is connected, we derive that $T_0=G$, almost everywhere. Thus, $S=G\times G$, almost everywhere, which clearly implies the conclusion.
\hfill$\blacksquare$

\section{Proof of Theorem \ref{thmC}}
Assume that a translation action $\Gamma\curvearrowright G$ has property (T). In this section, we use Theorem \ref{thmB} to describe the actions that are SOE to $\Gamma\curvearrowright G/\Sigma$, whenever $\Sigma<G$ is a discrete subgroup. In particular, by applying this description in the case $\Sigma=\{e\}$, we deduce Theorem \ref{thmC}.

\begin{theorem}\label{OESUP1}
Let $G$ be a simply connected l.c.s.c.  group and $\Gamma<G$  a countable dense subgroup. Assume that there exists a subgroup $\Gamma_1<\Gamma$ such that $g\Gamma_1 g^{-1}\cap\Gamma_1$ is dense in $G$, for all $g\in\Gamma$, and the translation  action $\Gamma_1\curvearrowright (G,m_G)$ has property (T).
Let $\Sigma<G$ be a discrete subgroup. Let $\Lambda\curvearrowright (Y,\nu)$ be a free ergodic nonsingular action of a countable group $\Lambda$ which is SOE to $\Gamma\curvearrowright (G/\Sigma,m_{G/\Sigma})$.

Then we can find a normal subgroup $\Delta<\Gamma\times\Sigma$, 
a subgroup $\Lambda_0<\Lambda$, and a $\Lambda_0$-invariant Borel subset $Y_0\subset Y$ with $\nu(Y_0)>0$ such that 
\begin{itemize}
\item $\Delta$ is discrete in $G\times G$,
\item the left-right multiplication action $\Delta\curvearrowright G$ admits a measurable fundamental domain,
\item the action $(\Gamma\times\Sigma)/\Delta\curvearrowright G/\Delta$ is conjugate to $\Lambda_0\curvearrowright Y_0$, and 
\item the action $\Lambda\curvearrowright Y$ is induced from $\Lambda_0\curvearrowright Y_0$.
\end{itemize}
\end{theorem}

{\it Proof.}  Let $\Lambda\curvearrowright (Y,\nu)$ be a free nonsingular action which is SOE to 
$\Gamma\curvearrowright (G/\Sigma,m_{G/\Sigma})$. Since the latter action is SOE to the action $\Gamma\times\Sigma\curvearrowright G$, we deduce that $\Lambda\curvearrowright (Y,\nu)$ is SOE  to $\Gamma\times\Sigma\curvearrowright G$. Let $A\subset G$, $B\subset Y$ be non-negligible measurable sets and $\theta:A\rightarrow B$ a nonsingular isomorphism such that $\theta((\Gamma\times\Sigma)x\cap A)=\Lambda\theta(x)\cap B$, for almost every $x\in A$. 
Since $\Gamma\curvearrowright G$ is ergodic, we may extend $\theta$ to a measurable map $\theta:G\rightarrow Y$ such that $\theta((\Gamma\times\Sigma)x)\subset\Lambda\theta(x)$, for almost every $x\in G$.

Define a cocycle $w:\Gamma\times G\rightarrow\Lambda$ by the formula $\theta(gx)=w(g,x)\theta(x)$, for all $g\in\Gamma$ and almost every $x\in G$. 
By applying Theorem \ref{thmB}, we can find a homomorphism $\delta:\Gamma\rightarrow\Lambda$ and a Borel map $\phi:G\rightarrow\Lambda$ such that $w(g,x)=\phi(gx)\delta(g)\phi(x)^{-1}$, for all $g\in\Gamma$ and almost every $x\in G$.

The map $\tilde\theta:G\rightarrow Y$ given by $\tilde\theta(x)=\phi(x)^{-1}\theta(x)$ therefore satisfies \begin{equation}\label{tt}\tilde\theta(gx)=\delta(g)\tilde\theta(x),\;\;\;\text{for all $g\in\Gamma$ and almost every $x\in G$.} \end{equation}

{\bf Claim}. $\delta:\Gamma\rightarrow\Lambda$ extends to a homomorphism $\delta:\Gamma\times\Sigma\rightarrow\Lambda$ such that  $\tilde\theta(gx\sigma^{-1})=\delta(g,\sigma)\tilde\theta(x)$, for all $g\in\Gamma$, $\sigma\in\Sigma$, and almost every $x\in G$.

{\it Proof of the claim.} Fix $\sigma\in\Sigma$.
Since $\tilde\theta(x\sigma^{-1})\in\Lambda\tilde\theta(x)$, for almost every $x\in G$, we can find a Borel map $v:G\rightarrow\Lambda$ such that $\tilde\theta(x\sigma^{-1})=v(x)\tilde\theta(x)$, for almost every $x\in G$.
By combining the fact that the actions of $\Gamma$ and $\Sigma$ on $G$  commute with \ref{tt} and the freeness of the $\Lambda$-action, it follows  that \begin{equation}\label{vv}v(gx)=\delta(g)v(x)\delta(g)^{-1},\;\;\;\text{for all $g\in\Gamma$ and almost every $x\in G$.} \end{equation}

Define $S$ to be the set of all $(x,y)\in G\times G$ such that $v(x)=v(y)$. Equation \ref{vv} gives that  $S$ is invariant under the diagonal action of $\Gamma$ on $G\times G$. Since $\Gamma<G$ is dense, $S$ must be invariant under the diagonal action of $G$ on $G\times G$. By repeating the argument from the proof of Lemma \ref{ergodique}, we derive that $S=G\times G$, almost everywhere.
This implies that $v:G\rightarrow\Lambda$ is a constant function. If we denote this constant by $\delta(\sigma)$, then $\tilde\theta(x\sigma^{-1})=\delta(\sigma)\tilde\theta(x)$, for almost every $x\in G$.  It is then clear that $\delta:\Sigma\rightarrow\Lambda$ is a homomorphism. Moreover, by equation \ref{vv} we get that $\delta(\sigma)$ commutes with $\delta(\Gamma)$. Altogether, the claim follows. \hfill$\square$

Let $\Delta:=\ker(\delta)$. Assume by contradiction that $\Delta$ is not discrete in $G\times G$.  Since $\Sigma<G$ is discrete, it follows there is a sequence $g_n\in\Gamma\setminus\{e\}$ such that $\lim\limits_{n\rightarrow\infty}g_n=e$ and $(g_n,e)\in\Delta$, for all $n$. Fix a Borel set $A_0\subset A$ with $0<m_G(A_0)<\infty$.
Then $\tilde\theta(g_nx)=\tilde\theta(x)$, for all $n$, and almost every $x\in A_0$. Since $\phi$ takes countably many values, we also have that $\lim\limits_{n\rightarrow\infty}m_G(\{x\in A_0|\phi(g_nx)=\phi(x)\})=m_G(A_0)$. By combining these facts, we get that $m_G(\{x\in A_0|\theta(g_nx)=\theta(x)\})>0$, for $n$ large enough. This contradicts the fact that the restriction of $\theta$ to $A$ is $1-1$.

Now, let $\rho:=\theta^{-1}:B\rightarrow A$. Since $\rho(\Lambda y\cap B)=(\Gamma\times\Sigma)\rho(y)\cap A$, for almost every $y\in B$, and 
 the action $\Lambda\curvearrowright Y$ is ergodic, we may extend $\rho$ to a measurable map $\rho:Y\rightarrow G$ such that $\rho(\Lambda y)\subset (\Gamma\times\Sigma)\rho(y)$, for almost every $y\in Y$. Then $\rho(\Lambda\tilde\theta(x))=\rho(\Lambda\theta(x))\subset (\Gamma\times\Sigma)x$, for almost every $x\in G$, and the rest of the assertions are a consequence of the following lemma.
\hfill$\blacksquare$

\begin{lemma}\label{conju}
Let $\Gamma\curvearrowright (X,\mu)$ be a nonsingular action of a countable group $\Gamma$ and $\Lambda\curvearrowright (Y,\nu)$  be a free nonsingular action of a countable group $\Lambda$. Assume that 
there exist nonsingular maps $\theta:X\rightarrow Y$, $\rho:Y\rightarrow X$ and a group homomorphism $\delta:\Gamma\rightarrow\Lambda$ such that $\rho(\Lambda\theta(x))\subset\Gamma x$, for almost every $x\in X$, and $\theta(gx)=\delta(g)\theta(x)$, for all $g\in\Gamma$ and almost every $x\in X$. 

Define $\Gamma_0:=\ker(\delta)$, $\Lambda_0:=\delta(\Gamma)$ and $Y_0:=\theta(X)$.

Then the action $\Gamma_0\curvearrowright X$ admits a measurable fundamental domain, the action $\Gamma/\Gamma_0\curvearrowright X/\Gamma_0$ is conjugate to $\Lambda_0\curvearrowright Y_0$, and the action $\Lambda\curvearrowright Y$ is induced from $\Lambda_0\curvearrowright Y_0$.
\end{lemma}

{\it Proof.} Consider the nonsingular map $\tau:=\rho\circ\theta:X\rightarrow X$. 
Let $X_1\subset X$ be a maximal measurable set such that $\mu(gX_1\cap X_1)=0$, for all $g\in\Gamma_0\setminus\{e\}$. Define $X_0=\cup_{g\in\Gamma_0}gX_1$. If $X_0=X$, then $X_0$ is a fundamental domain for  $\Gamma_0\curvearrowright X$. Assume by contradiction that $\mu(X\setminus X_0)>0$. Since 
$\tau(x)\in\Gamma x$, for almost every $x\in X$, there is a non-negligible subset $X_2\subset X\setminus X_0$ such that $\tau_{|X_2}$ is 1-1. 
Note that  $\tau(gx)=\tau(x)$, for all $g\in\Gamma_0$ and almost every $x\in X$. We deduce that $\mu(gX_2\cap X_2)=0$, for all $g\in\Gamma_0\setminus\{e\}$. 
Since $X_2\subset X\setminus X_0$ and $X_0$ is $\Gamma_0$-invariant, we get that $X_3=X_1\cup X_2$ also satisfies that $\mu(gX_3\cap X_3)=0$, for all $g\in\Gamma_0\setminus\{e\}$. This contradicts the maximality of $X_1$.

Since the action $\Gamma_0\curvearrowright X$ has a measurable fundamental domain, the quotient space $X/\Gamma_0$ endowed with the push-forward of $\mu$ is a standard measure space. Since $\theta(gx)=\theta(x)$, for all $g\in\Gamma_0$ and almost every $x\in X$,  letting $\bar{\theta}:X/\Gamma_0\rightarrow Y_0$ given by $\bar{\theta}(\Gamma_0 x)=\theta(x)$ defines an onto nonsingular map. Moreover, $\bar{\theta}$ is 1-1. Indeed, if $x,y\in X$ are such that $\theta(x)=\theta(y)$, then $\tau(x)=\tau(y)$, hence $\Gamma x=\Gamma y$. Let $g\in\Gamma$ such that $y=gx$. Since $\theta(x)=\theta(y)=\delta(g)\theta(x)$ and the action $\Lambda\curvearrowright Y$ is free, we deduce that $\delta(g)=e$. Hence $g\in\Gamma_0$ and therefore $\Gamma_0 x=\Gamma_0 y$. This shows that $\theta$ is 1-1. It is now clear that $\bar{\theta}$ is a conjugacy between $\Gamma/\Gamma_0\curvearrowright X/\Gamma_0$ and $\Lambda_0\curvearrowright Y_0$.

To see that the action $\Lambda\curvearrowright Y$ is induced from $\Lambda_0\curvearrowright Y_0$, let $h\in\Lambda$ such that $\nu(hY_0\cap Y_0)>0$. Then we can find $x_1, x_2\in X$ such that $y=\theta(x_1)$ and $h^{-1}y=\theta(x_2)$. But then we get that $\theta(x_2)\in\Lambda\theta(x_1)$ and by applying $\rho$ we deduce that $x_2\in\Gamma x_1$. Let $g\in\Gamma$ such that $x_2=gx_1$. Thus, $h^{-1}y=\theta(x_2)=\delta(g)\theta(x_1)=\delta(g)y$. This shows that $h=\delta(g)^{-1}\in\Lambda_0$, which finishes the proof.
\hfill$\blacksquare$

\subsection{Proof of Theorem \ref{thmC}} If $\Delta<\Gamma$ is a normal subgroup which is discrete in $G$, then since $\Gamma<G$ is dense and $G$ is connected,  $\Delta$ must be central in $G$. Using this fact, Theorem \ref{thmC} follows immediately by applying Theorem \ref{OESUP1} in the case $\Sigma=\{e\}$.
\hfill$\blacksquare$

\section{Proof of Proposition \ref{propD}}
In preparation for the proof of Proposition \ref{propD}, we recall the notion of weak compactness for  countable p.m.p. equivalence relations.
In \cite[Definition 3.1]{OP07}, N. Ozawa and S. Popa defined the notion of weak compactness for p.m.p. actions $\Gamma\curvearrowright (X,\mu)$. In \cite[Proposition 3.4]{OP07} they established that if the action $\Gamma\curvearrowright (X,\mu)$ is weakly compact, then the action of the full group of $\mathcal R(\Gamma\curvearrowright X)$ on $X$ is also weakly compact. Thus, one can define weak compactness for countable p.m.p. equivalence relations $\mathcal R$ by insisting that the associated action of the full group $[\mathcal R]$ is weakly compact:

\begin{definition}\label{weakcomp}\cite{OP07}
A countable p.m.p. equivalence relation $\mathcal R$ on a standard probability space $(X,\mu)$ is said to be {\it weakly compact} if there exists a net of vectors $\eta_n\in L^2(X\times X,\mu^{\otimes_2})$ such that $\eta_n\geqslant 0$ and $\|\eta_n\|_2=1$, for all $n$, and the following conditions are satisfied:
\begin{enumerate}
\item $\lim\limits_{n}\|\eta_n-(u\otimes\bar{u})\eta_n\|_2=0$, for all $u\in\mathcal U(L^{\infty}(X))$.
\item $\lim\limits_{n}\|\eta_n-\eta_n\circ(\theta\times\theta)\|_2=0$, for all $\theta\in [\mathcal R]$.
\item $\lim\limits_{n}\langle (v\otimes 1)\eta_n,\eta_n\rangle=\lim\limits_{n}\langle (1\otimes v)\eta_n,\eta_n\rangle=\displaystyle{\int_{X}v\;\text{d}\mu}$, for all $v\in L^{\infty}(X)$.
\end{enumerate}
\end{definition}

Here, for $u,v\in L^{\infty}(X)$, the function $u\otimes\bar{v}\in L^{\infty}(X\times X)$ is given by $(u\otimes\bar{v})(x,y)=u(x)\overline{v(y)}$. Also, $\mu^{\otimes 2}$ denotes the product measure $\mu\otimes\mu$ on $X\times X$.
Note that the above conditions are precisely conditions (1),(2) and (3') from \cite[Definition 3.1]{OP07} for the action $[\mathcal R]\curvearrowright (X,\mu)$.

\subsection{Proof of Proposition \ref{propD}}
Denote $\mathcal R=\mathcal R(\Gamma\curvearrowright G)$ and let $A\subset G$ be a Borel set with $0<m_G(A)<+\infty$. Then $\mu=m_G(A)^{-1}m_G$ is a Haar measure of $G$ such that $\mu(A)=1$. 

Our goal is to show that $\mathcal R|A$ is weakly compact. Note first that we may assume that $A$ is an open subset of $G$ such that $\bar{A}$ is compact. This is because we can find an open set $B\subset G$  such that $\bar{B}$ is compact and $m_G(A)=m_G(B)$. Then since $\mathcal R$ is ergodic and preserves $m_G$, there exists $\theta\in [\mathcal R]$ such that $\theta(A)=B$, hence $\mathcal R|A\cong\mathcal R|B$.

Since $G$ is second countable, we can find a left invariant compatible metric $d$ on $G$. For $x\in G$ and $r>0$, we denote by $B_{r}(x)=\{y\in G|d(x,y)<r\}$ the open ball of radius $r$ centered at $x$. We also let $B_{r}=B_{r}(e)$. Notice that $\mu(B_{r}(x))=\mu(xB_{r})=\mu(B_{r})$.

For $\varepsilon>0$, we define $S_{\varepsilon}=\{(x,y)\in A\times A|d(x,y)<\varepsilon\}$.  Then $\mu^{\otimes_2}(S_{\varepsilon})>0$. Indeed, otherwise we would get that 
$\mu(B_{\varepsilon}(x)\cap A)=0$, for almost every $x\in A$. Since $A$ is open, this would imply that there exists $x\in A$ and $\varepsilon'>0$ such that $\mu(B_{\varepsilon'}(x))=0$, which contradicts the fact that $\mu$ is a Haar measure of $G$. Thus, we may further define $$\eta_{\varepsilon}:=\frac{1_{S_{\varepsilon}}}{\sqrt{\mu^{\otimes_2}(S_{\varepsilon})}}\in L^2(A\times A,\mu^{\otimes_2})$$

Then $\eta_{\varepsilon}\geqslant 0$ and $\|\eta_{\varepsilon}\|_2=1$, for all $\varepsilon>0$. 
We will show that the net $(\eta_{\varepsilon})$ verifies conditions (1)-(3) from Definition \ref{weakcomp}. 
Firstly, we verify condition (3).  To this end, for $\varepsilon>0$, we define $\xi_{\varepsilon}:A\rightarrow [0,\infty)$ by letting $\xi_{\varepsilon}(x)=\displaystyle{\int_{A}\eta_{\varepsilon}^2(x,y)\;\text{d}\mu(y)}$.

{\bf Claim 1}. We have that $\lim\limits_{\varepsilon\rightarrow 0}\|\xi_{\varepsilon}\|_{\infty}=1$ and $\lim\limits_{\varepsilon\rightarrow 0}\xi_{\varepsilon}(x)=1$, for almost every $x\in A$.

{\it Proof Claim 1.} We define a function $r:A\rightarrow (0,\infty)$ by letting $r(x)=\sup\{r>0|B_r(x)\subset A\}$. Since $A$ is open, $r$ is a well-defined continuos function and $B_{r(x)}(x)\subset A$, for all $x\in A$.

Let $n\geqslant 1$. Then there exists $\varepsilon_n>0$ such that $A_n:=\{x\in A|r(x)\geqslant\varepsilon_n\}$ satisfies $\mu(A_n)>1-2^{-n-1}$. 
Let $\varepsilon\in (0,\varepsilon_n]$. Then for all $x\in A_n$ we have that $B_{\varepsilon}(x)\subset A$ and hence $$\mu^{\otimes_2}(S_{\varepsilon})=\int_{A}\mu(B_{\varepsilon}(x)\cap A)\;\text{d}\mu(x)\geqslant\int_{A_n}\mu(B_{\varepsilon}(x))\;\text{d}\mu(x)=\mu(B_{\varepsilon})\mu(A_n)>(1-2^{-n-1})\mu(B_{\varepsilon}).$$

As a consequence, for every $x\in A$ we have that $$\xi_{\varepsilon}(x)=\mu^{\otimes_2}(S_{\varepsilon})^{-1}\mu(B_{\varepsilon}(x)\cap A)\leqslant \mu^{\otimes_2}(S_{\varepsilon})^{-1}\mu(B_{\varepsilon})<(1-2^{-n-1})^{-1}<1+2^{-n}.$$

This shows that \begin{equation}\label{una}\|\xi_{\varepsilon}\|_{\infty}<1+2^{-n},\;\;\text{for all}\;\;\varepsilon\in (0,\varepsilon_n].\end{equation}

On the other hand,  we have that  $\mu^{\otimes_2}(S_{\varepsilon})=\displaystyle{\int_{A}\mu(B_{\varepsilon}\cap A)\;\text{d}\mu(x)\leqslant\mu(B_{\varepsilon})}$. Thus, if $x\in A_n$, then $\xi_{\varepsilon}(x)=\mu^{\otimes_2}(S_{\varepsilon})^{-1}\mu(B_{\varepsilon}(x)\cap A)=\mu^{\otimes_2}(S_{\varepsilon})^{-1}\mu(B_{\varepsilon})\geqslant 1$. Hence, we get that \begin{equation}\label{doua}\xi_{\varepsilon}(x)\geqslant 1,\;\;\text{for all}\;\;\varepsilon\in (0,\varepsilon_n]\;\;\text{and every}\;\;x\in A_n.\end{equation}

Since $\mu(A_n)\geqslant 1-2^{-n-1}$, for all $n$, it is easy to see that \ref{una} and \ref{doua} together imply the claim.
\hfill$\square$

If $v\in L^{\infty}(A)$, then $\langle(v\otimes 1)\eta_{\varepsilon},\eta_{\varepsilon}\rangle=\displaystyle{\int_{A}v(x)\xi_{\varepsilon}(x)\;\text{d}\mu(x)}$. By combining Claim 1 and the Lebesgue dominated convergence theorem, we get that $\lim\limits_{\varepsilon\rightarrow 0}\langle(v\otimes 1)\eta_{\varepsilon},\eta_{\varepsilon}\rangle=\displaystyle{\int_{A}v(x)\text{d}\mu(x)}$. Since $\eta_{\varepsilon}$ is symmetric, we have that $\langle (1\otimes v)\eta_{\varepsilon},\eta_{\varepsilon}\rangle=\langle (v\otimes 1)\eta_{\varepsilon},\eta_{\varepsilon}\rangle$ and altogether condition (3) follows.

Towards showing that the net $(\eta_{\varepsilon})$ verifies conditions (1) and (2), we first establish the following:

{\bf Claim 2.} We have that $\lim\limits_{\varepsilon\rightarrow 0}\|(v\otimes 1)\eta_{\varepsilon}-(1\otimes v)\eta_{\varepsilon}\|_2=0$, for all $v\in L^{\infty}(A)$.

{\it Proof of Claim 2.}  Let $\mathcal C$ be the set of functions $v\in L^{\infty}(A)$ for which there exists a continuous function $\tilde v:\bar{A}\rightarrow\mathbb C$ such that $v={\tilde{v}}_{|A}$. Then $\mathcal C$ is $\|.\|_2$-dense in $L^{\infty}(A)$. On the other hand, condition (3) implies that $\lim\limits_{\varepsilon\rightarrow 0}\|(v\otimes 1)\eta_{\varepsilon}\|_2=\lim\limits_{\varepsilon\rightarrow 0}\|(1\otimes v)\eta_{\varepsilon}\|_2=\|v\|_2$, for every $v\in L^{\infty}(A)$. Thus, in order to prove the claim, it suffices to show that it holds for every $v\in\mathcal C$.

Let $v\in\mathcal C$ and $\tilde v$ be a continuous extension of $v$ to $\bar{A}$. Let $\delta>0$. Since $\tilde{v}$ is continuous and $\bar{A}$ is compact,  $\tilde v$ is uniformly continuous. It follows that we can find $\varepsilon_0>0$ such that $|v(x)-v(y)|<\delta$, for all $x,y\in A$ such that $d(x,y)\leqslant\varepsilon_0$. From this we deduce that for all $\varepsilon\in (0,\varepsilon_0]$ we have

$$\|(v\otimes 1)\eta_{\varepsilon}-(1\otimes v)\eta_{\varepsilon}\|_2^2=\mu^{\otimes_2}(S_{\varepsilon})^{-1}\int_{S_{\varepsilon}}|v(x)-v(y)|^2\;\text{d}\mu^{\otimes_2}(x,y)<\delta^2.$$

Since $\delta>0$ is arbitrary, we conclude that  $\lim\limits_{\varepsilon\rightarrow 0}\|(v\otimes 1)\eta_{\varepsilon}-(1\otimes v)\eta_{\varepsilon}\|_2=0$.
\hfill$\square$

Now, if $u\in\mathcal U(L^{\infty}(A))$, then by Claim 2 we get $\|\eta_{\varepsilon}-(u\otimes\bar{u})\eta_{\varepsilon}\|_2=\|(1\otimes u)\eta_{\varepsilon}-(u\otimes 1)\eta_{\varepsilon}\|_2\rightarrow 0$, as $\varepsilon\rightarrow 0$, which proves condition (1).

 Finally, in order to prove condition (2), let $\theta\in [\mathcal R|A]$.
Then we can find a Borel map $\phi:A\rightarrow\Gamma$ such that $\theta(x)=\phi(x)x$, for almost every $x\in A$.
Notice that if $x,y\in A$ and $\phi(x)=\phi(y)$, then we have $d(\theta(x),\theta(y))=d(x,y)$. 
Using this observation we get that if $\varepsilon>0$ then \begin{equation}\label{eta1}\|\eta_{\varepsilon}-\eta_{\varepsilon}\circ(\theta\times\theta)\|_2^2=\mu^{\otimes_2}(S_{\varepsilon})^{-1}\int_{X\times X}|1_{S_{\varepsilon}}(x,y)-1_{S_{\varepsilon}}(\theta(x),\theta(y))|^2\;\text{d}\mu^{\otimes_2}(x,y)\leqslant\end{equation}
$$\mu^{\otimes_2}(S_{\varepsilon})^{-1}\mu^{\otimes_2}(\{(x,y)\in S_{\varepsilon}|\phi(x)\not=\phi(y)\}).$$

For $g\in\Gamma$, denote $A_g=\{x\in A|\phi(x)=g\}$. Then we have

\begin{equation}\label{eta2}\mu^{\otimes_2}(\{(x,y)\in S_{\varepsilon}|\phi(x)\not=\phi(y)\})\leqslant\sum_{g\in\Gamma}\mu^{\otimes_2}(\{(x,y)\in S_{\varepsilon}|x\in A_g,y\notin A_g\})\leqslant \end{equation} $$\mu^{\otimes_2}(S_{\varepsilon})\sum_{g\in\Gamma}\|(1_{A_g}\otimes 1)\eta_{\varepsilon}-(1\otimes 1_{A_g})\eta_{\varepsilon}\|_2^2$$

The combination of \ref{eta1} and \ref{eta2} further implies that \begin{equation}\label{eta3}\|\eta_{\varepsilon}-\eta_{\varepsilon}\circ(\theta\times\theta)\|_2^2\leqslant\sum_{g\in\Gamma}\|(1_{A_g}\otimes 1)\eta_{\varepsilon}-(1\otimes 1_{A_g})\eta_{\varepsilon}\|_2^2,\;\;\text{for all}\;\;\varepsilon>0.\end{equation}

Since $\lim\limits_{\varepsilon\rightarrow 0}\|\xi_{\varepsilon}\|_{\infty}=1$ by Claim 1, we can find $\varepsilon_0$ such that $\|\xi_{\varepsilon}\|_{\infty}\leqslant 2$, for all $\varepsilon\in (0,\varepsilon_0]$. Let $v\in L^{\infty}(A)$ and $\varepsilon\in (0,\varepsilon_0]$. Then $\|(v\otimes 1)\eta_{\varepsilon}\|_2^2=\displaystyle{\int_{A\times A}|v|^2\xi_{\varepsilon}\;\text{d}\mu}\leqslant 2\|v\|_2^2$.  Since $\eta_{\varepsilon}$ is symmetric, we get that $\|(1\otimes v)\eta_{\varepsilon}\|_2^2=\|(v\otimes 1)\eta_{\varepsilon}\|_2^2\leqslant 2\|v\|_2^2$. As a consequence, we have that \begin{equation}\label{eta4}\|(v\otimes 1)\eta_{\varepsilon}-(1\otimes v)\eta_{\varepsilon}\|_2^2\leqslant 8\|v\|_2^2,\;\;\text{for all}\;\; v\in L^{\infty}(A)\;\;\text{and every}\;\;\varepsilon\in (0,\varepsilon_0].\end{equation}

Let $\delta>0$. Then we can find a finite set $F\subset \Gamma$ such that $\sum_{g\in\Gamma\setminus F}\mu(A_g)\leqslant\displaystyle{\frac{\delta^2}{16}}$.
In combination with equation \ref{eta4} we further get \begin{equation}\label{eta5}\sum_{g\in\Gamma\setminus F}\|(1_{A_g}\otimes 1)\eta_{\varepsilon}-(1\otimes 1_{A_g})\eta_{\varepsilon}\|_2^2\leqslant 8\sum_{g\in\Gamma\setminus F}\|1_{A_g}\|_2^2\leqslant\frac{\delta^2}{2},\;\;\text{for all}\;\;\varepsilon\in (0,\varepsilon_0]. 
 \end{equation}

Next, by using condition (1), we can find $\varepsilon_1\in (0,\varepsilon_0]$ such that \begin{equation}\label{eta6}\sum_{g\in F}\|(1_{A_g}\otimes 1)\eta_{\varepsilon}-(1\otimes 1_{A_g})\eta_{\varepsilon}\|_2^2\leqslant\frac{\delta^2}{2},\;\;\text{for all}\;\;\varepsilon\in (0,\varepsilon_1].\end{equation}

Finally, the combination of equations \ref{eta3}, \ref{eta5} and \ref{eta6} gives that $\|\eta_{\varepsilon}-\eta_{\varepsilon}\circ(\theta\times\theta)\|_2\leqslant\delta$, for all $\varepsilon\in (0,\varepsilon_1]$. Since $\delta>0$ is arbitrary, we conclude that the net $(\eta_{\varepsilon})$ satisfies condition (2). \hfill$\blacksquare$

\section{Proof of Theorem \ref{thmE}}

This section is devoted to the proof of Theorem \ref{thmE}. In fact, we prove the following more general and more precise version of Theorem \ref{thmE}:

\begin{theorem}\label{genthmE} Let $G$ be a connected l.c.s.c. with trivial center. Assume that there is a l.c.s.c. simply connected group $\widetilde G$ together with a continuous onto homomorphism $p:\widetilde G\rightarrow G$ such that $\ker(p)<\widetilde G$ is discrete. Let $\Sigma<G$ be a discrete subgroup and $\Gamma<G$  a countable dense subgroup. 
Assume that the translation action $\Gamma\curvearrowright (G,m_G)$ is strongly ergodic.
Let $\bar{H}$ be a semisimple real algebraic group and denote $H=\bar{H}/Z(\bar{H})$.  Let $\Delta<H$ be a discrete subgroup and $\Lambda<H$ a countable subgroup.

Let $A\subset G/\Sigma$ and $B\subset H/\Delta$ be non-negligible measurable sets, and $\theta:A\rightarrow B$ be a nonsingular isomorphism such that  $\theta(\Gamma x\cap A)=\Lambda\theta(x)\cap B$, for almost every $x\in A$. 

Then we can find a Borel map $\phi:G/\Sigma\rightarrow\Lambda$, a topological isomorphism $\delta:G\rightarrow H$ and $h\in H$ such that $\delta(\Gamma)=\Lambda$, $\delta(\Sigma)=h\Delta h^{-1}$ and $\theta(x)=\phi(x)\delta(x)h\Delta$, for almost every $x\in A$.
\end{theorem}

Before proving Theorem \ref{genthmE}, we establish the following elementary result:

\begin{lemma}\label{liber} Let $G$ be a connected l.c.s.c. group, $\Gamma<G$  a countable subgroup and $\Lambda<G$ a discrete subgroup. 
If $\Gamma\cap\Lambda$ contains no non-trivial central element of $G$, then the action $\Gamma\curvearrowright G/\Lambda$ is free.
\end{lemma}

{\it Proof.} Let $g\in\Gamma$ such that the set $\{x\in G|gx\Lambda=x\Lambda\}$ has positive measure. Since $\Lambda$ is countable, we can find $h\in\Lambda$ such that the set $A=\{x\in G|gx=xh\}$ has positive measure. Since $A$ has positive measure, $AA^{-1}$ contains a neighborhood of $e\in G$. Since $g$ commutes with $AA^{-1}$ and $G$ is connected, we get that $g$ belongs to the center $Z(G)$ of $G$. Since $gx=xh$, for some $x\in G$, it follows that $g=h\in\Gamma\cap\Lambda\cap Z(G)$. Therefore, $g=e$, as claimed.
\hfill$\blacksquare$

\subsection{Proof of Theorem \ref{genthmE}} Since the action $\Gamma\curvearrowright G/\Sigma$ is ergodic, we may extend $\theta$ to a countable-to-1 map $\theta:G/\Sigma\rightarrow H/\Delta$ such that $\theta(\Gamma x)\subset\Lambda\theta(x)$, for almost every $x\in G/\Sigma$.

Denote by $\pi:G\rightarrow G/\Sigma$  the quotient. Then $\theta(\pi(gx))\in\Lambda\theta(\pi(x))$, for all $g\in\Gamma$ and almost every $x\in G$. Since $H$ has trivial center, by Lemma \ref{liber}, the action $\Lambda\curvearrowright H/\Delta$ is free. Thus,  we can define a cocycle $W:\Gamma\times G\rightarrow\Lambda$ through the formula $\theta(\pi(gx))=W(g,x)\theta(\pi(x))$.  The freeness of the $\Lambda$-action also implies that $W$ factors through the map $\Gamma\times G\rightarrow\Gamma\times G/\Sigma$.

Let  $r:H/\Delta\rightarrow H$ be a Borel map such that $r(x)\Delta=x$, for all $x\in H/\Delta$.
Further, we define $\Theta:G\rightarrow H$ by letting $\Theta(x)=r(\theta(\pi(x)))$. Then for all $g\in\Gamma$ and almost every $x\in G/\Sigma$ we have $$\Theta(gx)\Delta=r(\theta(\pi(gx)))\Delta=\theta(\pi(gx))=W(g,x)\theta(\pi(x))=W(g,x)\Theta(x)\Delta.$$

We can therefore find a Borel map $v:\Gamma\times G\rightarrow\Delta$ such that \begin{equation}\label{bigteta}\Theta(gx)=W(g,x)\Theta(x)v(g,x),\;\;\text{for all}\;\;g\in\Gamma\;\;\text{and almost every}\;\;x\in G/\Sigma.\end{equation}

The core of the proof is divided in two parts. 

\vskip 0.05in

{\bf The first part of the proof}. In this part we analyze the cocycle $W$ and show that there exists a Borel map $\alpha:G\rightarrow H$ such that  $W(g,x)=\alpha(gx)\alpha(x)^{-1}$, for all $g\in\Gamma$ and almost every $x\in G$.

\vskip 0.05in

To this end, for a subset $S\subset H$, we denote by $C(S)$ its centralizer in $H$. We denote by $\mathcal A$ the set of pairs $(h,k)\in H\times H$ such that $C(\{h,k\})=\{e\}$.  Since $H$ is a connected, semisimple real Lie group,  \cite[Theorem 4]{Wi02} implies that the subgroup generated by $h$ and $k$ is dense in $H$, for almost every $(h,k)\in H\times H$. Since $H$ has trivial center, for any such pair $(h,k)$ we have that $C(\{h,k\})=\{e\}$.  Altogether, we deduce that $\mathcal A$ is conull in $H\times H$.

Next, for $t\in G$, we define a Borel map $\rho_t:G\rightarrow H$ by letting $\rho_t(x)=\Theta(xt)\Theta(x)^{-1}$.

{\bf Claim 1.} For almost every $(s,t)\in G\times G$ we have that $(\rho_s(x),\rho_t(x))\in\mathcal A$, for almost every $x\in G$.

{\it Proof of Claim 1.} Let $x\in G$ and define $\mathcal A_x=\{(h\Theta(x),k\Theta(x))|(h,k)\in\mathcal A\}$. Since $\mathcal A$ is conull in $H\times H$, $\mathcal A_x$ is also conull in $H\times H$. Since $\theta:G/\Sigma\rightarrow H/\Delta$ is nonsingular, $\Theta:G\rightarrow H$ is also nonsingular: if $B\subset H$ satisfies $m_H(B)=0$, then $m_G(\Theta^{-1}(B))=0$. By combining these two facts we conclude that the set $\{(a,b)\in G\times G|(\Theta(a),\Theta(b))\in\mathcal A_x\}$ is conull. 
Further,  we get that $\{(s,t)\in G\times G|(\rho_s(x),\rho_t(x))\in\mathcal A\}=\{(s,t)\in G\times G|(\Theta(xs),\Theta(xt))\in\mathcal A_x\}$ is conull in $G\times G$.
Since this holds for all $x\in G$, we derive that the set $\{(s,t,x)\in G\times G\times G|(\rho_s(x),\rho_t(x))\in\mathcal A\}$ is conull in $G\times G\times G$. The claim now follows by applying Fubini's theorem. \hfill$\square$

Let $\mu$ be a Borel probability measure on $G$ which is equivalent to $m_G$.

{\bf Claim 2.} There exists a sequence of  Borel maps $\alpha_n:G\rightarrow H$ such that for all $g\in\Gamma$ we have $\lim\limits_{n\rightarrow\infty}\mu(\{x\in G|\alpha_n(gx)=W(g,x)\alpha_n(x)\})=1$.

{\it Proof of Claim 2.} Let $\varepsilon>0$ and $F\subset\Gamma$ finite be arbitrary. In order to prove the claim, it suffices to find a Borel map $\alpha:G\rightarrow H$ such that $\mu(\{x\in G|\alpha(gx)=W(g,x)\alpha(x)\})\geqslant 1-\varepsilon$, for all $g\in F$.

Since $\mu$ is $\Gamma$-quasi-invariant,  we can find $\varepsilon_0\in (0,\displaystyle{\frac{\varepsilon}{3}})$ such that if a Borel set $B\subset G$ satisfies $\mu(B)\geqslant 1-\varepsilon_0$, then $\mu(g^{-1}B)\geqslant 1-\displaystyle{\frac{\varepsilon}{3}}$, for all $g\in F$. Since the action $\Gamma\curvearrowright G$ is strongly ergodic, Lemma \ref{strong} provides $\delta>0$ and a finite set $S\subset\Gamma$ such that the following holds: if $\rho:G\rightarrow Y$ is a Borel map into a standard Borel space $Y$ which satisfies $\mu(\{x\in G|\rho(gx)=\rho(x)\})\geqslant 1-\delta$, for all $g\in S$, then we can find $y\in Y$ such that $\mu(\{x\in G|\rho(x)=y\})\geqslant 1-\varepsilon_0$.

Now, let $t\in G$. Using equation \ref{bigteta}, for all $g\in\Gamma$ and almost every $x\in G$ we have $$\rho_t(gx)=\Theta(gxt)\Theta(gx)^{-1}=W(g,xt)\Theta(xt)\;(v(g,xt)v(g,x)^{-1})\;\Theta(x)^{-1}W(g,x)^{-1}.$$ 

This implies that if $v(g,xt)=v(g,x)$ and $W(g,xt)=W(g,x)$, then $\rho_t(gx)=W(g,x)\rho_t(x)W(g,x)^{-1}$. Since $v$ and $W$ take values into countable groups, we have that $\mu(\{x\in G|v(g,xt)=v(g,x)\})\rightarrow 1$ and $\mu(\{x\in G|W(g,xt)=W(g,x)\})\rightarrow 1$, as $t\rightarrow e$. Altogether, we conclude that for all $g\in\Gamma$ we have $\mu(\{x\in G|\rho_t(gx)=W(g,x)\rho_t(x)W(g,x)^{-1}\})\rightarrow 1$, as $t\rightarrow e$.

This last fact in combination with Claim 1 implies that we can find $s,t\in G$ such that the following five conditions are satisfied:
\begin{enumerate}
\item $\mu(\{x\in G|\rho_s(gx)=W(g,x)\rho_s(x)W(g,x)^{-1}\})\geqslant 1-{\frac{\varepsilon}{6}}$, for all $g\in F$. 
\item $\mu(\{x\in G|\rho_s(gx)=W(g,x)\rho_s(x)W(g,x)^{-1}\})\geqslant 1-{\frac{\delta}{2}}$, for all $g\in S$. 
\item $\mu(\{x\in G|\rho_t(gx)=W(g,x)\rho_t(x)W(g,x)^{-1}\})\geqslant 1-{\frac{\varepsilon}{6}}$, for all $g\in F$. 
\item $\mu(\{x\in G|\rho_t(gx)=W(g,x)\rho_t(x)W(g,x)^{-1}\})\geqslant 1-{\frac{\delta}{2}}$, for all $g\in S$. 
\item $(\rho_s(x),\rho_t(x))\in\mathcal A$, for almost every $x\in G$.
\end{enumerate}

Next, consider the ``diagonal conjugation" action $\bar{H}\curvearrowright \bar{H}\times \bar{H}$ given by $h\cdot (k,l)=(hkh^{-1},hlh^{-1})$, for all $h\in \bar{H}$ and $(k,l)\in \bar{H}\times \bar{H}$. Since $\bar{H}$ is a real algebraic group, this action is smooth by Theorem \ref{smooth}.  Therefore, the quotient space $\bar{Y}$ is a standard Borel space. Analogously, consider the diagonal conjugation action $H\curvearrowright H\times H$ and denote by $Y$ the quotient space. Since the left multiplication action of  $Z(\bar{H})\times Z(\bar{H})$ on $\bar{H}\times\bar{H}$ commutes with the diagonal conjugation action of $\bar{H}$, we have a well-defined left multiplication action $Z(\bar{H})\times Z(\bar{H})\curvearrowright\bar{Y}$.  Then $Y\equiv \bar{Y}/(Z(\bar{H})\times Z(\bar{H}))$ and we conclude that $Y$ is a standard Borel space.
We denote by  $q:H\times H\rightarrow Y$ the quotient map.

By combining  (2) and (4) we get that $\mu(\{x\in G|q(\rho_s(gx),\rho_t(gx))=q(\rho_s(x),\rho_t(x))\})\geqslant 1-\delta$, for all $g\in S$. Since $Y$ is a standard Borel space we can apply the above consequence of strong ergodicity to the Borel map $G\ni x\rightarrow q(\rho_s(x),\rho_t(x))\in Y$. Thus, we get that there exists $y\in H\times H$ such that $\mu(\{x\in G|(\rho_s(x),\rho_t(x))\in H\cdot y\})\geqslant 1-\varepsilon_0$.

This fact and condition (5) imply  that we can find $x_0\in G$ such that $(\rho_s(x_0),\rho_t(x_0))\in \mathcal A\cap H\cdot y$. 
Since $(\rho_s(x_0),\rho_t(x_0))\in\mathcal A$, the stabilizer of $(\rho_s(x_0),\rho_t(x_0))$ in $H$ is trivial. Thus, the stabilizer of $y$ in $H$ is trivial. Hence, we can identify the orbit $H\cdot y$ with $H$.  

We define a Borel map $\alpha:G\rightarrow H$ by letting $$\alpha(x)=\begin{cases} (\rho_s(x),\rho_t(x)),\;\;\text{if}\;\;(\rho_s(x),\rho_t(x))\in H\cdot y \\ e,\;\;\text{otherwise.}\end{cases}$$

Finally, fix $g\in F$.  Let $C=\{x\in G|(\rho_s(gx),\rho_t(gx))=W(g,x)\cdot (\rho_s(x),\rho_t(x))\}$.  The combination of conditions (1) and (3)  yields that $\mu(C)\geqslant 1-\displaystyle{\frac{\varepsilon}{3}}$.  Also, if $D=\{x\in G|\alpha(x)=(\rho_s(x),\rho_t(x))\}$, then $\mu(D)\geqslant 1-\varepsilon_0>1-\displaystyle{\frac{\varepsilon}{3}}$. Since $g\in F$, we get that $\mu(g^{-1}D)\geqslant 1-\displaystyle{\frac{\varepsilon}{3}}$, and hence $\mu(C\cap D\cap g^{-1}D)\geqslant 1-\varepsilon$. If $x\in C\cap D\cap g^{-1}D$,  then we clearly have that $\alpha(gx)=W(g,x)\alpha(x)$, and the claim is proven.
\hfill$\square$

\vskip 0.05in
Since the action $\Gamma\curvearrowright G$ is strongly ergodic, combining Claim 2 and Lemma \ref{untwist} implies that we can find a Borel map $\alpha:G\rightarrow H$ such that $W(g,x)=\alpha(gx)\alpha(x)^{-1}$, for all $g\in\Gamma$ and almost every $x\in G$.  
This finishes the first part of the proof.

\vskip 0.05in
{\bf The second part of the proof}.
 Next, we use the conclusion of the first part to derive the theorem. This part relies on the following claim:
 
 {\bf Claim 3.} There exist a Borel map $\phi:G\rightarrow\Lambda$, a topological isomorphism $\delta:G\rightarrow H$, and $y\in H/\Delta$ such that $\delta(\Gamma)\subset\Lambda$ and $\theta(\pi(x))=\phi(x)\delta(x)y$, for almost every $x\in G$. 

{\it Proof of Claim 3}.  Denote $\widetilde\Gamma=p^{-1}(\Gamma)$. Since $H$ is a connected real Lie group, we can find a simply connected l.c.s.c. group $\widetilde H$  together with a continuous onto homomorphisms $q:\widetilde H\rightarrow H$ such that  $\ker(q)<\widetilde H$ is a discrete subgroup.

Let $\tilde W:\widetilde\Gamma\times\widetilde G\rightarrow\Lambda$ and $\tilde\alpha:\widetilde G\rightarrow H$ be given by $\tilde W(g,x)=W(p(g),p(x))$ and $\tilde\alpha(x)=\alpha(p(x))$. Then $\tilde W(g,x)=\tilde\alpha(gx)\tilde\alpha(x)^{-1}$, for all $g\in\widetilde\Gamma$ and almost every $x\in\widetilde G$. Since $\widetilde G$ is simply connected, by Corollary \ref{cocrig2}, we can find a Borel map $\phi:\widetilde G\rightarrow\Lambda$ and a homomorphism $\rho:\widetilde\Gamma\rightarrow\Lambda$ such that \begin{equation}\label{tildeW1}\tilde W(g,x)=\phi(gx)\rho(g)\phi(x)^{-1},\;\;\text{for all $g\in\widetilde\Gamma$ and almost every $x\in\widetilde G$}. \end{equation}

 Define $\beta:\widetilde G\rightarrow H$ by letting $\beta(x)=\phi(x)^{-1}\tilde\alpha(x)$. Then we get that $\beta(gx)=\rho(g)\beta(x)$, for all $g\in\widetilde\Gamma$ and almost every $x\in\widetilde G$. By Lemma \ref{ext}, $\rho$ extends to a continuous homomorphism $\rho:\widetilde G\rightarrow H$.

Our next goal is to show that $\rho$ is onto.
Define $\tilde\theta:\widetilde G\rightarrow H/\Delta$ by $\tilde\theta(x)=\theta(\pi(p(x)))$. Then $\tilde\theta(gx)=\tilde W(g,x)\tilde\theta(x)$, for all $g\in\widetilde\Gamma$ and almost every $x\in\widetilde G$.
Thus, if 
 $\hat{\theta}:\widetilde G\rightarrow H/\Delta$ is given by  $\hat{\theta}(x)=\phi(x)^{-1}\tilde\theta(x)$, then \ref{tildeW1} implies that \begin{equation}\label{hatheta}\hat{\theta}(gx)=\rho(g)\hat{\theta}(x),\;\;\text{for all $g\in\widetilde\Gamma$ and almost every $x\in\widetilde G$}. \end{equation}

Since  $\widetilde\Gamma<\widetilde G$ is dense and $\rho:\widetilde G\rightarrow H$ is continuos, equation \ref{hatheta} holds for every $g\in\widetilde G$. By Fubini's theorem, we can find $x_0\in H/\Delta$ such that $\hat{\theta}(gx_0)=\rho(g)\hat{\theta}(x_0)$, for almost every $g\in\widetilde G$. Thus, if $y=\rho(x_0)^{-1}\hat{\theta}(x_0)$, then $\hat{\theta}(g)=\rho(g)y$, for almost every $g\in\widetilde G$. From this we deduce that \begin{equation}\label{teta0}\tilde\theta(x)=\phi(x)\rho(x)y,\;\;\text{for almost every}\;\;x\in\widetilde G.\end{equation}

Since $\theta:A\rightarrow B$ is a nonsingular isomorphism, we get that $m_{H/\Delta}(\tilde\theta(\widetilde G))=m_{H/\Delta}(\theta(G/\Delta))>0$. Since $\Lambda$ and $\Delta$ are countable, by using \ref{teta0}, it follows that $m_H(\rho(\widetilde G))>0$.
This implies that $\rho(\widetilde G)<H$ is an open subgroup. Since $H$ is connected, we conclude that $\rho(\widetilde G)=H$, hence $\rho$ is onto.

Next, let us show that $\ker(\rho)$ is a discrete subgroup of $\widetilde G$.
To this end, let $g_n\in\ker(\rho)$ be a sequence which converges to the identity. Since $\phi(\widetilde G)\subset\Lambda$ and $\Lambda$ is countable, there exists $N\geqslant 1$ such that $m_{\widetilde G}(\{x\in\widetilde G|\phi(g_nx)=\phi(x)\;\text{and}\; \pi(p(x)), \pi(p(g_nx))\in A\})>0$, for all $n\geqslant N$. By using \ref{teta0} we further get that $m_{\widetilde G}(\{x\in\widetilde G|\tilde\theta(g_nx)=\tilde\theta(x)\;\text{and}\;\pi(p(x)), \pi(p(g_nx))\in A\})>0$, for all $n\geqslant N$.  Since $\theta$ is 1-1 on $A$, we get that  $m_{\widetilde G}(\{x\in\widetilde G|p(g_n)\pi(p(x))=\pi(p(g_nx))=\pi(p(x))\})>0$, for all $n\geqslant N$. Thus, $p(g_n)$ does not act freely on $G/\Sigma$, hence by Lemma \ref{liber} we get that $p(g_n)=e$, for all $n\geqslant N$. Since $\ker(p)$ is discrete, we deduce that $g_n=e$, for $n$ large enough.

Altogether, we have shown that $\rho:\widetilde G\rightarrow H$ is an onto continuous homomorphism with discrete kernel,  hence $\widetilde G$ is a covering group of $H$. The uniqueness of covering groups implies that we can find a topological isomorphism $\tau:\widetilde G\rightarrow\widetilde H$ such that $\rho=q\circ\tau$. Since $G$ and  $H$ have trivial center, we get that $\ker(p)=Z(\widetilde G)$ and $\ker(q)=Z(\widetilde H)$ and hence $\ker(\rho)=\tau^{-1}(\ker(p))=\tau^{-1}(Z(\widetilde H))=Z(\widetilde G)$.
Therefore, $\rho$ descends to a topological isomorphism $\delta:G\rightarrow H$. Since $\rho(\widetilde\Gamma)\subset\Lambda$, we have $\delta(\Gamma)\subset\Lambda$.

Finally, let $\sigma\in\ker(p)$. Then $\rho(\sigma)=e$ and $\tilde\theta(x\sigma)=\tilde\theta(x)$, for all $x\in\widetilde G$. By using \ref{teta0} we deduce that $\phi(x\sigma)\rho(x)y=\phi(x)\rho(x)y$, for almost every $x\in\widetilde G$. 
Since $\phi(x),\phi(x\sigma)\in\Lambda$ and $\Lambda$ is countable, Lemma \ref{liber} implies that $\phi(x\sigma)=\phi(x)$, for almost every $x\in\widetilde G$. Therefore, $\phi:\widetilde G\rightarrow\Lambda$ descends to a map $\phi:G\rightarrow\Lambda$. Together with equation \ref{teta0}, this proves the claim.
 \hfill$\square$

{\bf Claim 4.} $\phi:G\rightarrow\Lambda$ factors through the quotient $\pi:G\rightarrow G/\Sigma$.

{\it Proof of Claim 4.} Since $\phi:\widetilde G\rightarrow\Lambda$ factors through the map $p:\widetilde G\rightarrow G$ by the previous claim,  equation \ref{tildeW1} rewrites as $W(g,x)=\phi(gx)\delta(g)\phi(x)^{-1}$, for all $g\in\Gamma$ and almost every $x\in G$.

 Let $\sigma\in\Sigma$. 
 Since $W$ factors through the map $\Gamma\times G\rightarrow\Gamma\times G/\Sigma$, we get that $$\phi(gx\sigma)\delta(g)\phi(x\sigma)^{-1}=W(g,x\sigma)=W(g,x)=\phi(gx)\delta(g)\phi(x)^{-1}.$$ Hence, if $\Phi:G\rightarrow\Lambda$ is given by $\Phi(x)=\phi(x)^{-1}\phi(x\sigma)$, then $\Phi(gx)=\delta(g)\Phi(x)\delta(g)^{-1}$, for all $g\in\Gamma $ and almost every $x\in G$. Since $\delta:G\rightarrow H$ is continuous, it follows that  $\Phi(gx)=\delta(g)\Phi(x)\delta(g)^{-1}$, for all $g\in G$ and almost every $x\in G$. By Fubini's theorem, we can find $x_0\in G$ such that  $\Phi(gx_0)=\delta(g)\Phi(x_0)\delta(g)^{-1}$, for almost every $g\in G$. Thus, if we let $k=\delta(x_0)^{-1}\Phi(x_0)\delta(x_0)$, then $\Phi(g)=\delta(g)k\delta(g)^{-1}$, for almost every $g\in G$.

Since $\Lambda$ is countable, we can find $l\in\Lambda$ such that $C=\{g\in G|\Phi(g)=l\;\text{and}\;\Phi(g)=\delta(g)k\delta(g)^{-1}\}$ satisfies $m_G(C)>0$. Note that $k$ commutes with $\delta(C^{-1}C)$. Since $m_G(C)>0$ and $\delta$ is onto, we get that $m_H(\delta(C))>0$ and hence $m_H(\delta(C^{-1}C))>0$. Combining the last two facts we derive that the centralizer of $k$ in $H$ has positive measure,  hence is an open subgroup of $H$. Since $H$ is connected and has trivial center, we get that $k=e$.
This implies that $\Phi(x)=e$, for almost every $x\in G$. Thus, $\phi(x\sigma)=\phi(x)$, for almost every $x\in G$. Since $\sigma\in\Sigma$ is arbitrary, this proves the claim. \hfill$\square$

Let $h\in H$ such that $y=h\Delta$.
We end the proof of Theorem \ref{thmE} by showing the following:

{\bf Claim 5.} $\delta(\Gamma)=\Lambda$ and $\delta(\Sigma)=h\Delta h^{-1}$.

{\it Proof of Claim 5.} 
By Claim 3, we have that $\delta(\Gamma)\subset\Lambda$. To show the reverse inclusion, let $g\in G$ such that $\delta(g)\in\Lambda$. 
By Claim 3 we also have that $$\theta(\pi(gx))=\phi(gx)\delta(g)\delta(x)y\in\Lambda\delta(x)y=\Lambda\phi(x)\delta(x)y=\Lambda\theta(\pi(x)),\;\; \text{for almost every}\;\;x\in G.$$ Thus, $\theta(gx)\in\Lambda\theta(x)$, for almost every $x\in G/\Sigma$. Since $\theta:A\rightarrow B$ is a stable orbit equivalence, we deduce that $gx\in\Gamma x$, for almost every $x\in G$, hence $g\in\Gamma$. This proves that $\delta(\Gamma)=\Lambda$.

Next, we show that $\delta(\Sigma)\subset h\Delta h^{-1}$. For this, let $\sigma\in\Sigma$. By Claim 4 we have that $\phi(x\sigma)=\phi(x)$, for almost every $x\in G$. By using Claim 3 we get that $\theta(\pi(x\sigma))=\phi(x\sigma)\delta(x\sigma)h\Delta=\phi(x)\delta(x)\delta(\sigma)h\Delta$, for almost every $x\in G$. Since $\theta(\pi(x\sigma))=\theta(\pi(x))=\phi(x)\delta(x)h\Delta$, for almost every $x\in G$, we deduce that $\delta(\sigma)h\Delta=h\Delta$. Thus, $\delta(\sigma)\in h\Delta h^{-1}$, as desired.

Finally, we prove that $\delta^{-1}(h\Delta h^{-1})\subset\Sigma$. To see this, let $g\in G$ such that $\delta(g)\in h\Delta h^{-1}$.  By using Claim 3 we get that $\theta(\pi(xg))=\phi(xg)\delta(xg)h\Delta=\phi(xg)\delta(x)h\Delta=\phi(xg)\phi(x)^{-1}\theta(\pi(x))\in\Lambda\theta(\pi(x))$, for almost every $x\in G$. Equivalently, we have that $\theta(\pi(xg))\in\Lambda\theta(\pi(x))$, for almost every $x\in G$. Since $\theta:A\rightarrow B$ is an orbit equivalence, we conclude that $\pi(xg)\in\Gamma\pi(x)$, or equivalently, that $xg\in\Gamma x\Sigma$, for almost every $x\in G$. Since $\Gamma$ and $\Sigma$ are countable, we can find $\gamma\in\Gamma$ and $\sigma\in\Sigma$ such that $C=\{x\in G|xg=\gamma x\sigma\}$ satisfies $m_G(C)>0$. It is clear that $g\sigma^{-1}$ commutes with $C^{-1}C$. Thus, the centralizer of $g\sigma^{-1}$ in $G$ has positive measure,  hence is an open subgroup of $G$. Since $G$ is connected and has trivial center we conclude that $g\sigma^{-1}=e$. Hence, $g=\sigma\in\Sigma$, as claimed.
\hfill$\square$

This finishes the proof of Theorem \ref{genthmE}. \hfill$\blacksquare$

\section{Proof of Theorem \ref{thmF}} In this section we prove Theorem \ref{thmF}, in the following more precise form.

\begin{theorem}\label{genthmF}
 Let $\bar{G}$, $\bar{H}$ be connected real algebraic groups and  $\bar{K}<\bar{G}$, $\bar{L}<\bar{H}$  be $\mathbb R$-subgroups such that $\cap_{g\in \bar{G}}g\bar{K}g^{-1}=Z(\bar{G})$ and $\cap_{h\in\bar{H}}h\bar{L}h^{-1}=Z(\bar{H})$. 
 Denote $G=\bar{G}/Z(\bar{G})$, $K=\bar{K}/Z(\bar{G})$, $H=\bar{H}/Z(\bar{H})$ and $L=\bar{L}/Z(\bar{H})$.  Suppose that $K$ and $L$ are connected.
 
 Let $\Gamma<G$, $\Lambda<H$ be countable dense subgroups such that the translation actions $\Gamma\curvearrowright (G,m_G)$ and $\Lambda\curvearrowright (H,m_H)$ are strongly ergodic.  

Let $A\subset G/K$ and $B\subset H/L$ be non-negligible measurable sets and $\theta:A\rightarrow B$ be a nonsingular isomorphism such that  $\theta(\Gamma x\cap A)=\Lambda\theta(x)\cap B$, for almost every $x\in A$. 

Then we can find a Borel map $\phi:G/K\rightarrow\Lambda$, a topological isomorphism $\delta:G\rightarrow H$ and $h\in H$ such that $\delta(\Gamma)=\Lambda$, $\delta(K)=hLh^{-1}$ and $\theta(x)=\phi(x)\delta(x)hL$, for almost every $x\in A$.
\end{theorem}

{\it Proof.} 
We begin by recording a fact which we will use repeatedly. If $g\in\bar{G}\setminus Z(\bar{G})$, then Lemma \ref{norm} gives that the set of $x\in\bar{G}$ such that $gx\bar{K}=x\bar{K}$ has measure zero. This implies the following:

{\bf Fact}. If $g\in G\setminus\{e\}$, then  the set of $x\in G$ such that $gxK=xK$ has measure zero. 
Similarly, if $h\in H\setminus\{e\}$, then the set of $y\in H$ such that $hyL=yL$ has measure zero.

As the action $\Gamma\curvearrowright G/K$ is ergodic, we may extend $\theta$ to a countable-to-1 map $\theta:G/K\rightarrow H/L$ such that $\theta(\Gamma x)\subset\Lambda\theta(x)$, for almost every $x\in G/K$. Since the action $\Lambda\curvearrowright H/L$ ergodic, we may assume that $B\subset H/L$ is in fact an open set.

The above fact implies that the action $\Lambda\curvearrowright H/L$ is free and therefore we can define a cocycle $w:\Gamma\times G/K\rightarrow\Lambda$ by the formula ${\theta}(gx)=w(g,x){\theta}(x)$. 
Let $\pi:G\rightarrow G/K$ denote the quotient. We define $\Theta:G\rightarrow H/L$ and $W:\Gamma\times G\rightarrow\Lambda$ by $\Theta(x)=\theta(\pi(x))$ and $W(g,x)=w(g,\pi(x))$.  Note that $\Theta(gx)=W(g,x)\Theta(x)$, for all $g\in\Gamma$ and almost every $x\in G$.

\vskip 0.05in
 {\bf The first part of the proof}. In the first part of the proof we show that there exists a Borel map $\alpha:G\rightarrow H$ such that $W(g,x)=\alpha(gx)\alpha(x)^{-1}$, for all $g\in\Gamma$ and almost every $x\in G$.
To this end, we follow closely the proof of \cite[Theorem 4.4]{Io13}. 

\vskip 0.05in
Let $\mu$ be a Borel probability measure on $G$ which is equivalent to $m_G$.

 {\bf Claim 1.} There exists a sequence of Borel maps $\alpha_n:G\rightarrow H$ such that for all $g\in\Gamma$ we have that $\mu(\{x\in G|\;W(g,x)=\alpha_n(gx)\alpha_n(x)^{-1}\})\rightarrow 1$, as $n\rightarrow\infty$.
 
 {\it Proof of Claim 1.}  Let $\varepsilon>0$ and $F\subset\Gamma$ be finite. 
 Since $\bar{L}<\bar{H}$ are real algebraic groups and $\cap_{h\in\bar{H}}h\bar{L}h^{-1}=Z(\bar{H})$,  Lemma \ref{norm}  provides $m\geqslant 1$, an $\bar{H}$-invariant open conull  set $\Omega\subset (\bar{H}/\bar{L})^m$  and a Borel map $\tau:\Omega\rightarrow H=\bar{H}/Z(\bar{H})$ such that $\tau(hx)=h\tau(x)$, for all $h\in\bar{H}$ and $x\in\Omega$.
 Thus, if we identify $\bar{H}/\bar{L}$ with $H/L$, then we can view $\Omega$ as an $H$-invariant open conull subset of $(H/L)^m$ which admits an $H$-equivariant Borel map $\tau:\Omega\rightarrow H$.
 
Since $W$ takes values into a countable group, there exists a neighborhood $V_0$ of $e\in G$ such that \begin{equation}\label{ecuatie}\mu(\{x\in G|\;W(g,xt)=W(g,x)\})\geqslant 1-m^{-1}\varepsilon,\;\;\text{for all}\;\; g\in F\;\;\text{and every}\;\;t\in V_0. \end{equation}
 
Now, since $\theta:G/K\rightarrow H/L$ is a nonsingular map and  $\Omega\subset (H/L)^m$ is conull, we derive that the set $\{(x_1,..,x_m)\in (G/K)^m|\;(\theta(x_1),...,\theta(x_m))\in\Omega\}$ is conull in $(G/K)^m$. Equivalently, we get that the set $\{(x_1,...,x_m)\in G^m|\;(\Theta(x_1),...,\Theta(x_m))\in\Omega\}$ is conull in $G^m$. Fubini's theorem implies that for almost every $(t_1,...,t_m)\in G^m$, the set $\{x\in G|\; (\Theta(xt_1),...,\Theta(xt_m))\in\Omega\}$ is conull in $G$.

In particular, we can find $t_1,...,t_m\in V_0$ such that the set $\{x\in G|\; (\Theta(xt_1),...,\Theta(xt_m))\in\Omega\}$ is conull in $G$. We define $\psi:G\rightarrow (H/L)^m$ by letting $\psi(x)=(\Theta(xt_1),...,\Theta(xt_m))$.
Further, we define $C$ to be the set of $x\in G$ satisfying the following three conditions:

\begin{itemize}
\item  $\psi(x)\in\Omega$ and $\psi(gx)\in\Omega$, for all $g\in F$,
\item  $\Theta(gxt_i)=w(g,xt_i)\Theta(xt_i)$, for all $g\in F$ and $1\leqslant i\leqslant m$, and
\item $W(g,x)=W(g,xt_i)$, for all $g\in F$ and $1\leqslant i\leqslant m$. 
\end{itemize}

Since $\psi(x)\in \Omega$ and $\Theta(gx)=W(g,x)\Theta(x)$, for all $g\in \Gamma$ and almost every $x\in G$, equation \ref{ecuatie} (applied to $t_1,...,t_m\in V_0$) implies that $\mu(C)\geqslant 1-\varepsilon$. 

Finally, we define $\alpha:G\rightarrow H$ by letting $\alpha(x)=\begin{cases} \tau(\psi(x)),\;\;\text{if}\;\; \psi(x)\in\Omega \\ e,\;\;\text{if}\;\;\psi(x)\not\in\Omega\end{cases}.$ 

Then for every $x\in C$ and $g\in F$ we have that $$\alpha(gx)=\tau(\psi(gx))=\tau(\Theta(gxt_1),...,\Theta(gxt_m))=\tau(W(g,xt_1)\Theta(xt_1),...,W(g,xt_m)\Theta(xt_m))=$$
$$\tau(W(g,x)\Theta(xt_1),...,W(g,x)\Theta(xt_m))=W(g,x)\tau(\Theta(xt_1),...,\Theta(xt_m))=W(g,x)\alpha(x).$$

Thus, $\alpha:G\rightarrow H$ is a Borel map which satisfies $\mu(\{x\in G|\alpha(gx)=W(g,x)\alpha(x)\})\geqslant 1-\varepsilon$, for all $g\in F$. Since $\varepsilon>0$ and the finite subset $F\subset\Gamma$  are arbitrary, the claim follows.
 \hfill$\square$
 
Since the action $\Gamma\curvearrowright G$ is strongly ergodic, Claim 1 and Lemma \ref{untwist} imply that we can find a Borel map $\alpha:G\rightarrow H$ such that $W(g,x)=\alpha(gx)\alpha(x)^{-1}$, for all $g\in\Gamma$ and almost every $x\in G$.

\vskip 0.05in
{\bf The second part of the proof}. In the second part of proof we obtain the conclusion by using a strategy similar to the one employed in the second part of the proof of Theorem \ref{genthmE}.
\vskip 0.05in
We start by adapting part of the proof of Theorem \ref{genthmE} to our present context.
Since $G, H$ are connected real algebraic groups, we can find simply connected l.c.s.c. groups $\widetilde G$, $\widetilde H$  together with continuous onto homomorphisms $p:\widetilde G\rightarrow G$, $q:\widetilde H\rightarrow H$ such that $\ker(p)<\widetilde G$, $\ker(q)<\widetilde H$ are discrete subgroups. Denote $\widetilde\Gamma=p^{-1}(\Gamma)$.

Let $\tilde W:\widetilde\Gamma\times\widetilde G\rightarrow\Lambda$ and $\tilde\alpha:\widetilde G\rightarrow\Lambda$ be given by $\tilde W(g,x)=W(p(g),p(x))$ and $\tilde\alpha(x)=\alpha(p(x))$. Then $\tilde W(g,x)=\tilde\alpha(gx)\tilde\alpha(x)^{-1}$, for all $g\in\widetilde\Gamma$ and almost every $x\in\widetilde G$. Since $\widetilde G$ is simply connected,  Corollary \ref{cocrig2} yields a Borel map $\phi:\widetilde G\rightarrow\Lambda$ and a homomorphism $\rho:\widetilde\Gamma\rightarrow\Lambda$ such that \begin{equation}\label{tildeW}\tilde W(g,x)=\phi(gx)\rho(g)\phi(x)^{-1},\;\;\text{for all $g\in\widetilde\Gamma$ and almost every $x\in\widetilde G$}. \end{equation}

 Define $\beta:\widetilde G\rightarrow H$ by letting $\beta(x)=\phi(x)^{-1}\tilde\alpha(x)$. Then we get that $\beta(gx)=\rho(g)\beta(x)$, for all $g\in\tilde\Gamma$ and almost every $x\in\tilde G$. By Lemma \ref{ext}, $\rho$ extends to a continuous homomorphism $\rho:\widetilde G\rightarrow H$.
Define $\tilde\theta:\widetilde G\rightarrow H/L$ by letting $\tilde\theta(x)=\theta(\pi(p(x)))$. By using Fubini's theorem as in the proof of Claim 3 in the proof of Theorem \ref{genthmE}, we find  $z\in H/L$ such that \begin{equation}\label{theta}\tilde\theta(x)=\phi(x)\rho(x)z,\;\;\text{for almost every $x\in G$.} \end{equation}
We are now ready to prove the following claim, which is the core of the second part of the proof:

{\bf Claim 2}. $\rho$ is onto, $\rho(\widetilde\Gamma)=\Lambda$, and $\ker(\rho)<\widetilde G$ is discrete.

{\it Proof of Claim 2}. 
Since $\rho$ is continuous and $\phi$ takes countably many values,  equation \ref{tildeW} implies that $\lim\limits_{g\in\widetilde\Gamma, g\rightarrow e}\tilde W(g,x)=e$, for almost every $x\in\widetilde G$.  This further implies that $\lim\limits_{g\in\Gamma, g\rightarrow e} w(g,x)=e$, for almost every $x\in G/K$. Recall that $\theta:A\rightarrow B$ satisfies $\theta(\Gamma x\cap A)=\Lambda\theta(x)\cap B$, for almost every $x\in A$.  Thus, if $x\in A$, then since $B$ is open we get that $h\theta(x)\in B$, for all $h\in\Lambda$ which are close enough to $e\in H$. Thus, for all $h\in\Lambda$ close enough to $e\in H$ there is $v(h,x)\in\Lambda$ such that $h\theta(x)=\theta(v(h,x)x)$. By repeating the above, with the roles of the two actions reversed (note that the assumptions are symmetric), we get that $\lim\limits_{h\in\Lambda, h\rightarrow e}v(h,x)=e$, for almost every $x\in A$.

 Let $D\subset \widetilde G$ be a non-negligible compact set such that $D_0:=\pi(p(D))\subset A$ and $\phi$ is constant on $D$. Let $\lambda\in\Lambda$ such that $\phi(y)=\lambda$, for all $y\in D$.  
Note that if $y\in D$, then $\tilde\theta(y)=\theta(\pi(p(y)))\in\theta(A)=B$. Since $B$ is open, it follows that $\lim\limits_{h\in H,h\rightarrow e}m_{\widetilde G}(\{y\in D|\; h\tilde\theta(y)\in B\})=m_{\widetilde G}(D)$.
 
Let $y\in D$ and $h\in\Lambda$ such that $h\tilde\theta(y)\in B$. Denote $x=\pi(p(y))$.
Since $h\theta(x)=h\tilde\theta(y)\in B$,  by the first paragraph there exists $v(h,x)\in\Gamma$ such that $h\theta(x)=\theta(v(h,x)x)$. Let $\omega(h,y)\in\widetilde\Gamma$ such that $p(\omega(h,y))=v(h,x)$. 
Since $\lim\limits_{h\in\Lambda, h\rightarrow e}v(h,x)=e$ and $\ker(p)<\widetilde G$ is discrete, we may choose $\omega(h,y)$ such that $\lim\limits_{h\in\Lambda,h\rightarrow e}\omega(h,y)=e$. Then $h\tilde\theta(y)=h\theta(x)=\theta(v(h,x)x)=\tilde\theta(\omega(h,y)y)$ and \ref{theta} gives that $$h\phi(y)\;\rho(y)z=h\tilde\theta(y)=\tilde\theta(\omega(h,y)y)=\phi(\omega(h,y)y)\rho(\omega(h,y)y)z=\phi(\omega(h,y)y)\;\rho(\omega(h,y))\rho(y)z.$$ Since $h\phi(y), \phi(\omega(h,y)y)\rho(\omega(h,y))\in\Lambda$,  $\Lambda$ is countable and
 acts freely on $H/L$, and $\rho(D)z\subset H/L$ is non-negligible,  we conclude that $h\phi(y)=\phi(\omega(h,y)y)\rho(\omega(h,y))$, whenever $h\tilde\theta(y)\in B$. 
 
 Next, since  $\lim\limits_{h\in H,h\rightarrow e}m_{\widetilde G}(\{y\in D|\; h\tilde\theta(y)\in B\})=m_{\widetilde G}(D)$ and $\lim\limits_{h\in\Lambda,h\rightarrow e}\omega(h,y)=e$, we derive that
  $$\lim\limits_{h\in\Lambda,h\rightarrow e}m_{\widetilde G}(\{y\in D|\;h\tilde\theta(y)\in B\;\text{and}\;\omega(h,y)y\in D\})=m_{\widetilde G}(D).$$ 
  
  In particular, we can find a neighborhood $V_1$ of $e\in H$ such that for every $h\in\Lambda\cap V_1$, there is $y\in D$ such that $h\tilde\theta(y)\in B$ and $\omega(h,y)y\in D$. Since $\phi(y)=\phi(\omega(h,y))=\lambda$ (as $\phi\equiv\lambda$ on $D$), we get that $h\lambda=\lambda\rho(\omega(h,y))$ and $\omega(h,y)\in Dy^{-1}\subset DD^{-1}$. Hence $h\in\lambda\rho(\widetilde\Gamma\cap DD^{-1})\lambda^{-1}$. Since $h\in\Lambda\cap V_1$ is arbitrary, we deduce that  $\lambda^{-1}(\Lambda\cap V_1)\lambda\subset\rho(\widetilde\Gamma\cap DD^{-1})$.
 
 Since $\Lambda<H$ is dense and $H$ is connected,  $\Lambda\cap V_1$ generates $\Lambda$. We deduce that $\Lambda=\lambda^{-1}\Lambda\lambda\subset\rho(\widetilde\Gamma)$. Since $\rho(\widetilde\Gamma)\subset\Lambda$, we conclude that $\rho(\widetilde\Gamma)=\Lambda$.
  Moreover, since  $\Lambda<H$ is dense, $\rho$ is continuous, and $DD^{-1}$ is compact, we get that $\lambda^{-1}V_1\lambda\subset\rho(DD^{-1})$. Since $H$ is connected, $V_1$ generates $H$, and it follows that $\rho$ is onto. 
  
  Finally, to see that $\ker(\rho)$ is discrete, let $g_n\in\ker(\rho)$ be a sequence such that $\lim\limits_{n\rightarrow\infty}g_n=e$. Since $\phi$ takes countably many values, we have that $\lim\limits_{n\rightarrow\infty}m_{\widetilde G}(\{x\in D|\;\phi(g_nx)=\phi(x)\})=m_{\widetilde G}(D)$. By using equation \ref{theta} we get that $\lim\limits_{n\rightarrow\infty}m_{\widetilde G}(\{x\in D|\;\tilde\theta(g_nx)=\tilde\theta(x)\})=m_{\widetilde G}(D)$. Since $D_0\subset G/H$ is compact and $\lim\limits_{n\rightarrow\infty}p(g_n)=e$, we  also have that $\lim\limits_{n\rightarrow\infty}m_{G/K}(\{x\in D_0|\;p(g_n)x\in D_0)\})=m_{G/K}(D_0)$.  Since $\pi$ and $p$ are nonsingular, we get that $\lim\limits_{n\rightarrow\infty}m_{\widetilde G}(\{x\in D|\; \pi(p(g_nx))=p(g_n)\pi(p(x))\in D_0\})=m_{\widetilde G}(D)$. 
  
  Altogether, the set of $x\in D$ satisfying $\pi(p(g_nx))\in D_0$ and $\tilde\theta(g_nx)=\tilde\theta(x)$ is non-negligible, for large enough $n$. Since $\theta$ is 1-1 on $D_0\subset A$, we get that $p(g_n)\pi(p(x))=\pi(p(g_nx))=\pi(p(x))$, for all such $x\in D$. By using the fact from the beginning of the proof, we conclude that $p(g_n)=e$, for large enough $n$. Since $\ker(p)<\widetilde G$ is discrete, we must have that $g_n=e$, for large enough $n$.
  \hfill$\square$
 
 The rest of the proof is divided between the following three claims.
 
 {\bf Claim 3}.  The map $\phi:\widetilde G\rightarrow\Lambda$ factors through $p:\widetilde G\rightarrow G$ and there exists a topological isomorphism $\delta:G\rightarrow H$ such that $\delta(\Gamma)=\Lambda$ and $\theta(\pi(x))=\phi(x)\delta(x)z$, for almost every $x\in G$. 
   
  {\it Proof of Claim 3.} 
  By Claim 2, $\rho:\widetilde G\rightarrow H$ is an onto continuous homomorphism with discrete kernel,  hence $\widetilde G$ is a covering group of $H$. The uniqueness of covering groups implies that we can find a topological isomorphism $\tau:\widetilde G\rightarrow\widetilde H$ such that $\rho=q\circ\tau$. Since $G$ and  $H$ have trivial center, we get that $\ker(p)=Z(\widetilde G)$ and $\ker(q)=Z(\widetilde H)$ and hence $\ker(\rho)=\tau^{-1}(\ker(p))=\tau^{-1}(Z(\widetilde H))=Z(\widetilde G)$.
Therefore, $\rho$ descends to a topological isomorphism $\delta:G\rightarrow H$. Since $\rho(\widetilde\Gamma)=\Lambda$, we have $\delta(\Gamma)=\Lambda$.
  
Let $\sigma\in\ker(p)$. Then $\rho(\sigma)=e$ and $\tilde\theta(x\sigma)=\tilde\theta(x)$, for all $x\in\widetilde G$. By using \ref{theta} we deduce that $\phi(x\sigma)\rho(x)z=\phi(x)\rho(x)z$, for almost every $x\in\widetilde G$. 
Since $\phi(x),\phi(x\sigma)\in\Lambda$ and the action $\Lambda\curvearrowright H/L$ is free, we get that $\phi(x\sigma)=\phi(x)$, for almost every $x\in\widetilde G$. Therefore, $\phi:\widetilde G\rightarrow\Lambda$ descends to a map $\phi:G\rightarrow\Lambda$. Together with equation \ref{theta}, this proves the claim. 
  \hfill$\square$

 {\bf Claim 4.} The map $\phi:G\rightarrow\Lambda$ factors through $\pi:G\rightarrow G/K$.
 
 {\it Proof of Claim 4.} Since $\phi$ factors through $p:\widetilde G\rightarrow G$ by Claim 3, equation \ref{tildeW} can be rewritten as $W(g,x)=\phi(gx)\delta(g)\phi(g)^{-1}$, for all $g\in\Gamma$ and almost every $x\in G$.
 Thus, if  $k\in K$, then we have that   $\phi(gxk)\delta(g)\phi(xk)^{-1}=W(g,xk)=W(g,x)=\phi(gx)\delta(g)\phi(x)^{-1}$. 
 
 Hence, if we define $\lambda_k(x)=\phi(x)^{-1}\phi(xk)$, then $\lambda_k(gx)=\delta(g)\lambda(x)\delta(g)^{-1}$, for all $g\in\Gamma$ and almost every $x\in G$. This implies that $C_k=\{x\in G|\lambda_k(x)=e\}$ is $\Gamma$-invariant, for every $k\in K$. Since $\Lambda$ is countable we can find a neighborhood $V_2$ of $e\in G$ such that $m_G(C_k)>0$, for all $k\in K\cap V_2$. Since the action $\Gamma\curvearrowright G$ is ergodic, we deduce that 
 $C_k=G$, almost everywhere, for every $k\in K\cap V_2$. Thus, the set $K_0$ of $k\in K$ such that $C_k=G$, almost everywhere, is an open subgroup of $K$. Since $K$ is connected, we conclude that $K_0=K$. This clearly implies the claim.
 \hfill$\square$
 
 Let $h\in H$ such that $z=hL$.
 
 {\bf Claim 5}. $\delta(K)=hLh^{-1}$.

 {\it Proof of Claim 5.} Let $k\in K$. Then by Claim 4, for almost every $x\in G$ we have $\phi(xk)=\phi(x)$, hence $\phi(x)\delta(x)yL=\theta(x)=\theta(xk)=\phi(x)\delta(xk)yL$. Thus,  $\delta(k)yL=yL$ and therefore $\delta(k)\in yLy^{-1}$.
 
 To show the reverse inclusion, let $g_n\in\delta^{-1}(hLh^{-1})$ be a sequence such that $\lim\limits_{n\rightarrow\infty}g_n=e$. We claim that $g_n\in K$, for $n$ large enough. 
Indeed, the set of $x\in \pi^{-1}(A)$ such that $\phi(g_nx)=\phi(x)$ and $g_nx\in\pi^{-1}(A)$ is non-negligible, for large enough $n$. Since $\delta(g_n)yL=yL$,  by Claim 3 for almost every such $x$ we have that $\theta(xg_nK)=\phi(xg_n)\delta(x)\delta(g_n)yL=\phi(x)\delta(x)yL=\theta(xK)$. Since $xK,xg_nK\in A$ and $\theta$ is 1-1 on $A$, we get that $xg_nK=xK$ and hence $g_n\in K$.

 The previous paragraph implies that there exists a neighborhood $V_3$ of $e\in G$ such that we have $\delta^{-1}(hLh^{-1})\cap V_3\subset K$. Since $\delta(K)\subset hLh^{-1}$, we derive that $K\subset\delta^{-1}(hLh^{-1})$ is an open subgroup. Since $L$ is connected, we get that $K=\delta^{-1}(hLh^{-1})$ and thus $\delta(K)=hLh^{-1}$.
 \hfill$\blacksquare$

\section{Proofs of Propositions \ref{propG} and \ref{propH}}

\subsection {Proof of Proposition \ref{propG}} Let $\alpha$ denote the  action $\Gamma\curvearrowright (G_1,m_{G_1})$ given by $gx=p(g)x$.

 (1) 
 In order to show that  $\alpha$ is strongly ergodic, by Lemma \ref{lattice} it suffices to show that the induced action $G\curvearrowright^{\tilde\alpha} (G/\Gamma\times G_1,m_{G/\Gamma}\times m_{G_1})$ is strongly ergodic.
Note that $\alpha$ is the restriction to $\Gamma$ of the action $G\curvearrowright (G_1,m_{G_1})$ given by $gx=p(g)x$, for all $g\in G,x\in G_1$. It is well-known (see e.g. \cite[Proposition 4.2.22]{Zi84}) that $\tilde\alpha$ is isomorphic to the product action $G\curvearrowright^{\beta} (G/\Gamma\times G_1,m_{G/\Gamma}\times m_{G_1})$ given by $g(x,y)=(gx,p(g)y)$, for all $g\in G,x\in G/\Gamma$ and $y\in G_1$.

To show that $\beta$ (and hence $\tilde\alpha$) is strongly ergodic,  fix a Borel probability measure $\mu$ on $G_1$ which is equivalent to $m_{G_1}$. Then $\tilde\mu=m_{G/\Gamma}\times\mu$ is a  Borel probability measure on $G/\Gamma\times G_1$ which is equivalent to $m_{G/\Gamma}\times m_{G_1}$. Let $\{A_n\}$ be a sequence of measurable subsets of $G/\Gamma\times G_1$ satisfying \begin{equation}\label{ai}\lim\limits_{n\rightarrow\infty}\sup_{g\in K}\tilde\mu(gA_n\Delta A_n)=0,\;\;\text{for every compact set $K\subset G$}.\end{equation}

Next, recall that the representation $\rho:G_2\rightarrow\mathcal U(L^2(G/\Gamma,m_{G/\Gamma})\ominus\mathbb C1)$ has spectral gap.
Note that  the restriction of $\beta$ to $G_2$ preserves $\tilde\mu$. Let $\pi:G_2\rightarrow\mathcal U(L^2(G/\Gamma\times G_1,\tilde\mu))$ be the associated  Koopman representation. Then $\pi(G_2)$ acts trivially on the subspace $L^2(G_1,\mu)\subset L^2(G/\Gamma\times G_1,\tilde\mu)$. Moreover, the restriction of $\pi$ to $L^2(G/\Gamma\times G_1,\tilde\mu)\ominus L^2(G_1,\mu)$ is unitarily equivalent to $\oplus_{i=1}^{\infty}\rho$ and therefore has spectral gap.
We denote by $P:L^2(G/\Gamma\times G_1,\tilde\mu)\rightarrow L^2(G_1,\mu)$ the orthogonal projection.

Finally,  equation \ref{ai} gives that we have $\sup_{g\in K}\|\pi(g)(1_{A_n})-1_{A_n}\|_{L^2(\tilde\mu)}\rightarrow 0$, for every compact set $K\subset G_2$. By using the spectral gap property described in the previous paragraph we deduce that $\|1_{A_n}-P(1_{A_n})\|_{L^2(\tilde\mu)}\rightarrow 0$. This  easily implies that there exists a sequence $\{B_n\}$ of measurable subsets of $G_1$ such that $\tilde\mu(A_n\Delta (G/\Gamma\times B_n))\rightarrow 0$. In combination with \ref{ai} this further implies that $\sup_{g\in K}\mu(gB_n\Delta B_n)\rightarrow 0$, for every compact set $K\subset G_1$.  Since the action $G_1\curvearrowright (G_1,m_{G_1})$ is strongly ergodic by Lemma \ref{strlc}, we get that $\lim\limits_{n\rightarrow\infty}\tilde\mu(A_n)(1-\tilde\mu(A_n))=\lim\limits_{n\rightarrow\infty}\mu(B_n)(1-\mu(B_n))=0$.

(2) Since $G_2$ has property (T), by \cite[Proposition 2.4]{Zi81} the action $G_2\curvearrowright (G/\Gamma,m_{G/\Gamma})$ has property (T). Applying \cite[Proposition 3.5]{PV08} gives that the action $\Gamma\curvearrowright (G/G_2,m_{G/G_2})$ has property (T).    Since  the action $\Gamma\curvearrowright (G/G_2,m_{G/G_2})$ is isomorphic to $\alpha$, we conclude that $\alpha$ has property (T).
\hfill$\blacksquare$

\subsection{Proof of Proposition \ref{propH}}
Let  $\nu$ be a Borel probability measure of $X:=K\backslash G$ which is quasi-invariant under the right $G$-action. Let $\pi:X\rightarrow G$ be a Borel map satisfying $K\pi(x)=x$, for all $x\in X$. The map $\theta:K\times X\rightarrow G$ defined by $\theta(k,x)=k\pi(x)$ is a $K$-equivariant Borel isomorphism (where we consider the action $K\curvearrowright K\times X$ given by $k(k'x)=(kk',x)$). We identify $G$ with $K\times X$ via $\theta$. Since the push-forward of $m_K\times\nu$ through $\theta$ is equivalent to $m_G$, we may view $\mu=m_K\times\nu$ as a probability measure on $G$  which is equivalent to $m_G$.

Assume that the action $\Gamma\cap K\curvearrowright (K,m_{K})$ has spectral gap. Our goal is to show that the action $\Gamma\curvearrowright (G,m_G)$ is strongly ergodic.
Let $\{A_n\}$ be a sequence of measurable subsets of $G$ such that $\lim\limits_{n\rightarrow\infty}\mu(gA_n\Delta A_n)=0$, for all $g\in\Gamma$. For all $n$ and $x\in X$, let $A_n^x=\{k\in K|(k,x)\in A_n\}$. 

Then  $\mu(gA_n\cap A_n)=\displaystyle{\int_{X}m_K(gA_n^x\Delta A_n^x)\;\text{d}\nu(x)}$, for all $n$ and $g\in K$.
Thus, after replacing $\{A_n\}$ with a subsequence, we have that $\lim\limits_{n\rightarrow\infty}m_K(gA_n^x\Delta A_n^x)=0$, for all $g\in\Gamma\cap 
K$ and almost every $x\in X$.

Since the  action $\Gamma\cap K\curvearrowright (K,m_{K})$ has spectral gap, we get that $\lim\limits_{n\rightarrow\infty}m_K(A_n^x)(1-m_K(A_n^x))=0$, for almost every $x\in X$.
Further, this implies that \begin{equation}\label{Kinv}\lim\limits_{n\rightarrow\infty}\sup_{g\in K}\mu(gA_n\Delta A_n)=0. \end{equation}

Next, we need the following result whose proof we postpone until the end of this section.

\begin{lemma}\label{SO}
Let $G=SL_n(\mathbb R)$ and $K=SO_n(\mathbb R)$ for some $n\geqslant 2$. Let $g\in G\setminus K$.

Then the group generated by $g$ and $K$ is equal to $G$. Moreover, if $Y\subset G$ is a compact set, then we can find $m\geqslant 1$ and $g_1,...,g_m\in\{g,g^{-1}\}$ such that $Y\subset g_1Kg_2...g_{m-1}Kg_m$.
\end{lemma}

\begin{remark}\label{alireza} Lemma \ref{SO} implies that if $g\in G\setminus K$, then there exists $N\geqslant 1$ such that the compact groups $\{g^nKg^{-n}\}_{|n|\leqslant N}$ generate $G$.
I am grateful to  Alireza Salehi-Golsefidy for pointing out to me that the following general result holds: if $G$ is a semisimple algebraic group of real rank $r$ over a local field $k$, then $G$ is generated by $r+1$ compact subgroups.
\end{remark}

Going back to the proof of Proposition \ref{propH}, 
let $g\in\Gamma\setminus K$ and $Y\subset G$ be a compact set. By Lemma \ref{SO} we can find  $g_1,...,g_m\in\{g,g^{-1}\}$ such that $Y\subset g_1Kg_2...g_{m-1}Kg_m$. Since $g\in\Gamma$ we have that $\lim\limits_{n\rightarrow\infty}\mu(gA_n\Delta A_n)=\lim\limits_{n\rightarrow\infty}\mu(g^{-1}A_n\Delta A_n)=0$.
In combination with equation \ref{Kinv}, these facts  imply that $\lim\limits_{n\rightarrow\infty}\sup_{g\in Y}\mu(gA_n\Delta A_n)=0$.
 Since $Y\subset G$ is an arbitrary compact set and the  action $G\curvearrowright (G,m_G)$ is strongly ergodic by Lemma \ref{strlc}, we get that $\lim\limits_{n\rightarrow\infty}\mu(A_n)(1-\mu(A_n))=0$. This shows that the sequence $\{A_n\}$ is trivial. \hfill$\blacksquare$

\subsection{Proof of Lemma \ref{SO}} We denote by diag$(x_1,....,x_n)$ the diagonal $n\times n$ matrix whose diagonal entries are $x_1,...,x_n$.
We also denote by $Tr:\mathbb M_n(\mathbb C)\rightarrow\mathbb C$ the usual trace.

Let $H$ be the group generated by $g$ and $K$. 
Let $A<G$ be the subgroup of positive diagonal matrices. Write $g=k_1dk_2$, where $k_1,k_2\in K$ and $d\in A$. Write $d=\text{diag}(\lambda_1,...,\lambda_n)$.  Since $g\notin K$ we have that $d\not=I$, hence not all the $\lambda_i$'s are equal.  We assume for simplicity that $\lambda_1\not=\lambda_2$.

Let us first show that $H=G$ in the case $n=2$. In this case, $d=\text{diag}(\lambda_1,\lambda_2)$. For every $0\leqslant\alpha\leqslant 1$, we define the $2\times 2$ matrix $$k_{\alpha}:=\begin{pmatrix}\alpha&\sqrt{1-\alpha^2}\\-\sqrt{1-\alpha^2}&\alpha\end{pmatrix}\in K.$$ A direct computation shows that  $Tr(dk_{\alpha}dk_{\alpha}^*d)=\alpha^2(\lambda_1^3+\lambda_2^3-\lambda_1-\lambda_2)+(\lambda_1+\lambda_2)$. Let $a\in A$ with $\lambda_1+\lambda_2\leqslant Tr(a)\leqslant\lambda_1^3+\lambda_2^3$. Then we can find $\alpha\in [0,1]$ such that $Tr(a)=Tr(dk_{\alpha}dk_{\alpha}^*d)$. Since $a$ and $dk_{\alpha}dk_{\alpha}^*d$ are positive $2\times 2$ matrices,  it follows that there exists $k\in K$ such that $a=k(dk_{\alpha}dk_{\alpha}^*d)k^*$. Hence $a\in H$. This shows that $H\cap A$ contains $\{a\in A|\lambda_1+\lambda_2\leqslant Tr(a)\leqslant\lambda_1^3+\lambda_2^3\}$. Since $\lambda_1\not=\lambda_2$, we derive that $H\cap A$ is an open subgroup of $A$. Since $A$ is connected we deduce that $H\cap A=A$, therefore $A\subset H$. Since $G=KAK$ and $K\subset H$, we conclude that $H=G$.

Before proceeding to the general case, let us derive an additional fact in the case $n=2$. Note that we can find $a_1,...,a_4\in A$ such that $I=a_1a_2a_3a_4$ and $\lambda_1+\lambda_2\leqslant Tr(a_i)\leqslant\lambda_1^3+\lambda_2^3$, for all $1\leqslant i\leqslant 4$.  The previous paragraph implies that we can find $k_0,k_1,...,k_{12}\in K$ such that $I=k_0dk_1...k_{11}dk_{12}$.
 
 Now, assume that $n\geqslant 2$ is arbitrary.  For $1\leqslant i<j\leqslant n$, we denote by $G_{i,j}$ the subgroup of $G$ consisting of all matrices $g=(g_{k,l})\in G$ with the property that $g_{k,l}=\delta_{k,l}$, for all $1\leqslant k,l\leqslant n$ with $\{k,l\}\not\subset\{i,j\}$. We clearly have that $G_{i,j}\cong SL_2(\mathbb R)$. We claim that $G_{1,2}\subset H$.

 Let $g=\displaystyle{\text{diag}(\sqrt{\frac{\lambda_1}{\lambda_2}},\sqrt{\frac{\lambda_2}{\lambda_1}})}\in SL_2(\mathbb R)$ and view $g\in G_{1,2}$.
 By applying the above fact to $g$, we can find $k_0,k_1,...,k_{12}\in G_{1,2}\cap K$ such that $I=k_0gk_1...k_{11}gk_{12}$.
 Let $h=dg^{-1}=\text{diag}(\sqrt{\lambda_1\lambda_2},\sqrt{\lambda_1\lambda_2},\lambda_3,...,\lambda_n)$. Then $d=gh$ and $h$ commutes with $G_{1,2}$. 
 It is then immediate that $$k_0dk_1...k_{11}dk_{12}=(k_0gk_1...k_{11}gk_{12})h^{12}=h^{12}.$$ Since $d\in H$, we get that $h^{12}\in H$ and further that $g^{12}=d^{12}h^{-12}\in H$. 
 Thus, $g^{12}\in G_{1,2}\cap H$. Since $g^{12}\not\in K$ and $G_{1,2}\cap K\subset H$, by using the case $n=2$, we get that $G_{1,2}\subset H$.

Let $1\leqslant i<j\leqslant n$. Then there exists $k\in K$ such that $G_{i,j}=kG_{1,2}k^{-1}$ and hence $G_{i,j}\subset H$. Since the groups $G_{i,j}$ generate $G$, we conclude that $G=H$.

To see the moreover assertion, note that since $g\not\in K$,  the above gives that $g$ and $K$ generate $G$. Consequently, we can find 
$p\geqslant 1$ and $h_1,...,h_p\in\{g,g^{-1}\}$ such that $A=h_1Kh_2...h_{p-1}Kh_p$ is non-negligible. It follows that $A^{-1}A$ contains an open neighborhood $V$ of the identity in $G$. Since $G$ is connected we get that $G=\cup_{q\geqslant 1}V^q$. Further, since $Y$ is compact, we can find $q\geqslant 1$ such that $Y\subset V^q$. Therefore, we have that $Y\subset (A^{-1}A)^q$, which proves the claim. 
\hfill$\blacksquare$

\section{Proofs of Corollaries \ref{corI}, \ref{corJ} and \ref{corK}} 

\subsection{Proof of Corollary \ref{corI}} Assume that $\Sigma<SL_m(\mathbb R)$ and $\Delta<SL_n(\mathbb R)$ are both either 
\begin{enumerate}
\item discrete subgroups, or
\item connected real algebraic subgroups. 
\end{enumerate}

Suppose that the actions $SL_m(\mathbb Z[S^{-1}])\curvearrowright SL_m(\mathbb R)/\Sigma$ and $SL_n(\mathbb Z[T^{-1}])\curvearrowright SL_n(\mathbb R)/\Delta$ are SOE. 

Denote by $Z_m$ the center of $SL_m(\mathbb R)$ and by $PSL_m(\mathbb R)$ the quotient  group $SL_m(\mathbb R)/Z_m$. 
Following Example \ref{ex1}, the action $SL_m(\mathbb Z[S^{-1}])\curvearrowright SL_m(\mathbb R)/\Sigma$ is SOE to $PSL_m(\mathbb Z[S^{-1}])\curvearrowright PSL_m(\mathbb R)/\Sigma_0$, where $\Sigma_0=(\Sigma Z_m)/Z_m$. 
Moreover, Example \ref{ex3} gives that the action $SL_m(\mathbb Z[S^{-1}])\curvearrowright SL_m(\mathbb R)$ and hence the action $PSL_m(\mathbb Z[S^{-1}])\curvearrowright PSL_m(\mathbb R)$ is strongly ergodic.
Similarly, we deduce that $SL_n(\mathbb Z[T^{-1}])\curvearrowright SL_n(\mathbb R)/\Delta$ is SOE to $PSL_n(\mathbb Z[T^{-1}])\curvearrowright PSL_n(\mathbb R)/\Delta_0$, where $\Delta_0=(\Delta Z_n)/Z_n$, and that the action $PSL_n(\mathbb Z[T^{-1}])\curvearrowright PSL_n(\mathbb R)$ is strongly ergodic.

Now, in case (1), we have that $\Sigma_0<PSL_m(\mathbb R)$ and $\Delta_0<PSL_n(\mathbb R)$ are discrete subgroups.
In case (2), we have that $\Sigma Z_m<SL_m(\mathbb R)$ and $\Delta Z_n<SL_n(\mathbb R)$ are real algebraic subgroups. Moreover, since $\Sigma$ and $\Delta$ are connected, we get that $\Sigma_0\cong\Sigma/(\Sigma\cap Z_m)$ and $\Delta_0\cong\Delta/(\Delta\cap Z_n)$ are connected.

Altogether, we can apply Theorem \ref{genthmE} in case (1) and Theorem \ref{genthmF} in case (2) to conclude that there exists an isomorphism $\delta:PSL_m(\mathbb R)\rightarrow PSL_n(\mathbb R)$ such that $\delta(PSL_m(\mathbb Z[S^{-1}]))=PSL_n(\mathbb Z[T^{-1}])$. This implies that $m=n$ and there exists $g\in GL_m(\mathbb R)$ such that either $\delta(x)=gxg^{-1}$, for all $x\in PSL_m(\mathbb R)$, or $\delta(x)=g(x^{t})^{-1}g^{-1}$, for all $x\in PSL_m(\mathbb R)$. It now follows easily that $S=T$.
\hfill$\blacksquare$

\subsection{Proof of Corollary \ref{corJ}} In this subsection, we establish the following more precise version of Corollary \ref{corJ}.

\begin{theorem}
Let $G=SL_n(\mathbb R)$, for some $n\geqslant 2$. Let $\Gamma<G$ be a countable dense subgroup which contains the center of $G$ such that the translation action $\Gamma\curvearrowright (G,m_G)$ is strongly ergodic.

 If $\theta:G\rightarrow G$ is a nonsingular isomorphism satisfying $\theta(\Gamma x)=\Gamma\theta(x)$, for almost every $x\in G$,
then we can find a Borel map $\phi:G\rightarrow\Gamma$, $g\in GL_n(\mathbb R)$ and $h\in G$ such that either $\theta(x)=\phi(x)gxg^{-1}h$  or $\theta(x)=\phi(x)g(x^{-1})^{t}g^{-1}h$, for almost every $x\in G$.
 In particular,  $\theta$ preserves $m_G$.

\end{theorem}

{\it Proof.}  Let $\theta:G\rightarrow G$ be a nonsingular isomorphism with $\theta(\Gamma x)=\Gamma\theta(x)$, for almost every $x\in G$. 
 Let $Z$ denote the center of $G=SL_n(\mathbb R)$. Put $G_0=G/Z=PSL_n(\mathbb R)$ and $\Gamma_0=\Gamma/Z$. By Example \ref{ex1}, the actions $\Gamma\curvearrowright G$ and $\Gamma_0\curvearrowright G_0$ are orbit equivalent. Moreover, there exists an orbit equivalence $\tau:G\rightarrow G_0$  such that $\tau(x)\in\Gamma_0 xZ$, for almost every $x\in G$. It follows that $\tau^{-1}(xZ)\in\Gamma x$, for almost every $x\in G$.

 Let $\theta_0=\tau\theta\tau^{-1}:G_0\rightarrow G_0$. Then $\theta_0$ is a nonsingular isomorphism and $\theta_0(\Gamma_0 x)=\Gamma_0\theta_0(x)$, for almost every $x\in G_0$. 
 Since  $G_0$ has trivial center and admits a l.c.s.c. universal covering group, and the action $\Gamma_0\curvearrowright G_0$ is strongly ergodic, by Theorem \ref{OE} we can find an automorphism $\delta:G_0\rightarrow G_0$ satisfying $\delta(\Gamma_0)=\Gamma_0$, a Borel map $\phi:G_0\rightarrow\Gamma_0$, and $h\in G_0$ such that $\theta_0(x)=\phi(x)\delta(x)h$, for almost every $x\in G_0$. 
 Thus, there is $g\in GL_n(\mathbb R)$ such that either $\delta(x)=gxg^{-1}$, for all $x\in G_0$, or $\delta(x)=g(x^{t})^{-1}g^{-1}$, for all $x\in G_0$.  In particular, $\delta$ lifts to an automorphism to $G$ (defined by the same formulas). Since $\delta(\Gamma_0)=\Gamma_0$, it follows that  $\delta(\Gamma)=\Gamma$. Using the properties of $\tau$ listed in the previous paragraph, it is easy to check that $\theta$ has the desired form. 
 
 Moreover, since the conjugation action of $GL_n(\mathbb R)$ on $G$ preserves the Haar measure $m_G$, we deduce that $\theta$ is measure preserving.
 \hfill$\blacksquare$

\subsection{Proof of Corollary \ref{corK}} 
In this subsection, we prove an OE superrigidity result for actions $PSL_m(\mathbb Z[S^{-1}])\curvearrowright PSL_m(\mathbb R)/\Sigma$, where $\Sigma<PSL_m(\mathbb R)$ is an arbitrary discrete subgroup and $m\geqslant 3$, and then explain how Corollary \ref{corK} follows from this result.

\begin{theorem}\label{gencorK}
Let $m\geqslant 3$ be an integer, $G'=PSL_m(\mathbb R)$, and $\Sigma'<G'$ be a discrete subgroup. Let $S$ be a nonempty set of primes and denote $\Gamma'=PSL_m(\mathbb Z[S^{-1}])$. 

Then a free ergodic nonsingular action $\Lambda\curvearrowright (Y,\nu)$ of a countable group $\Lambda$ is SOE to the left translation action $\Gamma'\curvearrowright G'/\Sigma'$ if and only we can find a subgroup $\Lambda_0<\Lambda$, a finite normal subgroup $N<\Lambda_0$, and a normal subgroup $M'<\Sigma'$ such that 
\begin{itemize}
\item $\Lambda\curvearrowright Y$ is induced from some nonsingular action $\Lambda_0\curvearrowright Y_0$, and
 \item $\Lambda_0/N\curvearrowright Y_0/N$ is conjugate to the left-right multiplication action $\Gamma'\times\Sigma'/M'\curvearrowright G'/M'$ given by $(g,\sigma M')\cdot xM'=gx\sigma^{-1}M'$, for all $g\in\Gamma',\sigma\in\Sigma'$ and $x\in G'$.
\end{itemize}
\end{theorem}

{\it Proof.}
 Let $\Gamma=SL_m(\mathbb Z[S^{-1}])$ and $G=SL_m(\mathbb R)$. We denote by $\widetilde G$ the common universal cover of  $G$ and $G'$, and by $\pi:\widetilde G\rightarrow G$, $\pi':\widetilde G\rightarrow G'$ the covering homomorphisms. Let $\widetilde\Gamma=\pi^{-1}(\Gamma)=\pi'^{-1}(\Gamma')$.

{\bf Claim}. There is a subgroup $\widetilde\Gamma_1<\widetilde\Gamma$ such that $g\widetilde\Gamma_1 g^{-1}\cap\widetilde\Gamma_1$ is dense in $\widetilde G$, for every $g\in\widetilde\Gamma$, and the action $\widetilde\Gamma_1\curvearrowright\widetilde G$ has property (T).

{\it Proof of the claim.}
Fix $p\in S$ and denote $\Gamma_1=SL_m(\mathbb Z[\displaystyle{\frac{1}{p}}])$.  Since $m\geqslant 3$, by Example \ref{ex3}, the translation action $\Gamma_1\curvearrowright G$  has property (T).

Moreover, if $g\in\Gamma$, then $g\Gamma_1g^{-1}\cap\Gamma_1$ is dense in $G$. To see this, note first that we can find an integer $N\geqslant 1$ such that $p\nmid N$  and $Ng,Ng^{-1}\in GL_m(\mathbb Z[\displaystyle{\frac{1}{p}}])$. Consider the quotient homomorphism $\rho:\Gamma_1\rightarrow SL_m(\mathbb Z[\displaystyle{\frac{1}{p}}]/N^2\mathbb Z[\displaystyle{\frac{1}{p}}])$. Then $\Gamma_2:=\ker(\rho)$ satisfies $g^{-1}\Gamma_2g\subset\Gamma_1$, hence $\Gamma_2\subset g\Gamma_1g^{-1}\cap\Gamma_1$. Since $G$ is connected, $\Gamma_1<G$ is dense, and $\Gamma_2$ is a finite index subgroup of $\Gamma_1$, we conclude that $\Gamma_2$ is dense in $G$. This proves our assertion.

Let $\widetilde\Gamma_1=\pi^{-1}(\Gamma_1)$. By Example \ref{ex1} we have that $\widetilde\Gamma_1\curvearrowright\widetilde G$ is SOE to $\Gamma_1\curvearrowright G$ and therefore has property (T). Since $\widetilde{G}$ is connected and $\pi$ is a finite-to-$1$ map, by using the above assertion, it follows that $g\widetilde\Gamma_1 g^{-1}\cap\widetilde\Gamma_1$ is dense in $\widetilde G$, for every $g\in\widetilde\Gamma$. \hfill$\square$

We are now ready to prove the {\it only if} assertion of Theorem \ref{gencorK}. The {\it if} assertion follows easily (e.g. by using Example \ref{ex1}) and is left to the reader.  
Let $\Lambda\curvearrowright (Y,\nu)$ be a nonsingular action which is SOE to $\Gamma'\curvearrowright G'/\Sigma'$.
 Let $\widetilde\Sigma=\pi'^{-1}(\Sigma)$ and note that we can identify $\widetilde G/\widetilde\Sigma$ with $G'/\Sigma'$ via the map $x\widetilde\Sigma\mapsto\pi'(x)\Sigma'$. Then, under this identification, the actions $\widetilde\Gamma\curvearrowright\widetilde G/\widetilde\Sigma$ and $\Gamma'\curvearrowright G'/\Sigma'$
 have the same orbits. We therefore deduce that $\Lambda\curvearrowright Y$  is SOE to $\widetilde\Gamma\curvearrowright \widetilde G/\widetilde\Sigma$. 
 
 Since $\widetilde G$ is simply connected, the claim allows us to apply Theorem \ref{OESUP1}.
Thus, we can find  a normal subgroup $\Delta<\widetilde\Gamma\times\widetilde\Sigma$, 
a subgroup $\Lambda_0<\Lambda$, and a $\Lambda_0$-invariant measurable subset $Y_0\subset Y$ such that 
\begin{itemize}
\item $\Delta$ is discrete in $\widetilde G\times\widetilde G$,
\item the left-right multiplication action $\Delta\curvearrowright\widetilde G$ admits a measurable fundamental domain,
\item the left-right multiplication action $(\widetilde\Gamma\times\widetilde\Sigma)/\Delta\curvearrowright \widetilde G/\Delta$ is conjugate to $\Lambda_0\curvearrowright Y_0$, and 
\item the action $\Lambda\curvearrowright Y$ is induced from $\Lambda_0\curvearrowright Y_0$.
\end{itemize}

Next, we claim that $\Delta\subset Z(\widetilde G)\times\widetilde\Sigma$. To see this, let $\Delta_1=\{x\in\widetilde\Gamma|(x,e)\in\Delta\}$. Then $\Delta_1$ is discrete in $\widetilde G$ and normal in $\widetilde\Gamma$.
Since $\widetilde\Gamma<\widetilde G$ is dense and $\widetilde G$ is connected, it follows that $\Delta_1\subset Z(\widetilde G)$. Now, let $(x,y)\in\Delta$. Since $\Delta$ is normal in $\widetilde\Gamma\times\widetilde\Sigma$, for every $g\in\widetilde\Gamma$ we have that $(gxg^{-1},y)\in\Delta$, and hence $gxg^{-1}x^{-1}\in\Delta_1$. Since $\Delta_1$ is finite, we derive that $gxg^{-1}=x$, for every $g\in\widetilde\Gamma$ that is sufficiently close to $e$. Using again the fact that $\widetilde\Gamma<\widetilde G$ is dense and $\widetilde G$ is connected, we get that $x\in Z(\widetilde G)$. This shows that $\Delta\subset Z(\widetilde G)\times\widetilde\Sigma$, as claimed.

Since $Z(\widetilde G)$ is finite and $\Delta<\widetilde\Gamma\times\widetilde\Sigma$ is normal, we can find a normal subgroup $\widetilde M<\widetilde\Sigma$ such that $\Delta\subset Z(\widetilde G)\times\widetilde M$ and the inclusion $\Delta\subset Z(\widetilde G)\times\widetilde M$ has finite index. We may clearly assume that $\widetilde M$ contains $Z(\widetilde G)$, since the latter group is finite. Then $M'=\pi'(\widetilde M)$ is a normal subgroup of $\Sigma'$. Moreover, since the kernel of  $\pi':\widetilde G\rightarrow G'$ is equal to $Z(\widetilde G)$, we have that $\widetilde\Gamma/Z(\widetilde G)\cong\Gamma'$. Also, we have that $\widetilde\Sigma/\widetilde M\cong\Sigma'/M'$.

 Let $\delta:(\widetilde\Gamma\times\widetilde\Sigma)/\Delta\rightarrow\Lambda_0$ be the group isomorphism provided by the above conjugacy of actions.
 Then $N:=\delta((Z(\widetilde G)\times\widetilde M)/\Delta)$ is a finite normal subgroup of $\Lambda_0$ and the action $\Lambda_0/N\curvearrowright Y_0/N$ is conjugate to the left-right multiplication action of $[(\widetilde\Gamma\times\widetilde\Sigma)/\Delta]/[(Z(\widetilde G)\times\widetilde M)/\Delta]\cong\Gamma'\times(\Sigma'/M')$ on $(\widetilde G/\Delta)/[(Z(\widetilde G)\times\widetilde M)/\Delta]$. 
 Since the latter space can be identified with $G'/M'$, we are done.
 \hfill$\blacksquare$
 
 {\bf Proof of Corollary \ref{corK}}. Let us briefly indicate how Theorem \ref{gencorK} implies Corollary \ref{corK}. By applying Theorem \ref{gencorK} in the case $\Sigma'=\{e\}$, the first part  of Corollary \ref{corK} follows.  For the second part of Corollary \ref{corK}, suppose (in the notation from Theorem \ref{gencorK}) that $\Sigma'<\Gamma'$ is a lattice. Assume that a free ergodic nonsingular action $\Lambda\curvearrowright Y$ is SOE to $\Gamma'\curvearrowright G'/\Sigma'$. Then there exist a subgroup $\Lambda_0<\Lambda$ and normal subgroups $N<\Lambda_0$, $M'<\Sigma'$ such that the conclusion of Theorem \ref{gencorK} holds true.

Since $\Sigma'<G'$ is a lattice, $G'$ has trivial center and $\mathbb R$-rank$(G')=m-1\geqslant 2$, Margulis' normal subgroup theorem (see e.g. \cite[Theorem 8.1.2]{Zi84}) implies that either $M'=\{e\}$ or $\Sigma'/M'$ is finite. If $M'=\{e\}$, then we get that $\Lambda_0/N\curvearrowright Y_0/N$ is conjugate to the left-right multiplication action $\Gamma'\times\Sigma'\curvearrowright G'$. If $\Sigma'/M'$ is finite, let $\delta: \Gamma'\times\Sigma'/M'\rightarrow\Lambda_0/N$ be the group homomorphism that witnesses the conjugacy between $\Gamma'\times\Sigma'/M'\curvearrowright G'/M'$ and 
$\Lambda_0/N\curvearrowright Y_0/N$ given by Theorem \ref{gencorK}. Let $N'<\Lambda_0$ be a finite normal subgroup which contains $N$ and satisfies $\delta(\{e\}\times\Sigma'/M')=N'/N$. It is now easy to see that the actions $\Gamma'\curvearrowright G'/\Sigma'$ and $\Lambda_0/N'\curvearrowright Y_0/N'$ are conjugate. \hfill$\blacksquare$

\end{document}